\documentclass[11pt]{article}
\usepackage{amssymb}
\usepackage{amsmath}
\usepackage{times}
\usepackage[all]{xy}

\title{Acyclicity and reduction\indent}
\author{Dominique LECOMTE}
\date{\today}

\def\ufootnote#1{\let\savedthfn\thefootnote\let\thefootnote\relax
\footnote{#1}\let\thefootnote\savedthfn\addtocounter{footnote}{-1}}

\newcommand{\Ana}{{\it\Sigma}^{1}_{1}}

\newcommand{\Borel}{{\it\Delta}^{1}_{1}}
\newcommand{\ana}{{\bf\Sigma}^{1}_{1}}

\newcommand{\boraone}{{\bf\Sigma}^{0}_{1}}
\newcommand{\boratwo}{{\bf\Sigma}^{0}_{2}}
\newcommand{\borathree}{{\bf\Sigma}^{0}_{3}}

\newcommand{\boraxi}{{\bf\Sigma}^{0}_{\xi}}

\newcommand{\borone}{{\bf\Delta}^{0}_{1}}
\newcommand{\bortwo}{{\bf\Delta}^{0}_{2}}

\newcommand{\bormone}{{\bf\Pi}^{0}_{1}}

\newcommand{\bormtwo}{{\bf\Pi}^{0}_{2}}
\newcommand{\bormthree}{{\bf\Pi}^{0}_{3}}

\newcommand{\bormxp}{{\bf\Pi}^{0}_{\xi +1}}

\newcommand{\boraxp}{{\bf\Sigma}^{0}_{\xi +1}}

\newcommand{\bormxi}{{\bf\Pi}^{0}_{\xi}}

\newtheorem{thm} {Theorem} [section]
\newtheorem{defi} [thm] {Definition}
\newtheorem{cor} [thm] {Corollary}
\newtheorem{lem} [thm] {Lemma}
\newtheorem{prop} [thm] {Proposition}

\begin{document}

\maketitle

\centerline{$\bullet$ Universit\' e Paris 6, Institut de Math\'ematiques de Jussieu, Projet Analyse Fonctionnelle}

\centerline{Couloir 16-26, 4\`eme \'etage, Case 247, 4, place Jussieu, 75 252 Paris Cedex 05, France}

\centerline{dominique.lecomte@upmc.fr}\bigskip

\centerline{$\bullet$ Universit\'e de Picardie, I.U.T. de l'Oise, site de Creil,}

\centerline{13, all\'ee de la fa\"\i encerie, 60 107 Creil, France}\bigskip\bigskip\bigskip\bigskip\bigskip

\ufootnote{{\it 2010 Mathematics Subject Classification.}~Primary: 03E15, Secondary: 54H05}

\ufootnote{{\it Keywords and phrases.}~acyclic, Borel, reduction, dichotomy}

\noindent {\bf Abstract.} The literature provides dichotomies involving homomorphisms (like the 
$\mathbb{G}_0$ dichotomy) or reductions (like the characterization of sets potentially in a Wadge class of Borel sets, which holds on a subset of a product). However, part of the motivation behind the latter result was to get reductions on the whole product, like in the classical notion of Borel reducibility considered in the study of analytic equivalence relations. This is not possible in general. We show that, under some acyclicity (and also topological) assumptions, this is widely possible. In particular, we prove that, for any non-self dual Borel class $\bf\Gamma$, there is a concrete finite $\sqsubseteq_c$-antichain basis for the class of Borel relations, whose closure has acyclic symmetrization, and which are not potentially in $\bf\Gamma$. Along similar lines, we  provide a sufficient condition for $\sqsubseteq_c$-reducing $\mathbb{G}_0$. We also prove a similar result giving a minimum set instead of an antichain if we allow rectangular reductions.
\vfill\eject
 
\section{$\!\!\!\!\!\!$ Introduction}

$\bullet$ In [K-S-T], the authors characterize the analytic graphs having a Borel countable coloring. In order to do this, they introduce a graph ${\cal G}_0$ on the Cantor space $2^\omega$. We will consider the dissymetrized version $\mathbb{G}_0$ of ${\cal G}_0$, so that ${\cal G}_0$ is the symmetrization $s(\mathbb{G}_0)$ of the oriented graph $\mathbb{G}_0$. The following result, often called the $\mathbb{G}_0$ dichotomy, is essentially proved in [K-S-T]. All our relations will be binary.

\begin{thm} \label{KST} (Kechris, Solecki, Todor\v cevi\'c) Let $X$ be a Polish space, and $A$ be an analytic  relation on $X$. Then exactly one of the following holds:\smallskip  

(a) there is a Borel {\bf countable coloring} of $(X,A)$, i.e., a Borel function 
$c\! :\! X\!\rightarrow\!\omega$ such that $c(x)\!\not=\! c(y)$ if $(x,y)\!\in\! A$,\smallskip  

(b) there is a continuous {\bf homomorphism} from $(2^{\omega},\mathbb{G}_0)$ into $(X,A)$, i.e., a continuous function ${f\! :\! 2^\omega\!\rightarrow\! X}$ such that 
$\big( f(\alpha ),f(\beta )\big)\!\in\! A$ if $(\alpha ,\beta )\!\in\!\mathbb{G}_0$ (or, equivalently, 
$\mathbb{G}_0\!\subseteq\! (f\!\times\! f)^{-1}(A)$).\end{thm} 

 The authors conjecture the injectivity of the continuous homomorphism when (b) holds. In [L3], it is proved that this is not possible in general, considering a counter-example with countable vertical sections. However, the authors show that the injectivity is possible in several cases, in particular for acyclic graphs with $s(\mathbb{G}_0)$. In practice, we will consider acyclicity only for symmetric relations since this is what matters in our Cantor-like constructions. We will say that an arbitrary relation is \bf Acyclic\rm\ (with a capital A) if its symmetrization is acyclic. The following is also essentially proved in [K-S-T].

\begin{thm} \label{KSTinjdigr} (Kechris, Solecki, Todor\v cevi\'c) Let $X$ be a Polish space, and 
$A$ be an analytic digraph on $X$. We assume that $A$ is Acyclic. Then exactly one of the following holds:\smallskip  

(a) there is a Borel countable coloring of $(X,A)$,\smallskip  

(b) there is an injective continuous homomorphism from $(2^{\omega},\mathbb{G}_0)$ into 
$(X,A)$.\end{thm} 

\noindent $\bullet$ It is natural to ask for a reduction instead of a homomorphism in (b). Recall that if $X,Y$ are topological spaces, and $A$ (resp., $B$) is a relation on $X$ (resp., $Y$), then\bigskip
 
\leftline{$(X,A)\sqsubseteq_c(Y,B)~\Leftrightarrow$}\smallskip

\rightline{$\exists f\! :\! X\!\rightarrow\! Y\mbox{ injective continuous such that }
\big( f(x),f(y)\big)\!\in\! B\mbox{ if and only if }(x,y)\!\in\! A.$}\bigskip

\noindent In this case, we say that $f$ is an injective continuous \bf reduction\rm\ from $(X,A)$ into $(Y,B)$. If $f$ is only Borel, then we say that $(X,A)$ is {\bf Borel reducible} to $(Y,B)$ (notion widely studied when $A$ and $B$ are analytic equivalence relations). In [L3], we can find the following result:

\begin{thm} \label{Mil} (Miller) Let $X$ be a Polish space, and $A$ be an analytic oriented graph on $X$. We assume that $A$  is locally countable and Acyclic. Then exactly one of the following holds:\smallskip  

(a) there is a Borel countable coloring of $(X,A)$,\smallskip  

(b) there is an injective continuous reduction from $(2^\omega ,\mathbb{G}_0)$ into $(X,A)$.
\end{thm} 

 There is a more general version of this result in [L-M] (see Theorem 15), with the same kind of assumptions.

\vfill\eject

\noindent $\bullet$ In [L3], Theorem \ref{KST} is applied to the theory of potential complexity (notion defined in [Lo2]).
 
\begin{defi} (Louveau) Let $X$, $Y$ be Polish spaces, $B$ be a Borel subset of 
$X\!\times\! Y$, and ${\bf\Gamma}$ be a class of sets closed under continuous pre-images. We say that $B$ is \bf potentially in\it\ $\bf\Gamma$ $\big($denoted 
${B\!\in\!\mbox{pot}(\bf{\Gamma})\big)}$ if there are finer Polish topologies $\sigma$ and $\tau$ on $X$ and $Y$, respectively, such that $B$, viewed as a subset of the product 
$(X,\sigma )\!\times\! (Y,\tau )$, is in $\bf\Gamma$.\end{defi}
 
 One of the motivations for introducing this notion was that it is a natural invariant for the Borel reducibility, in the sense that a relation Borel reducible to a relation potentially in $\bf\Gamma$ has also to be potentially in $\bf\Gamma$. Theorem \ref{KST} was used in the first proof of the following result.\bigskip

\noindent\bf Notation.\rm ~The letters $X$, $Y$ will refer to some sets. We set 
$\Delta (X)\! :=\!\{ (x,y)\!\in\! X^2\mid x\! =\! y\}$.

\begin{thm} \label{potclo} Let $X,Y$ be Polish spaces, and $A,B$ be disjoint analytic subsets of 
$X\!\times\! Y$. Then exactly one of the following holds:\smallskip  

(a) the set $A$ is separable from $B$ by a potentially closed set,\smallskip  

(b) there are $f\! :\! 2^\omega\!\rightarrow\! X$, $g\! :\! 2^\omega\!\rightarrow\! Y$ continuous such that the inclusions $\mathbb{G}_0\!\subseteq\! (f\!\times\! g)^{-1}(A)$ and 
$\Delta (2^\omega )\!\subseteq\! (f\!\times\! g)^{-1}(B)$ hold.\smallskip

 Moreover, we can neither have a reduction on the whole product, nor ensure that $f$ and $g$ are injective.\end{thm} 

 This result was generalized to all non self-dual Borel classes in [L4], and to all Wadge classes of Borel sets in [L5]. For instance, the following is proved in [L4].
 
\begin{thm} \label{DeLe} (1) (Debs-Lecomte) Let $\xi\!\geq\! 1$ be a countable ordinal. Then there is a Borel relation $S$ on $2^\omega$ such that for any Polish spaces $X,Y$, and for any disjoint analytic subsets $A,B$ of $X\!\times\! Y$, exactly one of the following holds:\smallskip

(a) the set $A$ is separable from $B$ by a $\mbox{pot}(\bormxi )$ set,\smallskip
 
(b) there are $f\! :\! 2^\omega\!\rightarrow\! X$, $g\! :\! 2^\omega\!\rightarrow\! Y$ continuous such that the inclusions $S\!\subseteq\! (f\!\times\! g)^{-1}(A)$ and 
$\overline{S}\!\setminus\! S\!\subseteq\! (f\!\times\! g)^{-1}(B)$ hold.\smallskip

\noindent (2) (Debs) We cannot replace $\overline{S}\!\setminus\! S$ with $\neg S$ in (b).
\end{thm} 

 There are cycle problems behind the last assertion of Theorem \ref{potclo}, proved in [L3], and also behind Theorem \ref{DeLe}.(2). This leads to assume Acyclicity to get reduction results on the whole product, which is the goal of this paper. However, note that the Acyclicity property holds in the domain side in Theorems \ref{potclo} and \ref{DeLe}. In this paper, we will assume Acyclicity on the range side.\bigskip
 
\noindent $\bullet$ As in Theorem \ref{Mil}, we are looking for minimum sets. However, for some classes of sets, there is no minimum set but a family of minimal sets. This leads to the following.
 
\begin{defi} Let ${\cal C}$ be a class, and $\leq$ be a quasi-order (i.e., a reflexive transitive relation) on $\cal C$. We say that ${\cal B}\!\subseteq\! {\cal C}$ is\smallskip 

\noindent (1) a {\bf basis} for $\cal C$ if for any element $a$ of ${\cal C}$ there is $b$ in $\cal B$ with $b\leq a$,\smallskip

\noindent (2) an {\bf antichain} if the elements of $\cal B$ are pairwise $\leq$-incomparable.\smallskip

 If moreover ${\cal B}$ is a singleton $\{ b\}$, then we say that $b$ is {\bf minimum} among elements of ${\cal C}$.\end{defi}
 
 Intuitively, we are looking for basis as small as possible for the inclusion, i.e., for antichain basis.  In practice, $\cal C$ will always be a class of pairs of the form $(X,A)$, where is a Polish space and $A$ is a relation on $X$. The elements of our basis will be of the form $(2^\omega ,B)$ (except where indicated), and $\leq$ will always be $\sqsubseteq_c$, so that we will not mention Polish spaces, $2^\omega$ and $\sqsubseteq_c$. For example, Theorem \ref{Mil} says that 
$\mathbb{G}_0$ is minimum among analytic locally countable Acyclic oriented graphs without Borel countable coloring.\bigskip

\noindent $\bullet$ We prove the following sufficient condition for reducing $\mathbb{G}_0$.

\begin{thm} \label{G0} $\{ (1,1^2),\mathbb{G}_0,s(\mathbb{G}_0)\}$ is an antichain basis for the class of analytic relations, contained in a pot$(F_\sigma )$ symmetric acyclic relation, without Borel countable coloring. In particular,\smallskip

(i) $\mathbb{G}_0$ is minimum among analytic oriented graphs, contained in a pot$(F_\sigma )$ acyclic graph, without Borel countable coloring,\smallskip

(ii) $s(\mathbb{G}_0)$ is minimum among analytic graphs, contained in a pot$(F_\sigma )$ acyclic graph, without Borel countable coloring.\end{thm} 

 Note that this extends Theorem \ref{Mil}. Indeed, under the assumptions of Theorem \ref{Mil}, the reflexion theorem gives a Borel locally countable Acyclic digraph $B$ containing $A$. It remains to note that $B$ is $\mbox{pot}(F_\sigma )$ since a Borel set with countable vertical sections has $F_\sigma$ vertical sections and is therefore $\mbox{pot}(F_\sigma )$ (see [Lo1]). We will see that this is a real extension, in the sense that we can find a $F_\sigma$ acyclic graph $D$ on 
$2^\omega$ and Borel oriented subgraphs of $D$, without Borel countable coloring, of arbitrarily high potential complexity (see Proposition \ref{exa}). Theorem \ref{G0} applies to analytic relations whose closure is Acyclic. More generally, all the dichotomy results in this paper work for Borel relations whose closure is an Acyclic oriented graph, and for Borel graphs whose closure is an acyclic graph. We always prove more than that, in different directions.\bigskip

\noindent $\bullet$ In order to state our main theorem, we need some notation.\bigskip
 
\noindent\bf Notation.\rm ~If $s\!\in\! 2^{<\omega}$, then 
$N_s\! :=\!\{\alpha\!\in\! 2^\omega\mid s\!\subseteq\!\alpha\}$ is the associated basic clopen set.\bigskip

\noindent - The \bf dual class\rmÊ of $\bf\Gamma$ is 
$\check {\bf\Gamma}\! :=\!\{\neg A\mid A\!\in\! {\bf\Gamma}\}$. If 
${\bf\Gamma}\!\not=\!\check {\bf\Gamma}$ is a Borel class, we say that $\bf\Gamma$ is a 
{\bf non self-dual} Borel class (this means that $\bf\Gamma$ is of the form $\boraxi$ or $\bormxi$).\bigskip

\noindent - If $R$ is a relation on $2^\omega$, then $R^=\! :=\! R$, 
$R^\square\! :=\! R\cup\Delta (2^\omega )$, $R^\sqsubset\! :=\! R\cup\Delta (N_0)$ and 
$R^\sqsupset\! :=\! R\cup\Delta (N_1)$.\bigskip

\noindent - Let $A\!\subseteq\! X\!\times\! Y$. We consider the bipartite oriented graph $G_A$ on $X\!\oplus\! Y$ defined by
$$\big( (\varepsilon ,z),(\varepsilon',z')\big)\!\in\! G_A\ \Leftrightarrow\  
(\varepsilon ,\varepsilon')\! =\! (0,1)\ \wedge\ (z,z')\!\in\! A.$$
$\bullet$ We introduce a bipartite version of $\mathbb{G}_0$. We set 
$\mathbb{B}_0\! :=\!\{ (0\alpha ,1\beta )\mid (\alpha ,\beta )\!\in\!\mathbb{G}_0\}$. In particular, with a slight abuse of notation, $\mathbb{B}_0\! =\! G_{\mathbb{G}_0}$. We will repeat this abuse of notation.\bigskip

 Now we can state our main positive result.

\begin{thm} \label{generalpos} Let $\bf\Gamma$ be a non self-dual Borel class. Then there is a  concrete relation $R$ on $2^\omega$, contained in $N_0\!\times\! N_1$, satisfying the following properties.\smallskip

\noindent (1) $R$ is complete for the class of sets which are the intersection of a 
$\check {\bf\Gamma}$ set with a closed set.\smallskip

\noindent (2) If ${\bf\Gamma}\!\not=\!\boraone$, then the set\smallskip

\centerline{${\cal A}\! :=\!\big\{ A^e\mid A\!\in\!\{ R,R\cup\overline{R}^{-1},
R\cup (\overline{R}^{-1}\!\setminus\! R^{-1})\}\ \wedge\ 
e\!\in\!\{ =,\square ,\sqsubset ,\sqsupset\}\big\}\cup\big\{ s(R)^e\mid 
e\!\in\!\{ =,\square ,\sqsubset\}\big\}$}\smallskip

\noindent is an antichain made of non-pot$({\bf\Gamma})$ Acyclic relations.\smallskip

\noindent (3) If $\bf\Gamma$ is of rank at least two, then\smallskip

(i) $\cal A$ is a basis for the class of non-pot$({\bf\Gamma})$ Borel subsets of a pot$(F_\sigma )$ Acyclic relation $G$,\smallskip

(ii) $R$ is minimum among non-pot$({\bf\Gamma})$ Borel subsets of a pot$(F_\sigma )$ Acyclic oriented graph $G$,\smallskip

(iii) $s(R)$ is minimum among non-pot$({\bf\Gamma})$ Borel graphs contained in a 
pot$(F_\sigma )$ acyclic graph $G$.\smallskip

(iv) $R\cup\Delta (2^\omega )$ is minimum among non-$\mbox{pot}({\bf\Gamma})$ Borel quasi-orders (or partial orders) contained in a $\mbox{pot}(F_\sigma )$ Acyclic relation $G$.\smallskip

\noindent (4) If ${\bf\Gamma}\! =\!\bormtwo$, then\smallskip

(i) the set $\big\{ R^e\mid e\!\in\!\{ =,\square ,\sqsubset ,\sqsupset\}\big\}\cup\big\{ s(R)^e\mid 
e\!\in\!\{ =,\square ,\sqsubset\}\big\}$ is a basis for the class of non-pot$({\bf\Gamma})$ Borel locally countable Acyclic relations,\smallskip

(ii) $R$ is minimum among non-pot$({\bf\Gamma})$ Borel locally countable Acyclic oriented graphs,\smallskip

(iii) $s(R)$ is minimum among non-pot$({\bf\Gamma})$ Borel locally countable acyclic graphs.\smallskip

(iv) $R\cup\Delta (2^\omega )$ is minimum among non-$\mbox{pot}({\bf\Gamma})$ Borel locally countable Acyclic quasi-orders (or partial orders).\smallskip

\noindent (5) If ${\bf\Gamma}\! =\!\bormone$, then $R\! =\!\mathbb{B}_0$ and\smallskip

(i) the conclusions of (3).(ii), (3).(iii) and (3).(iv) remain true if $G$ is potentially closed,\smallskip

(ii) the set ${\cal A}\cup\{\mathbb{G}_0,s(\mathbb{G}_0)\}$ is an antichain basis for the class of 
non-pot$({\bf\Gamma})$ Borel subsets of a potentially closed Acyclic relation.\smallskip

\noindent (6) If ${\bf\Gamma}\! =\!\boraone$, then 
$R\! =\!\{ (0\alpha ,1\alpha )\mid\alpha\!\in\! 2^\omega\}$ and the conclusions of (3).(ii) and (3).(iii)  remain  true if the potential complexity of $G$ is arbitrary. In fact, $\{\Delta (2^\omega ),R,s(R)\}$ is an antichain basis for the class of non-pot$({\bf\Gamma})$ Borel Acyclic relations, and 
$\Delta (2^\omega )$ is minimum among non-$\mbox{pot}({\bf\Gamma})$ Borel Acyclic quasi-orders (or partial orders).\end{thm}

 Recall that a set is in the class $D_2(\boraone )$ if it is the difference of two open sets. We will see that any $\mbox{pot}\big(\check D_2(\boraone )\big)$ Acyclic relation is in fact potentially closed (see Proposition \ref{simpler}). In particular, we can replace the assumption ``potentially closed" in (5) with ``$\mbox{pot}\big(\check D_2(\boraone )\big)$". An immediate consequence of Theorem \ref{generalpos} is the following.

\begin{cor} Let $\bf\Gamma$ be a non self-dual Borel class. Then there is a concrete finite 
$\sqsubseteq_c$-antichain basis for the class of non-$\mbox{pot}({\bf\Gamma})$ Borel relations whose closure has acyclic symmetrization.\end{cor}

\noindent $\bullet$ We will state our main negative result, showing the optimality of some of the assumptions in Theorem \ref{generalpos}.\bigskip

\noindent\bf Notation.\rm ~If ${\bf\Gamma}\!\not=\!\check {\bf\Gamma}$ is a Borel class, then we denote by 
$${\bf\Gamma}\oplus\check {\bf\Gamma}\! :=\!\{ (A\cap C)\cup (B\!\setminus\! C)\mid 
A\!\in\! {\bf\Gamma},B\!\in\!\check {\bf\Gamma},C\!\in\!\borone\}$$ 
the successor of ${\bf\Gamma}$ in the Wadge quasi-order.
 
\begin{thm} \label{generalneg} Let $\bf\Gamma$ be a non self-dual Borel class.\smallskip

\noindent (1) If ${\bf\Gamma}\!\not=\!\boraone$, then there is no Acyclic oriented graph which is minimum among non-$\mbox{pot}({\bf\Gamma})$ Borel Acyclic oriented graphs.\smallskip

\noindent (2) If $\bf\Gamma$ is of rank at least two, then there is no relation which is minimum among non-$\mbox{pot}({\bf\Gamma})$ Borel subsets of a 
$\mbox{pot}({\bf\Gamma}\oplus\check {\bf\Gamma})$ Acyclic oriented graph.\smallskip

\noindent (3) If ${\bf\Gamma}\! =\!\bormone$, then there is no relation which is minimum among non-$\mbox{pot}({\bf\Gamma})$ Borel locally countable subsets of a 
$\mbox{pot}\big( D_2(\boraone )\big)$ Acyclic oriented graph.\end{thm}

 Let us precise our optimality considerations in Theorem \ref{generalpos}.\bigskip
 
\noindent (2) The assumption is optimal, because of (6). For instance, 
$\Delta (2^\omega )\sqsubseteq_c\{ (0\alpha ,1\alpha )\mid\alpha\!\in\! 2^\omega\}^\square$, but the converse fails.\bigskip
 
\noindent (3).(ii) By Theorem \ref{generalneg}.(2), the assumption ``$G$ is pot$(F_\sigma )$" is optimal for ${\bf\Gamma}\! =\!\boratwo$. We do not know whether this assumption is optimal if the rank of $\bf\Gamma$ is at least three (Theorem \ref{generalneg}.(2) just says that we cannot replace $F_\sigma$ with ${\bf\Gamma}\oplus\check {\bf\Gamma}$).\bigskip
 
\noindent (3).(i) and (3).(iii) We do not know whether the assumption on $G$ is optimal.\bigskip
 
\noindent (5) By Theorem \ref{generalneg}.(3), the class $\check D_2(\boraone )$ is optimal.\bigskip

\noindent $\bullet$ A common strategy is used to prove Theorems \ref{G0} and \ref{generalpos}.(3). In both cases, we want to build a reduction. Using some known results about injective homomorphisms (Theorem \ref{KSTinjdigr}) and injective reductions (Corollary 1.12 in [L4] and its injective version due to Debs), we work in the domain space only, with some concrete examples instead of the abstract notions of Borel chromatic number or potential Borel class. However, the injective version due to Debs is not true if the rank of $\bf\Gamma$ is at most two, because of cycle problems again. We use some injective versions in the style of Debs's one for the first Borel classes, in the acyclic case (see [L-Z]).\bigskip
  
\noindent $\bullet$ The fact of considering Borel locally countable Acyclic relations in Theorem 
\ref{generalpos}.(4) is natural if we look at Theorem \ref{Mil}, and also the assumption of Theorems \ref{G0} and \ref{generalpos}.(3). We would like to find, for each non self-dual Borel class 
$\bf\Gamma$, an antichain basis for the class of non-pot$({\bf\Gamma})$ Borel locally countable Acyclic relations. Recall that a Borel locally countable set is $\mbox{pot}(\boratwo )$. Theorem 
\ref{generalpos}.(6) solves the case ${\bf\Gamma}\! =\!\boraone$. We use an injective version of Corollary 1.12 in [L4] for ${\bf\Gamma}\! =\!\bormtwo$ in the locally countable case which improves Theorem 7 in [L2] (see [L-Z]). As a consequence, we get Theorem 
\ref{generalpos}.(4), which solves the case ${\bf\Gamma}\! =\!\bormtwo$. It remains to study the case ${\bf\Gamma}\! =\!\bormone$. Note that it is essential here to assume some acyclicity. Indeed, Theorem 5 in [L3] gives a $\sqsubseteq_c$-antichain of size continuum made of 
$D_2(\boraone )$ oriented graphs with locally countable closure which are $\sqsubseteq_c$-minimal among non-pot$(\bormone )$ Borel relations. Moreover, Theorem 19 in [L-M] shows that there is no antichain basis for the class of non-pot$(\bormone )$ $D_2(\boraone )$ oriented graphs with locally countable closure.\bigskip

\noindent $\bullet$ In order to try to extend Theorem \ref{generalpos}.(5), we introduce the following examples:
$$\mathbb{T}_0\! :=\!\big\{\big(\varepsilon\alpha ,(1\! -\!\varepsilon )\beta\big)\mid\varepsilon\!\in\! 2\ \wedge\ (\alpha ,\beta )\!\in\!\mathbb{G}_0\big\}\mbox{,}$$
$$\mathbb{U}_0\! :=\! G_{s(\mathbb{G}_0)}\cup\mathbb{T}_0.$$
Note that $s(\mathbb{T}_0)\! =\! s(\mathbb{U}_0)\! =\! s(G_{s(\mathbb{G}_0)})$. We prove the following additional dichotomy results.

\begin{thm} \label{addicho} The set ${\cal A}'\! :=\! {\cal A}\cup\{\mathbb{G}_0,s(\mathbb{G}_0)\}\cup\big\{ A^e\mid A\!\in\!\{ G_{s(\mathbb{G}_0)},\mathbb{U}_0\}\ \wedge\ 
e\!\in\!\{ =,\square ,\sqsubset ,\sqsupset\}\big\}\ \cup$\smallskip

\rightline{$\big\{ A^e\mid A\!\in\{\mathbb{T}_0,s(\mathbb{T}_0)\}\ \wedge\ 
e\!\in\!\{ =,\square ,\sqsubset\}\big\}$}\smallskip

\noindent is a $\sqsubseteq_c$-antichain made of $D_2(\boraone )$ Acyclic relations, with locally countable closure, which are $\sqsubseteq_c$-minimal among non-$\mbox{pot}(\bormone )$ relations.\end{thm}

\noindent\bf Question.\rm\ Is ${\cal A}'$ a basis for the class of non-pot$(\bormone )$ Borel Acyclic relations with locally countable closure?\bigskip
 
\noindent $\bullet$ Note that we cannot hope for a single minimum set in Theorem 
\ref{generalpos}.(3), since the pre-image of a symmetric set by a square map is symmetric. However, a positive result holds with  rectangular maps.
 
\begin{thm} \label{rect} Let $\bf\Gamma$ be a non self-dual Borel class of rank at least two. There is a $\check {\bf\Gamma}$ relation $S$ on $2^\omega$, contained in a closed set $C$ with 
$G_C$ Acyclic, such that for any Polish spaces $X,Y$, and for any Borel subset $B$ of 
$X\!\times\! Y$ contained in a $\mbox{pot}(F_\sigma )$ set $F$ with $G_F$ Acyclic, exactly one of the following holds:\smallskip  

(a) the set $B$ is $\mbox{pot}({\bf\Gamma})$,\smallskip  

(b) there are $f\! :\! 2^\omega\!\rightarrow\! X$ and $g\! :\! 2^\omega\!\rightarrow\! Y$ injective continuous such that $S\! =\! (f\!\times\! g)^{-1}(B)$.\end{thm}

 This result holds for ${\bf\Gamma}\! =\!\bormone$ when $F$ is 
$\mbox{pot}\big(\check D_2(\boraone )\big)$ (except that $S$ is not open, we can take 
$S\! =\!\mathbb{G}_0$, and the class $\check D_2(\boraone )$ is optimal), and 
${\bf\Gamma}\! =\!\boraone$, in which case $F$ does not have to be $\mbox{pot}(F_\sigma )$.\bigskip
 
\noindent $\bullet$ The paper is organized as follows. In Section 2, we prove Theorem \ref{G0}. In Section 3, we give some material concerning potential Borel classes useful for the sequel. In Section 4, we prove some general results about our antichain basis. In Sections 5-7, we prove Theorems \ref{generalpos}, \ref{generalneg}, \ref{addicho} and \ref{rect} when the rank is at least three, two and one respectively.

\section{$\!\!\!\!\!\!$ Countable Borel chromatic number}

\underline{\bf Basic facts and notions}\bigskip
 
 The reader should see [K] for the standard descriptive set theoretic notation used in this paper. 
 
\begin{defi} Let $A$ be a relation on $X$. We set $A^{-1}\! :=\!\{ (x,y)\!\in\! X^2\mid (y,x)\!\in\! A\}$, and the \bf symmetrization\it\ of $A$ is 
$s(A)\! :=\! A\cup A^{-1}$. We say that $A$ is\smallskip

(a) \bf symmetric\it\ if $A\! =\! A^{-1}$,\smallskip

(b) \bf antisymmetric\it\ if $A\cap A^{-1}\!\subseteq\!\Delta (X)$, and a \bf partial order\it\ if $A$ is an antisymmetric quasi-order,\smallskip

(c) \bf irreflexive\it\ , or a \bf digraph\it , if $A$ does not meet $\Delta (X)$, a \bf graph\it\ if $A$ is irreflexive and symmetric, an \bf oriented graph\it\ if $A$ is irreflexive and antisymmetric,\smallskip

(d) \bf acyclic\it\ if there is no injective $A$-path $(x_i)_{i\leq n}$ with $n\!\geq\! 2$ and 
$(x_n,x_0)\!\in\! A$ ($(x_i)_{i\leq n}$ is an $A$\bf -path\it\ if $(x_i,x_{i+1})\!\in\! A$ for each 
$i\! <\! n$),\smallskip

(e) \bf connected\it\ if for each $x,y\!\in\! X$ there is an $A$-path $(x_i)_{i\leq n}$ with $x_0\! =\! x$ and $x_n\! =\! y$,\smallskip

(f) \bf bipartite\it\ if there are disjoint subsets $S_0,S_1$ of $X$ such that 
$A\!\subseteq\! (S_0\!\times\! S_1)\cup (S_1\!\times\! S_0)$,\smallskip

(g) \bf locally countable\it\ if $A$ has countable horizontal and vertical sections (this also makes sense in a rectangular product $X\!\times\! Y$).\end{defi}

\vfill\eject
  
 We start with a simple algebraic fact about connected acyclic graphs.
 
\begin{lem} \label{iso} Let $G$ (resp., $H$) be an acyclic graph on $X$ (resp., $Y$), and $h$ be an injective homomorphism from $(X,G)$ into $(Y,H)$. We assume that $G$ is connected. Then $h$ is an isomorphism of graphs from $(X,G)$ onto $\big( h[X],H\cap (h[X])^2\big)$.\end{lem}

\noindent\bf Proof.\rm ~Assume that $(x,y)\!\notin\! G$. We have to see that 
$\big( h(x),h(y)\big)\!\notin\! H$. As $G$ is connected, there is $(x_i)_{i\leq n}$ injective with 
$x_0\! =\! x$, $x_n\! =\! y$, and $(x_i,x_{i+1})\!\in\! G$ if $i\! <\! n$. As $(x,y)\!\notin\! G$, 
$n\!\not=\! 1$. We may assume that $n\!\geq\! 2$. As $h$ is an injective homomorphism, 
$\big( h(x_i)\big)_{i\leq n}$ is injective and $\big( h(x_i),h(x_{i+1})\big)\!\in\! H$ if $i\! <\! n$. The acyclicity of $H$ gives the result.\hfill{$\square$}\bigskip

\noindent\bf Notation.\rm ~We have to introduce a minimum digraph without Borel countable coloring, namely $\mathbb{G}_0$.\bigskip

\noindent $\bullet$ Let $\psi\! :\!\omega\!\rightarrow\! 2^{<\omega}$ be a natural bijection. More precisely, $\psi (0)\! :=\!\emptyset$ is the sequence of length $0$, $\psi (1)\! :=\! 0$, $\psi (2)\! :=\! 1$ are the sequences of length $1$, and so on. Note that $|\psi (n)|\!\leq\! n$ if $n\!\in\!\omega$. Let 
$n\!\in\!\omega$. As $|\psi (n)|\!\leq\! n$, we can define $s_n\! :=\!\psi (n)0^{n-|\psi (n)|}$. The crucial properties of the sequence $(s_{n})_{n\in\omega}$ are the following:\bigskip

- $(s_n)_{n\in\omega}$ is \bf dense\rm\ in $2^{<\omega}$. This means that for each $s\!\in\! 2^{<\omega}$, there is $n\!\in\!\omega$ such that $s_n$ extends $s$ (denoted $s\!\subseteq\! s_n$).\smallskip

- $|s_n|\! =\! n$.\bigskip

\noindent $\bullet$ We put $\mathbb{G}_0\! :=\!\{ (s_n0\gamma ,s_n1\gamma )\mid n\!\in\!{\omega}
\ \wedge\ \gamma\!\in\! 2^{\omega}\}\!\subseteq\! 2^\omega\!\times\! 2^\omega$. Note that 
$\mathbb{G}_0$ is analytic (in fact a difference of two closed sets) since the map 
$(n,\gamma )\!\mapsto\! (s_n0\gamma ,s_n1\gamma )$ is continuous.\bigskip

\noindent $\bullet$ We identify $(2\!\times\! 2)^{<\omega}$ with 
$\bigcup_{l\in\omega}~(2^l\!\times\! 2^l)$, set 
${\cal T}\! :=\!\{ (s,t)\!\in\! (2\!\times\! 2)^{<\omega}\mid s\!\not=\! t\ \wedge\ 
(N_s\!\times\! N_t)\cap\mathbb{G}_0\!\not=\!\emptyset\}$ and, for $l\!\in\!\omega$, 
${\cal T}_l\! :=\! {\cal T}\cap (2^l\!\times\! 2^l)$. The set ${\cal T}\cup\Delta (2^{<\omega})$ is a tree with body $\overline{\mathbb{G}_0}\! =\!\mathbb{G}_0\cup\Delta (2^\omega )$. 

\begin{prop} \label{conn} Let $l\!\geq\! 1$. Then $s({\cal T}_l)$ is a connected acyclic graph on $2^l$. In particular, $\overline{\mathbb{G}_0}$ is Acyclic.\end{prop}

\noindent\bf Proof.\rm ~This comes from Proposition 18 in [L3].\hfill{$\square$}\bigskip

\noindent\bf Notation.\rm ~If $s,t\!\in\! 2^l$, then 
$p^{s,t}\! :=\! (u^{s,t}_i)_{i\leq L^{s,t}}$ is the unique injective $s({\cal T}_l)$-path from $s$ to $t$.\bigskip

 Here is another basic algebraic result about acyclicity.

\begin{lem} \label{suffacy} Let $A$ be a relation on $X$.\smallskip

(a) We assume that $A$ is irreflexive or antisymmetric, and that $A$ is Acyclic. Then $G_A$ is Acyclic.\smallskip

(b) We assume that there are disjoint subsets $X_0,X_1$ of $X$ such that 
$A\!\subseteq\! X_0\!\times\! X_1$, and that $G_A$ is Acyclic. Then $A$ is Acyclic.\end{lem}

\noindent\bf Proof.\rm ~(a) Assume first that $A$ is irreflexive. We argue by contradiction, which gives $n\!\geq\! 2$ and an injective $s(G_A)$-path $\big( (\varepsilon_i,z_i)\big)_{i\leq n}$ such that $\big( (\varepsilon_0,z_0),(\varepsilon_n,z_n)\big)\!\in\! s(G_A)$. As $A$ is Acyclic, there is 
$k\!\geq\! 1$ minimal for which there is $i\! <\! n$ such that $z_i\! =\! z_{i+k}$. As $A$ is irreflexive, $k\!\geq\! 3$. It remains to note that the $s(A)$-path $z_i,...,z_{i+k}$ contradicts the Acyclicity of 
$A$.\bigskip

 Assume now that $A$ is antisymmetric. We argue by contradiction, which gives $n\!\geq\! 2$ and an injective $s(G_A)$-path $\big( (\varepsilon_i,z_i)\big)_{i\leq n}$ such that 
$\big( (\varepsilon_0,z_0),(\varepsilon_n,z_n)\big)\!\in\! s(G_A)$. This implies that 
$\varepsilon_i\!\not=\!\varepsilon_{i+1}$ if $i\! <\! n$ and $n$ is odd. Thus $(z_i)_{i\leq n}$ is a 
$s(A)$-path such that $(z_{2j})_{2j\leq n}$ and $(z_{2j+1})_{2j+1\leq n}$ are injective and $(z_0,z_n)\!\in\! s(A)$. As $s(A)$ is acyclic, the sequence 
$(z_i)_{i\leq n}$ is not injective. We erase $z_{2j+1}$ from this sequence if $z_{2j+1}\!\in\!\{ z_{2j},z_{2j+2}\}$ and $2j\! +\!1\!\leq\! n$, which gives a sequence $(z'_i)_{i\leq n'}$ which is still a $s(A)$-path with $(z'_0,z'_{n'})\!\in\! s(A)$, and moreover satisfies $z'_i\!\not=\! z'_{i+1}$ if $i\! <\! n'$.\bigskip

 If $n'\! <\! 2$, then $n\! =\! 3$, $z_0\! =\! z_1$ and $z_2\! =\! z_3$. As $A$ is antisymmetric and 
$\varepsilon_3\! =\!\varepsilon_1\!\not=\!\varepsilon_2\! =\!\varepsilon_0$, we get $z_0\! =\! z_2$, which is absurd. If $n'\!\geq\! 2$, then $(z'_i)_{i\leq n'}$ is not injective again. We choose a subsequence of it with at least three elements,  made of consecutive elements, such that the first and the last elements are equal, and of minimal length with these properties. The Acyclicity of $A$ implies that this subsequence has exactly three elements, say $(z'_i,z'_{i+1},z'_{i+2}\! =\! z'_i)$.\bigskip

 If $z'_i\! =\! z_{2j+1}$, then $z'_{i+1}\! =\! z_{2j+2}$, $z'_{i+2}\! =\! z_{2j+4}$ and $z_{2j+3}\! =\! z_{2j+2}$.  As $A$ is antisymmetric and 
$\varepsilon_{2j+3}\! =\!\varepsilon_{2j+1}\!\not=\!\varepsilon_{2j+2}\! =\!\varepsilon_{2j+4}$, we get $z_{2j+2}\! =\! z_{2j+4}$, which is absurd. If 
$z'_i\! =\! z_{2j}$, then $z'_{i+1}\! =\! z_{2j+2}$, and $z'_{i+2}\! =\! z_{2j+3}$. As $A$ is antisymmetric and 
$\varepsilon_{2j+3}\! =\!\varepsilon_{2j+1}\!\not=\!\varepsilon_{2j+2}\! =\!\varepsilon_{2j}$, we get $z_{2j}\! =\! z_{2j+2}$, which is absurd.\bigskip

\noindent (b) Let $(z_i)_{i\leq n}$ be an injective $s(A)$-path such that $(z_0,z_n)\!\in\! s(A)$. As 
$A\!\subseteq\! X_0\!\times\! X_1$, $n$ is odd and $\big( (\varepsilon ,z_0),(1\! -\!\varepsilon ,z_1),(\varepsilon ,z_2),(1\! -\!\varepsilon ,z_3),...,(\varepsilon ,z_{n-1}),(1\! -\!\varepsilon ,z_n)\big)$ is an injective $s(G_A)$-path such that $\big( (\varepsilon ,z_0),(1\! -\!\varepsilon ,z_n)\big)\!\in\! s(G_A)$ for some $\varepsilon\!\in\! 2$.\hfill{$\square$}\bigskip

\noindent\bf Remark.\rm ~Proposition \ref{conn} says that 
$s(\overline{\mathbb{G}_0})\! =\! s\big( s(\overline{\mathbb{G}_0})\big)$ is acyclic. But 
$s(\overline{\mathbb{G}_0})$ is reflexive, and the sequence 
$\big( (0,0^\infty ), (1,0^\infty ), (0,10^\infty ), (1,10^\infty )\big)$ is a 
$s(G_{s(\overline{\mathbb{G}_0})})$-cycle. This shows that the assumption that $A$ is irreflexive or antisymmetric is useful.\bigskip

 The next result implies that the Acyclic reasonably definable relations are very small.

\begin{lem} \label{meager} Let $A$ be a $\sigma (\ana )$ relation on a Polish space $X$ such that $G_A$ is Acyclic, and $C,D$ be Cantor subsets of $X$. Then $A\cap (C\!\times\! D)$ is meager in $C\!\times\! D$.\end{lem}

\noindent\bf Proof.\rm ~We argue by contradiction, which gives homeomorphisms 
$\varphi\! :\! 2^\omega\!\rightarrow\! C$ and $\psi\! :\! 2^\omega\!\rightarrow\! D$. Then 
$(\varphi\!\times\!\psi )^{-1}(A)$ is not meager in $2^\omega\!\times\! 2^\omega$ and has the Baire property. By 19.6 in [K] we get Cantor sets $C'\!\subseteq\! C$ and $D'\!\subseteq\! D$ such that 
$C'\!\times\! D'\!\subseteq\! A$, and we may assume that they are disjoint. Take $\alpha_0\!\in\! C'$, $\alpha_1\!\in\! D'$, $\alpha_2\!\in\! C'\!\setminus\!\{\alpha_0\}$, and 
$\alpha_3\!\in\! D'\!\setminus\!\{\alpha_1\}$. Then 
$\big( (0,\alpha_0),(1,\alpha_1),(0,\alpha_2),(1,\alpha_3)\big)$ is an injective $s(G_A)$-path with 
$\big( (0,\alpha_0),(1,\alpha_3)\big)\!\in\! s(G_A)$, which contradicts the Acyclicity of $G_A$.
\hfill{$\square$}\bigskip

\noindent\underline{\bf Proof of Theorem \ref{G0}}\bigskip

 The next result will help us to prove Theorem \ref{G0} and will also be used later. 
 
\begin{thm} \label{Fsigma} Let $S$ be a $F_\sigma$ Acyclic digraph on $2^\omega$ containing 
$\mathbb{G}_0$. Then there is $f\! :\! 2^\omega\!\rightarrow\! 2^\omega$ injective continuous such that $\mathbb{G}_0\!\subseteq\! (f\!\times\! f)^{-1}(\mathbb{G}_0)\!\subseteq\! 
(f\!\times\! f)^{-1}(S)\!\subseteq\! s(\mathbb{G}_0)$.\end{thm} 

\noindent\bf Proof.\rm ~By Lemmas \ref{suffacy} and \ref{meager}, $S$ is meager, which gives a decreasing sequence $(O_n)_{n\in\omega}$ of dense open subsets of 
$2^\omega\!\times\! 2^\omega$ with $\neg S\! =\!\bigcap_{n\in\omega}~O_n$. We define 
$\varphi_n\! :\! N_{s_n0}\!\rightarrow\! N_{s_n1}$ by $\varphi_n(s_n0\gamma )\! :=\! s_n1\gamma$, so that $\mathbb{G}_0\! =\!\bigcup_{n\in\omega}~\mbox{Gr}(\varphi_n)$.

\vfill\eject

\noindent $\bullet$ We construct $\Psi\! :\! 2^{<\omega}\!\rightarrow\! 2^{<\omega}$ and 
$\delta\!\in\!\omega^\omega$ strictly increasing satisfying the following conditions:
$$\begin{array}{ll}
& (1)~\forall s\!\in\! 2^{<\omega}~~\forall\varepsilon\!\in\! 2~~
\Psi (s)\!\subsetneqq\!\Psi (s\varepsilon )\cr
& (2)~\forall l\!\in\!\omega ~~\exists k_l\!\in\!\omega ~~\forall s\!\in\! 2^l~~\vert\Psi (s)\vert\! =\! k_l\cr
& (3)~\forall n\!\in\!\omega ~~\forall v\!\in\! 2^{<\omega}~~\exists w\!\in\! 2^{<\omega}~~
\big(\Psi (s_n0v),\Psi (s_n1v)\big)\! =\! (s_{\delta (n)}0w,s_{\delta (n)}1w)\cr
& (4)~
\forall (s,t)\!\in\! (2\!\times\! 2)^{<\omega}\!\setminus\!\big( s({\cal T})\cup\Delta (2^{<\omega})\big) ~~ 
N_{\Psi (s)}\!\times\! N_{\Psi (t)}\!\subseteq\! O_{\vert s\vert}
\end{array}$$
$\bullet$ Assume that this is done. We define $f\! :\! 2^\omega\!\rightarrow\! 2^\omega$ by 
$\{ f(\alpha )\}\! =\!\bigcap_{n\in\omega}~N_{\Psi (\alpha\vert n)}$, and $f$ is continuous. In order to see that $f$ is injective, it is enough to check that $\Psi (s0)\!\not=\!\Psi (s1)$ if 
$s\!\in\! 2^{<\omega}$. Assume that $s\!\in\! 2^l$. We fix, for each $i\! <\! L\! :=\! L_{s,s_l}$, 
$n_i\! :=\! n^{s,s_l}_i\!\in\!\omega$ and $\varepsilon_i\! :=\!\varepsilon^{s,s_l}_i\!\in\! 2$ such that 
$u_{i+1}^{s,s_l}0^\infty\! =\!\varphi_{n_i}^{\varepsilon_i}(u_i^{s,s_l}0^\infty )$, so that 
$\Psi (s1)0^\infty\! =\!\varphi_{\delta (n_0)}^{-\varepsilon_0}...\varphi_{\delta (n_{L-1})}^{-\varepsilon_{L-1}}\varphi^{}_{\delta (l)}\varphi_{\delta (n_{L-1})}^{\varepsilon_{L-1}}...
\varphi_{\delta (n_0)}^{\varepsilon_0}\big(\Psi (s0)0^\infty\big)$. Thus $\Psi (s0)\!\not=\!\Psi (s1)$ since $k_{l+1}\! >\!\delta (l)\!\geq\!\mbox{sup}_{i<L}~\big( 1\! +\!\delta (n_i)\big)$. Note that 
$$\varphi_{\delta (n)}\big( f(s_n0\gamma )\big)\!\in\!
\varphi_{\delta (n)}[\bigcap_{p\in\omega}N_{\Psi (s_n0(\gamma\vert p))}]\!\subseteq\!
\bigcap_{p\in\omega}\varphi_{\delta (n)}[N_{\Psi (s_n0(\gamma\vert p))}]\! =\!
\bigcap_{p\in\omega}N_{\Psi (s_n1(\gamma\vert p))}\! =\!\{ f(s_n1\gamma )\}
\mbox{,}$$
so that $\mathbb{G}_0\!\subseteq\! (f\!\times\! f)^{-1}(\mathbb{G}_0)$.\bigskip

 Conversely, $\Delta (2^\omega )\!\subseteq\! (f\!\times\! f)^{-1}\big(\Delta (2^\omega )\big)
\!\subseteq\! (f\!\times\! f)^{-1}(\neg S)$. If 
$(\alpha ,\beta )\!\notin\! s(\mathbb{G}_0)\cup\Delta (2^\omega )$, then there is 
$n_0\!\in\!\omega$ such that 
$(\alpha\vert n,\beta\vert n)\!\notin\! s({\cal T})\cup\Delta (2^{<\omega})$ if $n\!\geq\! n_0$, so that 
$$\big( f(\alpha ),f(\beta )\big)\!\in\!
\bigcap_{n\geq n_0}~N_{\Psi (\alpha\vert n)}\!\times\! N_{\Psi (\beta\vert n)}\!\subseteq\!
\bigcap_{n\geq n_0}~O_n\!\subseteq\!\neg S.$$
$\bullet$ It remains to prove that the construction is possible. We first set 
$\Psi (\emptyset )\! :=\!\emptyset$. Assume that $\Psi [2^{\leq l}]$ satisfying (1)-(4) has been constructed, which is the case for $l\! =\! 0$. Note that $\Psi_{\vert 2^l}$ is an injective homomorphism from $s({\cal T}_l)$ into $s({\cal T}_{k_l})$, and therefore an isomorphism of graphs onto its range by Lemma \ref{iso}. Moreover, $\delta (n)\! <\! k_l$ if $n\! <\! l$. Let 
$\delta (l)\!\geq\!\mbox{sup}_{n<l}~\big( 1\! +\!\delta (n)\big)$ such that 
$\Psi (s_l)\!\subseteq\! s_{\delta (l)}$. We define temporary versions $\tilde\Psi (u\varepsilon )$ of the $\Psi (u\varepsilon )$'s by 
$\tilde\Psi (u\varepsilon )\! :=\!\Psi (u)(s_{\delta (l)}\varepsilon\! -\! s_{\delta (l)}\vert k_l)$, ensuring Conditions (1), (2) and (3).\bigskip

 For Condition (4), note that $L\! :=\! L^{s,t}\!\geq\! 2$. Here again, $\tilde\Psi_{\vert 2^{l+1}}$ is an isomorphism of graphs onto its range. This implies that $\big(\tilde\Psi (u^{s,t}_i)\big)_{i\leq L}$ is the injective $s({\cal T})$-path from $\tilde\Psi (s)$ to $\tilde\Psi (t)$. Thus 
$\big(\tilde\Psi (u^{s,t}_i)0^\infty\big)_{i\leq L}$ is the injective $s(\mathbb{G}_0)$-path (and also 
$s(S)$-path) from $\tilde\Psi (s)0^\infty$ to $\tilde\Psi (t)0^\infty$. Therefore 
$\big(\tilde\Psi (s)0^\infty ,\tilde\Psi (t)0^\infty\big)\!\in\!\neg s(S)\!\subseteq\! O_{l+1}$ 
since $L\!\geq\! 2$. This gives $m\!\in\!\omega$ with  
${N_{\tilde\Psi (s)0^m}\!\times\! N_{\tilde\Psi (t)0^m}\!\subseteq\! O_{l+1}}$. It remains to set 
$\Psi'(u\varepsilon )\! :=\!\tilde\Psi (u\varepsilon )0^m$, which ensures  the inclusion 
$N_{\Psi'(s)}\!\times\! N_{\Psi'(t)}\!\subseteq\! O_{l+1}$.\hfill{$\square$}

\begin{cor} \label{corFsigma} Let $X$ be a Polish space, $A$ be an analytic subset of a 
$\mbox{pot}(F_\sigma )$ Acyclic digraph $G$ on $X$. Then exactly one of the following holds:\smallskip

(a) there is a Borel countable coloring of $(X,A)$,\smallskip

(b) there is $f\! :\! 2^\omega\rightarrow\! X$ injective continuous with $\mathbb{G}_0
\!\subseteq\! (f\!\times\! f)^{-1}(A)\!\subseteq\! (f\!\times\! f)^{-1}(G)\!\subseteq\! s(\mathbb{G}_0)$.
\end{cor} 

\noindent\bf Proof.\rm ~By Theorem \ref{KST}, (a) and (b) cannot hold simultaneously. So assume that (a) does not hold. Let $\tau$ be a finer Polish topology on $X$ such that 
$G\!\in\! F_\sigma\big( (X,\tau )^2\big)$. Theorem \ref{KSTinjdigr} gives 
$g\! :\! 2^\omega\rightarrow\! (X,\tau )$ injective continuous with 
$\mathbb{G}_0\!\subseteq\! (g\!\times\! g)^{-1}(A)$. We now apply Theorem \ref{Fsigma} to 
$S\! :=\! (g\!\times\! g)^{-1}(G)$, which gives $h\! :\! 2^\omega\rightarrow\! 2^\omega$ injective continuous with $\mathbb{G}_0\!\subseteq\! (h\!\times\! h)^{-1}(\mathbb{G}_0)\!\subseteq\! 
(h\!\times\! h)^{-1}(S)\!\subseteq\! s(\mathbb{G}_0)$. It remains to set $f\! :=\! g\circ h$.
\hfill{$\square$}

\vfill\eject

\noindent\bf Proof of Theorem \ref{G0}.\rm ~By Theorem \ref{KST}, $1^2$, $\mathbb{G}_0$ and 
$s(\mathbb{G}_0)$ are in the context of Theorem \ref{G0}. Assume that $A$ is an analytic relation on a Polish space $X$, without Borel countable coloring, contained in a pot$(F_\sigma )$ symmetric acyclic relation $S$. If $A$ is not irreflexive, then let $(x,x)\!\in\! A$, and $0\!\mapsto\! x$ is a witness for $(1,1^2)\sqsubseteq_c(X,A)$. So we may assume that $A$ and $S$ are irreflexive. Corollary \ref{corFsigma} gives $f\! :\! 2^\omega\!\rightarrow\! X$ with 
$\mathbb{G}_0\!\subseteq\! A'\! :=\! (f\!\times\! f)^{-1}(A)\!\subseteq\! s(\mathbb{G}_0)$. By Theorem \ref{KSTinjdigr} again, two cases can happen.\bigskip

\noindent $\bullet$ Either there is a Borel countable coloring of 
$R\! :=\! A'\setminus\! <_{\mbox{lex}}$. This gives a non-meager $R$-discrete $G_\delta$ subset $G$ of $2^\omega$. Note that $A'\cap G^2$ is an analytic oriented graph on $G$ without Borel countable coloring and $(f\!\times\! f)^{-1}(S)\cap G^2$ is a $\mbox{pot}(F_\sigma )$ acyclic graph containing $A'\cap G^2$. Corollary \ref{corFsigma} gives $g\! :\! 2^\omega\!\rightarrow\! G$ injective continuous with 
$\mathbb{G}_0\!\subseteq\! (g\!\times\! g)^{-1}(A'\cap G^2)\!\subseteq\! s(\mathbb{G}_0)$. Thus 
$(2^\omega ,\mathbb{G}_0)\sqsubseteq_c(X,A)$ since $A'\cap G^2$ is an oriented graph.\bigskip

\noindent $\bullet$ Or there is $h\! :\! 2^\omega\!\rightarrow\! 2^\omega$ injective continuous with $\mathbb{G}_0\!\subseteq\! (h\!\times\! h)^{-1}(R)$. Note that 
$A''\! :=\! (h\!\times\! h)^{-1}(A')$ is analytic, contains $s(\mathbb{G}_0)$, and is contained 
$S'\! :=\! (h\!\times\! h)^{-1}\big( (f\!\times\! f)^{-1}(S)\big)$, which is a $\mbox{pot}(F_\sigma )$ acyclic graph.\bigskip

 Indeed, if $(\alpha ,\beta )\!\in\! s(\mathbb{G}_0)\!\setminus\!\mathbb{G}_0$, then 
$(\alpha ,\beta )\!\in\!\mathbb{G}_0^{-1}$, 
$\big( h(\beta ),h(\alpha )\big)\!\in\! A'\setminus\! <_{\mbox{lex}}\subseteq\! s(\mathbb{G}_0)\!\setminus\!\mathbb{G}_0\! =\!\mathbb{G}_0^{-1}$, and 
$\big( h(\alpha ),h(\beta )\big)\!\in\!\mathbb{G}_0\!\subseteq\! A'$. Corollary \ref{corFsigma} gives 
$i\! :\! 2^\omega\!\rightarrow\! 2^\omega$ with 
$$\mathbb{G}_0\!\subseteq\! 
(i\!\times\! i)^{-1}\big( s(\mathbb{G}_0)\big)\!\subseteq\! (i\!\times\! i)^{-1}(S')\!\subseteq\! 
s(\mathbb{G}_0).$$ 
Thus $s(\mathbb{G}_0)\!\subseteq\! (i\!\times\! i)^{-1}(A'')\!\subseteq\! s(\mathbb{G}_0)$ and 
$\big( 2^\omega ,s(\mathbb{G}_0)\big)\sqsubseteq_c(X,A)$.\hfill{$\square$}\bigskip

\noindent\bf Question.\rm ~Can we extend Theorem \ref{G0} to any acyclic graph?\bigskip

 The next remark essentially says that Theorem \ref{G0} applies to analytic relations whose tree has Acyclic levels. 

\begin{prop} \label{lev} Let $X$ be a Polish space, $C$ be a closed subset of the Baire space, 
$b\! :\! C\!\rightarrow\! X$ be a continuous bijection, and $A$ an analytic relation on $X$. We assume that the levels of the tree of $s\big(\overline{(b\!\times\! b)^{-1}(A)}\big)$ are acyclic. Then $A$ is contained in a $\mbox{pot}(\bormone )$ symmetric acyclic relation.\end{prop}

\noindent\bf Proof.\rm ~The levels of the tree of $s\big(\overline{(b\!\times\! b)^{-1}(A)}\big)$ are defined, for each $l\!\in\!\omega$, by 
$$L_l\! :=\!\{ (s,t)\!\in\!\omega^l\!\times\!\omega^l\mid (N_s\!\times\! N_t)\cap
s\big(\overline{(b\!\times\! b)^{-1}(A)}\big)\!\not=\!\emptyset\} .$$ 
As they are acyclic, $s\big(\overline{(b\!\times\! b)^{-1}(A)}\big)$ is acyclic too. Thus 
$s\big(\overline{(b\!\times\! b)^{-1}(A)}\big)$ is a closed symmetric acyclic relation containing 
$(b\!\times\! b)^{-1}(A)$. We are done since $b$ is a Borel isomorphism.\hfill{$\square$}

\section{$\!\!\!\!\!\!$ Potential Borel classes}

\bf Notation.\rm ~Fix some standard bijection $<.,.>:\!\omega^2\!\rightarrow\!\omega$, for example 
$$(n,p)\!\mapsto <n,p>:=\!\frac{(n\! +\! p)(n\! +\! p\! +\! 1)}{2}\! +\! p.$$ 
Let $I\! :\!\omega\!\rightarrow\!\omega^2$ be its inverse ($I$ associates $\big( (l)_0,(l)_1\big)$ with $l$).\bigskip

 We identify $(2^l)^2$ and $(2^2)^l$, for each $l\!\in\!\omega\! +\! 1$.
 
\begin{defi} Let ${\cal F}\!\subseteq\!\bigcup_{l\in\omega}~(2^l)^2\!\equiv\! (2^2)^{<\omega}$. We say that ${\cal F}$ is a \bf frame\it\ if\medskip

\noindent (1) $\forall l\!\in\!\omega~\exists ! (u_l,v_l)\!\in\! {\cal F}\!\cap\! (2^l)^2$,\smallskip

\noindent (2) $\forall p,q\!\in\!\omega~\forall w\!\in\! 2^{<\omega}~\exists N\!\in\!\omega~
(u_q0w0^N,v_q1w0^N)\!\in\! {\cal F}$ and $(|u_q0w0^N|)_0\! =\! p$,\smallskip

\noindent (3) $\forall l\! >\! 0~\exists q\! <\! l~\exists w\!\in\! 2^{<\omega}\ (u_l,v_l)\! =\! (u_q0w,v_q1w)$.\medskip

 If ${\cal F}\! =\!\{ (u_l,v_l)\mid l\!\in\!\omega\}$ is a frame, then we will call $T$ the tree on $2^2$ generated by $\cal F$: 
$$T\! :=\!\big\{ (u,v)\!\in\! (2\!\times\! 2)^{<\omega}\mid u\! =\!\emptyset\vee
\big(\exists q\!\in\!\omega ~\exists w\!\in\! 2^{<\omega}~(u,v)\! =\! (u_q0w,v_q1w)\big)\big\}.$$
\end{defi}

 The existence condition in (1) and the density condition (2) ensure that $\lceil T\rceil$ is big enough to contain sets of arbitrary high complexity. The uniqueness condition in (1) and condition (3) ensure that $\lceil T\rceil$ is small enough to make the reduction in Theorem \ref{DL} to come possible. The last part of condition (2) gives a control on the verticals which is very useful to construct complex examples. This definition is a bit different from Definition 2.1 in [L5], where 
$(|u_q0w0^N|\! -\! 1)_0$ is considered instead of $(|u_q0w0^N|)_0$ in Condition (2). This new notion is simpler and more convenient to study  the equivalence relations associated with ideals (see [C-L-M] for a use of this kind of equivalence relations). In most cases, our examples will be ideals (see Lemma \ref{exist}). Also, we do not need Condition (d) in [L5] ensuring that 
$T\cap (d^d)^l$ is $\Borel$ when $d\! =\!\omega$, which is clear when $d\! =\! 2$.\bigskip

\noindent\bf Notation.\rm ~We set, for $l\!\in\!\omega$, $M(l)\! :=\!\mbox{max}\{ m\!\in\!\omega\mid\frac{m(m\! +\! 1)}{2}\!\leq\! l\}$, so that 
$M(l)\! =\! (l)_{0}\! +\! (l)_{1}$.
 
\begin{lem} There is a frame.\end{lem}

\noindent\bf Proof.\rm ~We first set $(u_0,v_0)\! :=\! (\emptyset ,\emptyset )$. Note that  
$$(l)_{0}\! +\! (l)_{1}\! =\! M(l)\!\leq\!\frac{M(l)(M(l)\! +\! 1)}{2}\!\leq\! l\mbox{,}$$
for each $l\!\in\!\omega$. This allows us to define
$$(u_{l+1},v_{l+1})\! :=\! (u_{ ( (l)_1)_0}~0~\psi ( ( (l)_1)_1 )~0^{ l-( (l)_1)_0-| \psi ( ( (l)_1)_1 ) |},
v_{ ( (l)_1)_0}~1~\psi ( ( (l)_1)_1 )~0^{ l-( (l)_1)_0-| \psi ( ( (l)_1)_1 ) |}).$$ 
Note that $(u_l,v_l)$ is well defined and $|(u_l,v_l)|\! =\! l$, by induction on $l$. It remains to check that condition (2) in the definition of a frame is fulfilled. We set $n\! :=\!\psi^{-1}(w)$, and 
$l\! :=\!\big<\ p\! +\! 1,<q,n>\!\big>$. It remains to put $N\!:=\! l\! -\! q\! -\! |w|$: 
$(u_q0w0^N,v_q1w0^N)\! =\! (u_{l+1},v_{l+1})$, and 
$$(|u_q0w0^N|)_0\! =\! (l\! +\! 1)_0\! =\! (<\! p,<q,n>\! +1\! >)_0=\! p.$$
This finishes the proof.\hfill{$\square$}\bigskip

 In the sequel, $T$ will be the tree generated by a fixed frame $\cal F$. We set, for each $l\!\in\!\omega$, 
$$T_l\! :=\! T\cap (2^l\times 2^l).$$ 
The proof of Proposition 3.2 in [L4] shows that $s(G_{T_l})$ is an acyclic graph if $l\!\in\!\omega$, and Lemma \ref{suffacy} shows that $s(T_l)$ is acyclic if $l\!\geq\! 1$ since $\lceil T\rceil\!\subseteq\! N_0\!\times\! N_1$ (it is also connected, by induction on $l$). Using Theorem 1.10 in [L4], this gives the next result, without the injectivity complement due to Debs.
 
\begin{thm} \label{DL} Let $\bf\Gamma$ be a non self-dual Borel class, 
$S\!\in\!\check  {\bf\Gamma}(\lceil T\rceil )$, $X,Y$ be Polish spaces, and $A,B$ be disjoint analytic subsets of $X\!\times\! Y$.\smallskip

\noindent (1) (Debs-Lecomte) One of the following holds:\smallskip  

(a) the set $A$ is separable from $B$ by a $\mbox{pot}({\bf\Gamma})$ set,\smallskip  

(b) there are $f\! :\! 2^\omega\!\rightarrow\! X$ and $g\! :\! 2^\omega\!\rightarrow\! Y$ continuous such that the inclusions $S\!\subseteq\! (f\!\times\! g)^{-1}(A)$ and 
$\lceil T\rceil\!\setminus\! S\!\subseteq\! (f\!\times\! g)^{-1}(B)$ hold.\smallskip

\noindent (2) (Debs) If moreover $\bf\Gamma$ is of rank at least three, then we can have $f$ and 
$g$ injective in (b).\end{thm}

\noindent\bf Notation.\rm ~We use complex one-dimensional sets to build complex two-dimensional sets, using the symmetric difference. More precisely, recall that the {\bf symmetric difference} $\alpha\Delta\beta$ of $\alpha ,\beta\!\in\! 2^\omega$ is the element of $2^\omega$ defined by 
$(\alpha\Delta\beta )(m)\! =\! 1$ exactly when $\alpha (m)\!\not=\!\beta (m)$. We associate the following two-dimensional sets to the one-dimensional set 
${\cal I}\!\subseteq\! 2^\omega$. We set 
$$E_{\cal I}\! :=\!\{ (\alpha ,\beta )\!\in\! 2^\omega\!\times\! 2^\omega\mid\alpha\Delta\beta\!\in\! {\cal I}\}$$ 
and $S_{\cal I}\! :=\!\lceil T\rceil\cap E_{\cal I}$. If $\cal I$ is a nonempty {\bf ideal} (i.e., $\cal I$ is closed under taking subsets and finite unions), then 
$E_{\cal I}$ is the equivalence relation associated with $\cal I$. The following result ensures that $S_{\cal I}$ is complicated if $\cal I$ is.

\begin{defi} Let ${\cal I}\!\subseteq\! 2^\omega$, $2^\omega$ being identified with the power set of 
$\omega$. We say that $\cal I$ is \bf vertically\ invariant\it\ if, whenever 
$i\! :\!\omega\!\rightarrow\!\omega$ is injective such that $\big( i(m)\big)_0\! =\! (m)_0$ for each 
$m\!\in\!\omega$, then, for each $N\!\subseteq\!\omega$, 
$N\!\in\! {\cal I}\Leftrightarrow i[N]\!\in\! {\cal I}$.\end{defi}
 
 Recall that 
$\mathbb{E}_0\! :=\!\{ (\alpha ,\beta )\!\in\! 2^\omega\!\times\! 2^\omega\mid\exists m\!\in\!\omega ~~
\forall n\!\geq\! m~~\alpha (n)\! =\!\beta (n)\}$. 
 
\begin{thm} \label{compl} Let $\bf\Gamma$ be a non self-dual Borel class, 
${\cal I}\!\subseteq\! 2^\omega$ be a vertically and $\mathbb{E}_0$-invariant true 
$\check {\bf\Gamma}$ set, $(u,v)\!\in\! T$ and $G$ be a dense $G_\delta$ subset of $2^\omega$. Then $S_{\cal I}\cap\big( (N_u\cap G)\!\times\! (N_v\cap G)\big)$ is not separable from its complement in $\lceil T\rceil$ by a $\mbox{pot}({\bf\Gamma})$ set.\end{thm}

 This is essentially Lemma 2.6 in [L5], when $s\! =\!\emptyset$ and $G\! =\! 2^\omega$. The general proof is very similar, but we give it for completeness. The first part of the next definition gives the objects expressing the complexity of $S_{\cal I}$ on some generic vertical $(S_{\cal I})_\alpha$. The second part gives a condition on $\cal I$ which is sufficient to ensure the complexity of $S_{\cal I}$, together with a topological complexity condition.

\begin{defi} Let $n\!\in\!\omega\!\setminus\!\{ 0\}$, $\alpha\!\in\! 2^\omega$, 
$F\! :\! 2^\omega\!\rightarrow\! 2^\omega$, and ${\cal I}\!\subseteq\! 2^\omega$. We say that\medskip

\noindent (a) $(n,\alpha ,F)$ is a \bf transfer\ triple\it\ if, for any $\beta\!\in\! 2^\omega$, there is an injection $i\! :\!\omega\!\rightarrow\!\omega$ such that 
$$i[\{ m\!\in\!\omega\mid\beta (m)\! =\! 1\} ]\! =\!\big\{ k\!\geq\! n\mid\big(\alpha\Delta F(\beta )\big)(k)\! =\! 1\big\}\mbox{,}$$ 
and $\big( i(m)\big)_0\! =\! (m)_0$ if $m\!\in\!\omega$,\medskip

\noindent (b) ${\cal I}$ is \bf transferable\it\ if 
$\beta\!\in\! {\cal I}\ \Leftrightarrow\ \alpha\Delta F(\beta )\!\in\! {\cal I}$ for any transfer triple 
$(n,\alpha ,F)$ and any $\beta\!\in\! 2^\omega$,\medskip

\noindent (c) $\cal I$ is \bf weakly transferable\it\ if 
$\beta\!\in\! {\cal I}\ \Leftrightarrow\ \alpha\Delta F(\beta )\!\in\! {\cal I}$ for any transfer triple 
$(1,\alpha ,F)$ and any $\beta\!\in\! 2^\omega$. \end{defi}

 We could also mention $\{ m\!\in\!\omega\mid\beta (m)\! =\! 0\}$, but we really care about the value $1$. The reason why we wrote 
``$n\!\in\!\omega\!\setminus\!\{ 0\}$"  is that $(\gamma\Delta\delta )(0)\! =\! 1$ if $(\gamma ,\delta )\!\in\!\lceil T\rceil$. The following lemma is the key ingredient in the proof of Theorem \ref{compl}.

\begin{lem} \label{key} Let $(u,v)\!\in\! T$ and $G$ be a dense $G_{\delta}$ subset $2^\omega$. Then we can find $n\!\in\!\omega\!\setminus\!\{ 0\}$, $\alpha\!\in\! N_{u_n}\cap G$ and 
$F\! :\! 2^\omega\!\rightarrow\! N_{v_n}\cap G$ continuous such that\smallskip

\noindent (a) $(u,v)\!\subseteq\! (u_n,v_n)$,\smallskip

\noindent (b) for any $\beta\!\in\! 2^\omega$, $\big(\alpha ,F(\beta )\big)\!\in\!\lceil T\rceil$,\smallskip

\noindent (c) $(n,\alpha ,F)$ is a transfer triple.\smallskip

 If moreover $u\! =\!\emptyset$, then we can have $n\! =\! 1$.\end{lem}

\noindent\bf Proof.\rm ~We set 
$(u',v')\! :=\!\left\{\!\!\!\!\!\!\!\!
\begin{array}{ll}
& (0,1)\mbox{ if }u\! =\!\emptyset\mbox{,}\cr
& (u,v)\mbox{ otherwise.}
\end{array}
\right.$
Let $M\!\in\!\omega$ be such that $(u'0^M,v'0^M)\!\in\! {\cal F}$ and ${(\vert u'\vert\! +\! M)_0\! =\! (0)_0}$. We set 
$n\! :=\!\left\{\!\!\!\!\!\!\!\!
\begin{array}{ll}
& 1\mbox{ if }u\! =\!\emptyset\mbox{,}\cr
& \vert u'\vert\! +\! M\mbox{ otherwise.}
\end{array}
\right.$
Let $(O_q)_{q\in\omega}$ be a decreasing sequence of dense open subsets of $2^\omega$ whose intersection is $G$. We construct finite approximations of $\alpha$ and $F$. The idea is to linearize the binary tree $2^{<\omega}$. This is the reason why we will use the bijection $\psi$ defined in the introduction. In order to construct $F(\beta )$, we have to imagine, for each length 
$l$, the different possibilities for $\beta\vert l$. More precisely, we construct a map 
$l\! :\! 2^{<\omega}\!\rightarrow\!\omega\!\setminus\!\{ 0\}$. We want the map $l$ to satisfy the following conditions:
$$\begin{array}{ll}
& (1)~l(\emptyset )\! =\!\vert u'\vert\! +\! M\cr
& (2)~\forall w\!\in\! 2^{<\omega}\!\setminus\!\{\emptyset\}\ \ N_{u_{l(w)}}\cup N_{v_{l(w)}}
\!\subseteq\! O_{\vert w\vert}\cr
& (3)~\forall w\!\in\! 2^{<\omega}\ \ \forall\varepsilon\!\in\! 2\ \ 
\exists z_{w\varepsilon}\!\in\! 2^{<\omega}\ \ 
(u_{l(w\varepsilon )},v_{l(w\varepsilon )})\! =\! 
(u_{l(w)}0z_{w\varepsilon},v_{l(w)}\varepsilon z_{w\varepsilon})\cr
& (4)~\forall r\!\in\!\omega\ \ u_{l(\psi (r))}0\!\subseteq\! u_{l(\psi (r+1))}\cr
& (5)~\forall w\!\in\! 2^{<\omega}\ \ \big( l(w)\big)_0\! =\! (\vert w\vert)_0
\end{array}$$ 
$\bullet$ Assume that this construction is done. As $u_{l(0^q)}\subsetneqq u_{l(0^{q+1})}$ for each natural number $q$, we can define $\alpha\! :=\!\mbox{sup}_{q\in\omega}~u_{l(0^q)}$. Similarly, as $v_{l(\beta\vert q)}\subsetneqq v_{l(\beta\vert (q+1))}$, we can define 
$F(\beta )\! :=\!\mbox{sup}_{q\in\omega}~v_{l(\beta\vert q)}$, and $F$ is continuous. Note that 
$\alpha\!\in\!\bigcap_{q\in\omega}\ N_{u_{l(0^q)}}\!\subseteq\! N_{u_{l(\emptyset )}}\cap
\bigcap_{q>0}\ O_q\!\subseteq\! N_{u_n}\cap G$. Similarly, 
$F(\beta )\!\in\!\bigcap_{q\in\omega}\ N_{v_{l(\beta\vert q)}}\!\subseteq\! N_{v_{l(\emptyset )}}\cap\bigcap_{q>0}\ O_q\!\subseteq\! N_{v_n}\cap G$.\bigskip

\noindent (b) Note first that $l(w)\!\geq\!\vert w\vert$ since $l(w\varepsilon )\! >\! l(w)$. Fix 
$q\!\in\!\omega$. We have to see that $\big(\alpha ,F(\beta )\big)\vert q\!\in\! T$. Note that 
$u_{l(w)}\!\subseteq\!\alpha$ since 
$u_{l(0^{\vert w\vert})}\!\subseteq\! u^{}_{l(w)}\!\subseteq\! u_{l(0^{\vert w\vert +1})}$. Thus 
$\big(\alpha ,F(\beta )\big)\vert l(\beta\vert q)\! =\! (u_{l(\beta\vert q)},v_{l(\beta\vert q)})\!\in\! 
{\cal F}$. This implies that 
$\big(\alpha ,F(\beta )\big)\vert l(\beta\vert q)\!\in\! T$. We are done since $l(\beta\vert q)\!\geq\! q$.\bigskip

\noindent (c) Assume that $m\!\in\!\omega$ and $\beta (m)\! =\! 1$. We set $w\!:=\!\beta\vert m$, so that ${v_{l(w)}1\!\subseteq\! v_{l(w1)}\! =\! v_{l(\beta\vert (m+1))}\!\subseteq\! F(\beta )}$. As 
$\big( l(w)\big)_{0}\! =\! (m)_0$, $k\! :=\! l(w)\!\geq\! n$ and $(k)_0\! =\! (m)_0$. But 
$u_{l(w)}0\!\subseteq\! u_{l(w1)}\!\subseteq\!\alpha$, so that $\alpha\big( l(w)\big)$ is different from $F(\beta )\big( l(w)\big)$.\bigskip

 Assume that $k\!\geq\! n$ and $\alpha (k)\!\not=\! F(\beta )(k)$. Note that the only coordinates where $\alpha$ and $F(\beta )$ can differ are below $n$ or one of the $l(\beta\vert q)$'s. This gives $m$ with $k\! =\! l(\beta\vert m)$, and $(m)_{0}\! =\!\big( l(\beta\vert m)\big)_{0}\! =\! (k)_0$. Note that 
$\alpha\big( l(\beta\vert m)\big)\! =\! u_{l(\beta\vert (m+1))}\big( l(\beta\vert m)\big)\! =\! 0\!\not=\! 
F(\beta )\big( l(\beta\vert m)\big)\! =\! v_{l(\beta\vert (m+1))}\big( l(\beta\vert m)\big)\! =\!\beta (m)$. So $\beta(m)\! =\! 1$.\bigskip

 Now it is clear that the formula $i(m)\! :=\! l(\beta\vert m)$ defines the injection we are looking for.
 
\vfill\eject\bigskip

\noindent $\bullet$ So let us prove that the construction is possible. We construct $l(w)$ by induction on $\psi^{-1}(w)$.\bigskip

 We first choose $x\!\in\! 2^{<\omega}$ such that $N_{u_{l(\emptyset )}0x}\!\subseteq\! O_1$ and 
$y\!\in\! 2^{<\omega}$ such that $N_{v_{l(\emptyset )}0xy}\!\subseteq\! O_1$. Then we choose 
$L\!\in\!\omega$ with $(u_{l(\emptyset )}0xy0^L,v_{l(\emptyset )}0xy0^L)\!\in\! {\cal F}$ and 
$(|u_{l(\emptyset )}0xy0^L|)_0\! =\! (1)_0$. We put $z_0\! :=\! xy0^L$ and 
$l(0)\! :=\! l(\emptyset )\! +\! 1\! +\!\vert z_0\vert$. Assume that 
$\big(l(w)\big)_{\psi^{-1}(w)\leq r}$ satisfying (1)-(5) have been constructed, which is the case for $r\! =\! 1$.\bigskip

 Fix $s\!\in\! 2^{<\omega}$ and $\varepsilon\!\in\! 2$ such that $\psi (r\! +\! 1)\! =\! s\varepsilon$, with $r\!\geq\! 1$. Note that $\psi^{-1}(s)\! <\! r$, so that $l(s)\! <\! l\big(\psi (r)\big)$, by induction assumption. We set $t\! :=\!\big( u_{l(\psi (r))}\! -\! u_{l(\psi (r))}\vert (l(s)\! +\! 1)\big) 0$. We choose 
$x'\!\in\! 2^{<\omega}$ such that $N_{u_{l(s)}0tx'}\!\subseteq\! O_{\vert s\vert +1}$ and 
$y'\!\in\! 2^{<\omega}$ such that $N_{v_{l(s)}\varepsilon tx'y'}\!\subseteq\! O_{\vert s\vert +1}$. Then we choose $N\!\in\!\omega$ such that 
$(u_{l(s)}0tx'y'0^N,v_{l(s)}\varepsilon tx'y'0^N)\!\in\! {\cal F}$ and 
$\big( l(s)\! +\! 1\! +\!\vert tx'y'\vert\! +\! N\big)_{0}\! =\! (\vert s\vert\! +\! 1)_0$. We put 
$z_{s\varepsilon}\! :=\! tx'y'0^N$ and 
$l(s\varepsilon )\! :=\! l(s)\! +\! 1\! +\!\vert z_{s\varepsilon}\vert$.\hfill{$\square$}\bigskip

\noindent\bf Proof of Theorem \ref{compl}.\rm ~Let us prove that $\cal I$ is transferable. Let $(n,\alpha ,F)$ be a transfer triple, and $\beta$ in $2^\omega$. This gives an injection $i\! :\!\omega\!\rightarrow\!\omega$ with $\big( i(m)\big)_0\! =\! (m)_0$ if $m\!\in\!\omega$. We set 
$A\! :=\!\{ m\!\in\!\omega\mid\beta (m)\! =\! 1\}$. As $\cal I$ is vertically invariant, $A\!\in\! {\cal I}$ is equivalent to $i[A]\!\in\! {\cal I}$. But 
$i[A]\! =\!\{ k\!\geq\! n\mid\big(\alpha\Delta F(\beta )\big)(k)\! =\! 1\}$. As $\cal I$ is $\mathbb{E}_0$-invariant, $i[A]\!\in\! {\cal I}$ is equivalent to 
$\{ k\!\in\!\omega\mid\big(\alpha\Delta F(\alpha )\big)(k)\! =\! 1\}\!\in\! {\cal I}$, so that 
$$\beta\!\in\! {\cal I}\ \Leftrightarrow\ A\!\in\! {\cal I}\ \Leftrightarrow\ 
\{ k\!\in\!\omega\mid\big(\alpha\Delta F(\beta )\big)(k)\! =\! 1\}\!\in\! {\cal I}\ \Leftrightarrow\ \alpha\Delta F(\beta )\!\in\! {\cal I}.$$
Thus $\cal I$ is transferable.\bigskip

 We argue by contradiction. This gives $P\!\in\!\mbox{pot}({\bf\Gamma})$, and a dense $G_{\delta}$ subset $H$ of $2^\omega$ such that 
$P\cap H^2\!\in\! {\bf\Gamma}(H^2)$. Lemma \ref{key} provides $n\!\in\!\omega\!\setminus\!\{ 0\}$ such that $(u,v)\!\subseteq\! (u_n,v_n)$, 
$\alpha\!\in\! N_{u_n}\cap G\cap H$ and $F\! :\! 2^\omega\!\rightarrow\! N_{v_n}\cap G\cap H$ continuous. We set 
$S\! :=\! S_{\cal I}\cap\big( (N_{u_n}\cap G\cap H)\!\times\! (N_{v_n}\cap G\cap H)\big)$. Then 
$S\!\subseteq\! P\cap H^2\cap (N_{u_n}\!\times\! N_{v_n})\!\subseteq\!\neg \lceil T\rceil\cup S$. We set $D\! :=\!\big\{\beta\!\in\! 2^\omega\mid\big(\alpha,F(\beta )\big)\!\in\! P\cap H^2\big\}$. Then 
$D\!\in\! {\bf\Gamma}$. Let us prove that ${\cal I}\! =\! D$, which will contradict the fact that 
${\cal I}\!\notin\! {\bf\Gamma}$. Let $\beta\!\in\! 2^\omega$. As ${\cal I}$ is transferable, 
$\beta\!\in\! {\cal I}$ is equivalent to $\alpha\Delta F(\beta )\!\in\! {\cal I}$. Thus  
$$\beta\!\in\! {\cal I}\Rightarrow\alpha\Delta F(\beta )\!\in\! {\cal I}\Rightarrow
\big(\alpha , F(\beta )\big)\!\in\! S\!\subseteq\! P\cap H^2\Rightarrow\beta\!\in\! D.$$
Similarly, $\beta\!\notin\! {\cal I}\Rightarrow\beta\!\notin\! D$, and ${\cal I}\! =\! D$.\hfill{$\square$}\bigskip

\noindent\bf Notation.\rm ~In Theorem \ref{compl}, if $s\! =\!\emptyset$ and $G\! =\! 2^\omega$, then we do not need to assume that $\cal I$ is $\mathbb{E}_0$-invariant. It is enough to assume that $\cal I$ is invariant under the following map. Let $h_0\! :\! 2^\omega\!\rightarrow\! 2^\omega$ be the map defined by 
$h_0(\alpha )\! :=<1\! -\!\alpha (0),\alpha (1),\alpha (2),...>$. Note that $\mbox{Gr}(h_0)$ is a subgraph of $s(\mathbb{G}_0)$, so that it is acyclic. Similarly, we define $h_0(s)$ when $\emptyset\!\not=\! s\!\in\! 2^{<\omega}$.

\begin{cor} \label{corDL} Let $\bf\Gamma$ be a non self-dual Borel class of rank at least two, 
${\cal I}\!\subseteq\! 2^\omega$ be a vertically and $h_0$-invariant true $\check {\bf\Gamma}$ set, 
$X$ be a Polish space, and $A,B$ be disjoint analytic relations on $X$.\smallskip 

\noindent (1) Exactly one of the following holds:\smallskip  

(a) the set $A$ is separable from $B$ by a $\mbox{pot}({\bf\Gamma})$ set,\smallskip  

(b) there is $f\! :\! 2^\omega\!\rightarrow\! X$ continuous with  
$S_{\cal I}\!\subseteq\! (f\!\times\! f)^{-1}(A)$ and 
$\lceil T\rceil\!\setminus\! S_{\cal I}\!\subseteq\! (f\!\times\! f)^{-1}(B)$.\smallskip 

\noindent (2) If moreover $\bf\Gamma$ is of rank at least three, then we can have $f$ injective in (b).\smallskip 

\noindent (3) (Debs) We cannot replace $\lceil T\rceil\!\setminus\! S_{\cal I}$ with 
$\neg S_{\cal I}$ in (b).\end{cor} 

\noindent\bf Proof.\rm ~(1) We first prove the fact that Theorem \ref{compl} holds if $\cal I$ is only $h_0$-invariant, when $s\! =\!\emptyset$. The proof of Theorem \ref{compl} shows that $\cal I$ is weakly transferable if $\cal I$ is vertically and $h_0$-invariant. It remains to apply Lemma \ref{key} to $(u,v)\! :=\! (\emptyset ,\emptyset )$ and $G\! :=\! H$.\bigskip

 By Theorem \ref{compl}, (a) and (b) cannot hold simultaneously. Assume that $A$ is not separable from $B$ by a $\mbox{pot}({\bf\Gamma})$ set. This gives disjoint Borel subsets $C_0,C_1$ of $X$ such that $A\cap (C_0\!\times\! C_1)$ is not separable from 
$B\cap (C_0\!\times\! C_1)$ by a $\mbox{pot}({\bf\Gamma})$ set since the rank of $\bf\Gamma$ is at least two (consider a countable partition of the diagonal of $X$ into Borel rectangles with disjoint sides). We may assume that $C_0,C_1$ are clopen, refining the Polish topology if necessary. Theorem \ref{DL} gives, for each $\varepsilon\!\in\! 2$, 
$f_\varepsilon\! :\! 2^\omega\!\rightarrow\! C_\varepsilon$ continuous such that 
$S_{\cal I}\!\subseteq\! (f_0\!\times\! f_1)^{-1}\big( A\cap (C_0\!\times\! C_1)\big)$ and 
$\lceil T\rceil\!\setminus\! S_{\cal I}\!\subseteq\! 
(f_0\!\times\! f_1)^{-1}\big( B\cap (C_0\!\times\! C_1)\big)$. It remains to set 
$f(\alpha )\! :=\! f_\varepsilon (\alpha )$ if $\alpha\!\in\! N_\varepsilon$ since 
$\lceil T\rceil\!\subseteq\! N_0\!\times\! N_1$.\bigskip

\noindent (2) We apply Theorem \ref{DL} and the disjointness of $C_0$ and $C_1$.\bigskip

\noindent (3) See Theorem 1.13 in [L4].\hfill{$\square$}\bigskip

 We will construct some examples satisfying the assumptions of Theorem \ref{compl}.\bigskip

\noindent\bf Notation and definition.\rm ~We set $\mbox{FIN}\! :=\!\{\alpha\!\in\! 2^\omega\mid
\exists m\!\in\!\omega ~~\forall n\!\geq\! m~~\alpha (n)\! =\! 0\}$. Note that 
$\mathbb{E}_0\! =\! E_{\mbox{FIN}}$. We say that ${\cal I}\!\subseteq\! 2^\omega$ is \bf free\rm\ if 
${\cal I}\!\supseteq\!\mbox{FIN}$. 

\begin{prop} \label{transfer} Let ${\cal I}\!\subseteq\! 2^\omega$ be a free vertically invariant ideal. Then $\cal I$ is transferable.\end{prop}
 
\noindent\bf Proof.\rm ~Let $(n,\alpha ,F)$ be a transfer triple, and $\beta\!\in\! 2^\omega$. This gives an injection $i\! :\!\omega\!\rightarrow\!\omega$ such that $\big( i(m)\big)_0\! =\! (m)_0$ if 
$m\!\in\!\omega$. We set $N\! :=\!\{ m\!\in\!\omega\mid\beta (m)\! =\! 1\}$. As $\cal I$ is vertically invariant, $N\!\in\! {\cal I}$ is equivalent to $i[N]\!\in\! {\cal I}$. But $i[N]\! =\!\{ k\!\geq\! n\mid\big(\alpha\Delta F(\beta )\big)(k)\! =\! 1\}$. As $\cal I$ is a free ideal, $i[N]\!\in\! {\cal I}$ is equivalent to $\{ k\!\in\!\omega\mid\big(\alpha\Delta F(\beta )\big)(k)\! =\! 1\}\!\in\! {\cal I}$, so that 
$$\beta\!\in\! {\cal I}\ \Leftrightarrow\ N\!\in\! {\cal I}
\ \Leftrightarrow\ \{ k\!\in\!\omega\mid\big(\alpha\Delta F(\beta )\big)(k)\! =\! 1\}\!\in\! {\cal I}\ \Leftrightarrow\ \alpha\Delta F(\beta )\!\in\! {\cal I}.$$
This finishes the proof.\hfill{$\square$}\bigskip

\noindent\bf Notation.\rm\ We now introduce the operations that will be used to build our examples. They involve some bijection from $\omega^2$ onto 
$\omega$, which will not always be $<.,.>$. Indeed, in order to preserve the property of being vertically invariant, we will consider the bijection 
$\varphi\! :\!\omega^2\!\rightarrow\!\omega$ defined by
$$\varphi (n,p)\! :=\!\big<\! <n,(p)_0>,(p)_1\big>\mbox{,}$$
with inverse $q\!\mapsto\!\Big(\big( (q)_0\big)_0,<\big( (q)_0\big)_1,(q)_1>\!\Big)$.\bigskip

\noindent $\bullet$ Let $\alpha\!\in\! 2^\omega$ and $n\!\in\!\omega$. Recall that 
$(\alpha )_n\!\in\! 2^\omega$ is defined by $(\alpha )_n(p)\! :=\!\alpha (<n,p>)$. Similarly, we define 
${}^n(\alpha )\!\in\! 2^\omega$ by ${}^n(\alpha )(p)\! :=\!\alpha\big(\varphi (n,p)\big)$.\bigskip

\noindent $\bullet$ If $\alpha_0,...,\alpha_l\!\in\! 2^\omega$, then we define 
$\mbox{max}_{i\leq l}~\alpha_i\!\in\! 2^\omega$ by 
$(\mbox{max}_{i\leq l}~\alpha_i)(p)\! :=\!\mbox{max}_{i\leq l}~\alpha_i(p)$.

\bigskip

\noindent $\bullet$ If $\alpha ,\beta\!\in\! 2^\omega$, then we say that $\alpha\!\leq\!\beta$ when 
$\alpha (n)\!\leq\!\beta (n)$ for each $n\!\in\!\omega$. 

\begin{prop} \label{elem} Let $\alpha ,\beta ,\alpha_0,...,\alpha_l\!\in\! 2^\omega$ and $n\!\in\!\omega$. Then\smallskip

\noindent (1) $\alpha\!\leq\!\beta\Rightarrow (\alpha )_n\!\leq\! (\beta )_n$,\smallskip

\noindent (2) $(\mbox{max}_{i\leq l}~\alpha_i)_n\! =\!\mbox{max}_{i\leq l}~(\alpha_i)_n$,\smallskip

\noindent (3) $\alpha\!\in\! \mbox{FIN}\Rightarrow (\alpha )_n\!\in\! \mbox{FIN}$.\smallskip

These properties are also true with ${}^n(.)$ instead of $(.)_n$.\end{prop}
 
\noindent\bf Proof.\rm ~This is immediate.\hfill{$\square$}\bigskip

\noindent\bf Notation.\rm\ We now recall the operations of Lemma 1 in [Ca] (see also [Ka]). Let 
${\cal J},{\cal J}_0,{\cal J}_1,...\!\subseteq\! 2^\omega$.\medskip

\noindent $\bullet$ $\vec {\cal J}\! :=\! ({\cal J}_0,{\cal J}_1,...)$\smallskip

\noindent $\bullet$ 
$\vec {\cal J}^m\! :=\!\{\alpha\!\in\! 2^\omega\mid\forall n\!\in\!\omega ~~
{}^n(\alpha )\!\in\! {\cal J}_n\}$, and ${\cal J}^m\! :=\! ({\cal J},{\cal J},...)^m$\smallskip

\noindent $\bullet$ 
$\vec {\cal J}^a\! :=\!\{\alpha\!\in\! 2^\omega\mid\exists p\!\in\!\omega ~~\forall n\!\geq\! p~~
{}^n(\alpha )\!\in\! {\cal J}_n\}$, and ${\cal J}^a\! :=\! ({\cal J},{\cal J},...)^a$\bigskip

\noindent\bf Remark.\rm\ Proposition \ref{elem} implies that $\vec {\cal J}^m,\vec {\cal J}^a$ are ideals if the 
${\cal J}_n$'s are, free if the ${\cal J}_n$'s are.

\begin{lem} \label{twist} Let $n\!\in\!\omega$, ${\cal J}\!\subseteq\! 2^\omega$, and ${\cal I}\! :=\!\{\alpha\!\in\! 2^\omega\mid {}^n(\alpha )\!\in\! {\cal J}\}$. Then $\cal I$ is vertically invariant if $\cal J$ is.\end{lem}

\noindent\bf Proof.\rm ~Let $i\! :\!\omega\!\rightarrow\!\omega$ be injective such that 
$\big( i(m)\big)_0\! =\! (m)_0$ for each $m\!\in\!\omega$, and $N\!\subseteq\!\omega$ with characteristic function $\chi_N$. Then 
$$\begin{array}{ll}
N\!\in\! {\cal I}\!\!\!\!
& \Leftrightarrow\chi_N\!\in\! {\cal I}\Leftrightarrow {}^n(\chi_N)\!\in\! {\cal J}
\Leftrightarrow\{ p\!\in\!\omega\mid {}^n(\chi_N)(p)\! =\! 1\}\!\in\! {\cal J}\cr
& \Leftrightarrow\{ p\!\in\!\omega\mid\chi_N\big(\varphi (n,p)\big)\! =\! 1\}\!\in\! {\cal J}
\Leftrightarrow\{ p\!\in\!\omega\mid\varphi (n,p)\!\in\! N\}\!\in\! {\cal J}.
\end{array}$$
Similarly, $i[N]\!\in\! {\cal I}\Leftrightarrow\{ p\!\in\!\omega\mid\varphi (n,p)\!\in\! i[N]\}\!\in\! {\cal J}$. Recall that 
$\varphi (n,p)\! =\!\big<\! <n,(p)_0>,(p)_1\big>$. 
We define $I\! :\!\omega\!\rightarrow\!\omega$ by $I(p)\! :=\!\Big< (p)_0,\Big( i\big(\big<\! <n,(p)_0>,(p)_1\big>\big)\Big)_1\Big>$, so that 
$(p)_0\! =\!\big( I(p)\big)_0$ for each $p\!\in\!\omega$. Moreover, $I$ is injective. Indeed, $I(p)\! =\! I(p')$ implies successively that $(p)_0\! =\! (p')_0$, 
$$\begin{array}{ll}
i\big(\big<\! <n,(p)_0>,(p)_1\big>\big)\!\!\!\!
& =\!\Big<\Big( i\big(\big<\! <n,(p)_0>,(p)_1\big>\big)\Big)_0,\Big( i\big(\big<\! <n,(p)_0>,(p)_1\big>\big)\Big)_1\Big>\cr
& =\!\big<\! <n,(p)_0>,\big( I(p)\big)_1\big>\! =\!\big<\! <n,(p')_0>,\big( I(p')\big)_1\big>\cr
& =\ i\big(\big<\! <n,(p')_0>,(p')_1\big>\big)\mbox{,}
\end{array}$$
$\big<\! <n,(p)_0>,(p)_1\big>\! =\!\big<\! <n,(p')_0>,(p')_1\big>$, $(p)_1\! =\! (p')_1$ and $p\! =\! p'$. Now note that 
$$\begin{array}{ll}
\varphi\big( n,I(p)\big)\!\!\!\!\!  
& =\!\big<\! <n,\big( I(p)\big)_0>,\big( I(p)\big)_1\big>\! =\!\big<\! <n,(p)_0>,\big( I(p)\big)_1\big>\cr
& =\ i\big(\big<\! <n,(p)_0>,(p)_1\big>\big)\! =\ i\big( \varphi (n,p)\big) .
\end{array}$$
Thus
$$\begin{array}{ll}
\varphi (n,p)\!\in\! i[N]\!\!\!\!
& \Leftrightarrow\exists (n',p')\!\in\!\omega^2 ~~\varphi (n',p')\!\in\! N\ \wedge\ \varphi (n,p)\! =\! i\big(\varphi (n',p')\big)\cr
& \Leftrightarrow\exists (n',p')\!\in\!\omega^2 ~~\varphi (n',p')\!\in\! N\ \wedge\ \varphi (n,p)\! =\!\varphi\big( n',I(p')\big)\cr
& \Leftrightarrow\exists p'\!\in\!\omega ~~\varphi (n,p')\!\in\! N\ \wedge\ p\! =\! I(p')\cr
& \Leftrightarrow p\!\in\! I[\{ p'\!\in\!\omega\mid\varphi (n,p')\!\in\! N\} ].
\end{array}$$
Therefore $I[\{ p'\!\in\!\omega\mid\varphi (n,p')\!\in\! N\} ]\! =\!\{ p\!\in\!\omega\mid\varphi (n,p)\!\in\! i[N]\}$. As $\cal J$ is vertically invariant, 
$$\begin{array}{ll}
N\!\in\! {\cal I}\!\!\!\!
& \Leftrightarrow\{ p'\!\in\!\omega\mid\varphi (n,p')\!\in\! N\}\!\in\! {\cal J}\Leftrightarrow
I[\{ p'\!\in\!\omega\mid\varphi (n,p')\!\in\! N\} ]\!\in\! {\cal J}\cr
& \Leftrightarrow\{ p\!\in\!\omega\mid\varphi (n,p)\!\in\! i[N]\}\!\in\! {\cal J}\Leftrightarrow i[N]\!\in\! {\cal I}.
\end{array}$$
This finishes the proof.\hfill{$\square$}

\begin{cor} \label{transfervi} Let ${\cal J}_0,{\cal J}_1,...\!\subseteq\! 2^\omega$. Then $\vec {\cal J}^m,\vec {\cal J}^a$ are vertically invariant if the 
${\cal J}_n$'s are.\end{cor}

\noindent\bf Proof.\rm ~We set, for $n\!\in\!\omega$, 
${\cal I}_n\! :=\!\{\alpha\!\in\! 2^\omega\mid {}^n(\alpha )\!\in\! {\cal J}_n\}$, so that the ${\cal I}_n$'s are vertically invariant, by Lemma \ref{twist}. Let $i\! :\!\omega\!\rightarrow\!\omega$ be injective such that $\big( i(m)\big)_0\! =\! (m)_0$ for each $m\!\in\!\omega$, and $N\!\subseteq\!\omega$ with characteristic function $\chi_N$. Then 
$$\begin{array}{ll}
N\!\in\!\vec {\cal J}^m\!\!\!\!
& \Leftrightarrow\chi_N\!\in\!\vec {\cal J}^m
\Leftrightarrow\forall n\!\in\!\omega ~~{}^n(\chi_N)\!\in\! {\cal J}_n
\Leftrightarrow\forall n\!\in\!\omega ~~\chi_N\!\in\! {\cal I}_n\cr
& \Leftrightarrow\forall n\!\in\!\omega ~~N\!\in\! {\cal I}_n
\Leftrightarrow\forall n\!\in\!\omega ~~i[N]\!\in\! {\cal I}_n\Leftrightarrow i[N]\!\in\!\vec {\cal J}^m.
\end{array}$$
The proof is similar with $\vec {\cal J}^a$.\hfill{$\square$}\bigskip

 The next result is proved in [Ca] (see Lemmas 1 and 2).

\begin{lem} (Calbrix) \label{calbrix} Let ${\cal J}_0,{\cal J}_1,...\!\subseteq\! 2^\omega$ and $1\!\leq\!\xi\! <\!\omega_1$.\smallskip

\noindent (a) Assume that the ${\cal J}_n$'s are $\bormxi$-complete. Then $\vec {\cal J}^a$ is $\boraxp$-complete.\smallskip

\noindent (b) Assume that the ${\cal J}_n$'s are $\boraxi$-complete. Then $\vec {\cal J}^m$ is $\bormxp$-complete.\smallskip

\noindent (c) Assume that ${\cal J}_n$ is ${\bf\Sigma}^0_{2n+2}$-complete. Then $\vec {\cal J}^m$ is 
${\bf\Pi}^0_\omega$-complete.\smallskip

\noindent (d) Assume that $\lambda\! =\!\mbox{sup}_{n\in\omega}\uparrow\omega\! +\! 2\xi_n\! +\! 1$ is an infinite limit ordinal, and ${\cal J}_n$ is ${\bf\Sigma}^0_{\omega+2\xi_n+1}$-complete. Then 
$\vec {\cal J}^m$ is ${\bf\Pi}^0_\lambda$-complete.\end{lem}
 
 In the same spirit, we have the following.
 
\begin{lem} \label{transfercomplete} Let ${\cal J}_0,{\cal J}_1,...\!\subseteq\! 2^\omega$, and 
$\lambda\! =\!\mbox{sup}_{n\in\omega}\uparrow\xi_n$ be an infinite limit ordinal. We assume that 
${\cal J}_n$ is ${\bf\Pi}^0_{\xi_n}$-complete. Then $\vec {\cal J}^a$ is ${\bf\Sigma}^0_\lambda$-complete.\end{lem}

\noindent\bf Proof.\rm ~Assume that $A\! :=\!\bigcup_{n\in\omega}\uparrow ~A_n$, where 
$A_n\!\in\! {\bf\Pi}^0_{\xi_n}(2^\omega )$ (this is a typical ${\bf\Sigma}^0_\lambda$ set since 
$(\xi_n)_{n\in\omega}$ is strictly increasing). Let $f_n\! :\! 2^\omega\!\rightarrow\! 2^\omega$ continuous with $A_n\! =\! f_n^{-1}({\cal J}_n)$. We define $f\! :\! 2^\omega\!\rightarrow\! 2^\omega$ by 
$f(\alpha )(q)\! :=\! f_{((q)_0)_0}(\alpha )(<\big( (q)_0\big)_1,(q)_1>)$. Note that $f$ is continuous and 
${}^n\big( f(\alpha )\big)\! =\! f_n(\alpha )$ since
$${}^n\big( f(\alpha )\big)(p)\! =\! f(\alpha )\big(\varphi (n,p)\big)\! =\! f(\alpha )\big(\big<\! <n,(p)_0>,(p)_1\big>\big)\! =\! 
f_n(\alpha )(p).$$
Then 
$$\begin{array}{ll}
f(\alpha )\!\in\!\vec {\cal J}^a\!\!\!\!
& \Leftrightarrow\exists p\!\in\!\omega ~~\forall n\!\geq\! p~~{}^n\big( f(\alpha )\big)\!\in\! {\cal J}_n\Leftrightarrow\exists p\!\in\!\omega ~~\forall n\!\geq\! p~~f_n(\alpha )\!\in\! {\cal J}_n\cr
& \Leftrightarrow\exists p\!\in\!\omega ~~\forall n\!\geq\! p~~\alpha\!\in\! A_n
\Leftrightarrow\exists p\!\in\!\omega ~~\alpha\!\in\! A_p\Leftrightarrow\alpha\!\in\! A.
\end{array}$$
This finishes the proof.\hfill{$\square$}\bigskip

 We are now ready to introduce some examples.\bigskip

\noindent\bf Notation.\rm\ We set\smallskip

\noindent $\bullet$ $\mathbb{I}_3\! :=\!\{\alpha\!\in\! 2^\omega\mid\forall n\!\in\!\omega ~~(\alpha )_n\!\in\!\mbox{FIN}\}$,\smallskip

\noindent $\bullet$ $\mathbb{I}_{4+2n}\! :=\!\mathbb{I}_{3+2n}^a$ if $n\!\in\!\omega$,\smallskip
 
\noindent $\bullet$ $\mathbb{I}_{5+2n}\! :=\!\mathbb{I}_{4+2n}^m$ if $n\!\in\!\omega$,\smallskip
 
\noindent $\bullet$ $\mathbb{I}_\omega\! :=\! (\mathbb{I}_3,\mathbb{I}_5,...)^a$ and 
$\mathbb{J}_\omega\! :=\! (\mbox{FIN},\mathbb{I}_4,\mathbb{I}_6,...)^m$,\smallskip
 
\noindent $\bullet$ $\mathbb{I}_{\omega +2\xi +1}\! :=\!\mathbb{I}_{\omega +2\xi}^m$ and 
$\mathbb{J}_{\omega +2\xi +1}\! :=\!\mathbb{J}_{\omega +2\xi}^a$ if $\xi\! <\!\omega_1$,\smallskip
 
\noindent $\bullet$ $\mathbb{I}_{\omega +2\xi +2}\! :=\!\mathbb{I}_{\omega +2\xi +1}^a$ and 
$\mathbb{J}_{\omega +2\xi +2}\! :=\!\mathbb{J}_{\omega +2\xi +1}^m$ if $\xi\! <\!\omega_1$,\smallskip
 
\noindent $\bullet$ $\mathbb{I}_{\lambda}\! :=\! (\mathbb{I}_{\omega +2\xi_0+1},\mathbb{I}_{\omega +2\xi_1+1},...)^a$ and 
$\mathbb{J}_{\lambda}\! :=\! (\mathbb{J}_{\omega +2\xi_0+1},\mathbb{J}_{\omega +2\xi_1+1},...)^m$ if 
$$\lambda\! =\!\mbox{sup}_{n\in\omega}\uparrow\omega\! +\! 2\xi_n\! +\! 1$$ 
is an infinite limit countable ordinal.

\begin{cor} \label{firstex} All the sets previously defined are free and vertically invariant ideals, and in particular transferable. Moreover,\smallskip
 
\noindent $\bullet$ FIN is $\boratwo$-complete,\smallskip
 
\noindent $\bullet$ $\mathbb{I}_{2+2\xi +1}$ is ${\bf\Pi}^0_{2+2\xi +1}$-complete and $\mathbb{J}_{\omega +2\xi +1}$ is 
${\bf\Sigma}^0_{\omega +2\xi +1}$-complete,\smallskip
 
\noindent $\bullet$ $\mathbb{I}_{4+2\xi}$ is ${\bf\Sigma}^0_{4+2\xi}$-complete and $\mathbb{J}_{\omega +2\xi}$ is 
${\bf\Pi}^0_{\omega +2\xi}$-complete.\end{cor}

\noindent\bf Proof.\rm ~It is clear that\smallskip
 
\noindent - FIN and $\mathbb{I}_3$ are free ideals,\smallskip
 
\noindent - $\mbox{FIN}$ is vertically invariant and $\boratwo$-complete.\bigskip

\noindent $\bullet$ Let us prove that $\mathbb{I}_3$ is vertically invariant. Let 
$i\! :\!\omega\!\rightarrow\!\omega$ be injective such that $\big( i(m)\big)_0\! =\! (m)_0$ for each 
$m\!\in\!\omega$, and $N\!\subseteq\!\omega$ with characteristic function $\chi_N$. Then 
$$\begin{array}{ll}
N\!\in\!\mathbb{I}_3\!\!\!\!
& \Leftrightarrow\chi_N\!\in\!\mathbb{I}_3\Leftrightarrow\forall n\!\in\!\omega ~~(\chi_N)_n\!\in\!\mbox{FIN}\Leftrightarrow\forall n\!\in\!\omega ~~
\exists m\!\in\!\omega ~~\forall p\!\geq\! m~~(\chi_N)_n(p)\! =\! 0\cr
& \Leftrightarrow\forall n\!\in\!\omega ~~\exists m\!\in\!\omega ~~\forall p\!\geq\! m~<n,p>\ \notin\! N.
\end{array}$$
Thus $N\!\notin\!\mathbb{I}_3\Leftrightarrow\exists n\!\in\!\omega ~~\exists^\infty p\!\in\!\omega <n,p>\in\! N$. 
Similarly, 
$$\begin{array}{ll}
i[N]\!\notin\!\mathbb{I}_3\!\!\!\!
& \Leftrightarrow\exists n\!\in\!\omega ~~\exists^\infty p\!\in\!\omega ~<n,p>\in\! i[N]\cr
& \Leftrightarrow\exists n\!\in\!\omega ~~\exists^\infty p\!\in\!\omega ~~\exists (n',p')\!\in\!\omega^2~
<n',p'>\in\! N\mbox{ and }<n,p>=\! i(<n',p'>)\cr
& \Leftrightarrow\exists n\!\in\!\omega ~~\exists^\infty p\!\in\!\omega ~~\exists p'\!\in\!\omega ~<n,p'>\in\! N\mbox{ and }p\! =\!\big( i(<n,p'>)\big)_1\cr
& \Leftrightarrow\exists n\!\in\!\omega ~~\exists^\infty p'\!\in\!\omega ~<n,p'>\in\! N\cr
& \Leftrightarrow N\!\notin\!\mathbb{I}_3.
\end{array}$$
since $p'\!\mapsto\!\big( i(<n,p'>)\big)_1$ is injective because 
$\big( i(<n,p'>)\big)_1\! =\!\big( i(<n,p''>)\big)_1$ implies successively that 
$\big< n,\big( i(<n,p'>)\big)_1\big>\! =\!\big< n,\big( i(<n,p''>)\big)_1\big>$, $i(<n,p'>)\! =\! i(<n,p''>)$ and $p'\! =\! p''$.\bigskip

\noindent $\bullet$ $\mathbb{I}_3$ is $\bormthree$-complete by Lemma 1 in [Ca].\bigskip

\noindent $\bullet$ The rest follows from the remark before Lemma \ref{twist}, Corollary \ref{transfervi}, Proposition \ref{transfer}, and Lemmas \ref{calbrix} and \ref{transfercomplete}.\hfill{$\square$}\bigskip

 We now introduce some examples satisfying the assumptions of Theorem \ref{compl}.
 
\begin{lem} \label{exist}  Let $\bf\Gamma$ be a non self-dual Borel class of rank at least two. Then there is a vertically and $\mathbb{E}_0$-invariant true 
$\check {\bf\Gamma}$ set ${\cal I}\!\subseteq\! 2^\omega$ such that $S_{\cal I}$ and $S_{\neg {\cal I}}$ are dense in $\lceil T\rceil$. We can take ${\cal I}\! :=\!\mbox{FIN}$ if 
${\bf\Gamma}\! =\!\bormtwo$, and ${\cal I}\! :=\!\mathbb{I}_3\! :=\!
\{\gamma\!\in\! 2^\omega\mid\forall n\!\in\!\omega ~~(\gamma )_n\!\in\!\mbox{FIN}\}$ if 
${\bf\Gamma}\! =\!\borathree$.\end{lem}

\noindent\bf Proof.\rm ~If the rank of $\bf\Gamma$ is infinite or if ${\bf\Gamma}$ is in 
$\{ {\bf\Pi}^0_2,{\bf\Sigma}^0_3,{\bf\Pi}^0_4,{\bf\Sigma}^0_5,...\}$, then we apply Corollary \ref{firstex}, and in this case $\cal I$ can even be a free ideal, so that $E_{\cal I}$ is an equivalence relation. If ${\bf\Gamma}$ is in $\{ {\bf\Sigma}^0_2,{\bf\Pi}^0_3,{\bf\Sigma}^0_4,{\bf\Pi}^0_5,...\}$, then we take the complement of this ideal. It is also a vertically and $\mathbb{E}_0$-invariant true $\check {\bf\Gamma}$ set. It remains to see the density in $\lceil T\rceil$. So let $(u,v)\!\in\! T$. By Theorem \ref{compl}, $S_{\cal I}\cap (N_u\!\times\! N_v)$ is not pot$({\bf\Gamma})$ and 
$S_{\neg {\cal I}}\cap (N_u\!\times\! N_v)$ is not pot$(\check {\bf\Gamma})$, so that these sets are not empty.\hfill{$\square$}

\begin{prop} \label{exa} We can find a $D_2(\boraone )\!\subseteq\! F_\sigma$ acyclic graph $D$ on $2^\omega$ and Borel oriented subgraphs of $D$, without Borel countable coloring, of arbitrarily high potential complexity.\end{prop}

\noindent\bf Proof.\rm ~We set, for $\varepsilon\!\in\! 2$, 
$\psi_\varepsilon (\alpha )\! :=\!\varepsilon\alpha$, which defines homeomorphisms 
$\psi_\varepsilon\! :\! 2^\omega\!\rightarrow\! N_\varepsilon$. We set 
$D\! :=\! s\big( (\psi_0\!\times\!\psi_1)^{-1}(\lceil T\rceil )\big)\!\setminus\!\Delta (2^\omega )$, so that $D$ is a $D_2(\boraone )$ graph on $2^\omega$. Let us check that $D$ is acyclic. We argue by contradiction, which gives $n\!\geq\! 2$ and an injective $D$-path $(\gamma_i)_{i\leq n}$ with 
$(\gamma_0,\gamma_n)\!\in\! D$. This gives $(\varepsilon_i)_{i\leq n}$ such that 
$(\varepsilon_i\gamma_i,(1\! -\!\varepsilon_i)\gamma_{i+1})\!\in\! s(\lceil T\rceil )$ if $i\! <\! n$ and 
$(\varepsilon_n\gamma_0,(1\! -\!\varepsilon_n)\gamma_n)\!\in\! s(\lceil T\rceil )$. As 
$s(\lceil T\rceil )$ contains the couples of the form $(0\gamma ,1\gamma )$, this contradicts the acyclicity of $s(\lceil T\rceil )$.\bigskip

 Corollary \ref{firstex} gives a free vertically invariant ideal ${\cal I}\!\subseteq\! 2^\omega$ complete for a non self-dual Borel class $\bf\Gamma$ of arbitrarily high rank. Theorem \ref{compl} shows that $S_{\cal I}\!\notin\!\mbox{pot}(\check {\bf\Gamma})$. Note that the set 
$$G_{\cal I}\! :=\! (\psi_0\!\times\!\psi_1)^{-1}(S_{\cal I})$$ 
is Borel and not $\mbox{pot}(\check {\bf\Gamma})$. Thus $G_{\cal I}\!\setminus\!\Delta (2^\omega )\!\subseteq\! D\ \cap <_{\mbox{lex}}$ is a Borel oriented subgraph of $D$ and not $\mbox{pot}(\check {\bf\Gamma})$ in general. The freeness of $\cal I$ implies that there is no Borel countable coloring of 
$G_{\cal I}\!\setminus\!\Delta (2^\omega )$. This finishes the proof.\hfill{$\square$}

\section{$\!\!\!\!\!\!$ Some general facts}

$\underline{\mbox{\bf Antichains}}$\bigskip

 The following lemma gives a way of expanding antichains.

\begin{lem} \label{anti7} Let $\cal A$, $\cal B$ be $\sqsubseteq_c$-antichains made of nonempty subsets of $(N_0\!\times\! N_1)\cup (N_1\!\times\! N_0)$ such that each element $A$ of $\cal A$ has the property that $\big( 2^\omega ,A\cap (N_0\!\times\! N_1)\big)\not\sqsubseteq_c
\big( 2^\omega ,A\cap (N_1\!\times\! N_0)\big)$.\smallskip

(a) $\big\{ A^e\mid A\!\in\! {\cal B}\ \wedge\ e\!\in\!\{ =,\square ,\sqsubset\}\big\}$ is a 
$\sqsubseteq_c$-antichain.\smallskip 

(b) $\big\{ A^e\mid A\!\in\! {\cal A}\ \wedge\ e\!\in\!\{ =,\square ,\sqsubset ,\sqsupset\}\big\}$ is a 
$\sqsubseteq_c$-antichain.\smallskip

(c) If ${\cal A}\cup {\cal B}$ is a $\sqsubseteq_c$-antichain, then so is 
$$\big\{ A^e\mid A\!\in\! {\cal A}\ \wedge\ e\!\in\!\{ =,\square ,\sqsubset ,\sqsupset\}\big\}\cup
\big\{ A^e\mid A\!\in\! {\cal B}\ \wedge\ e\!\in\!\{ =,\square ,\sqsubset\}\big\} .$$\end{lem} 

\noindent\bf Proof.\rm ~(a) Let $A,B\!\in\! {\cal B}$ and $e,e'\!\in\!\{ =,\square ,\sqsubset\}$ such that $A^e\sqsubseteq_cB^{e'}$ with witness $f$. Then $f$ is also a witness for 
$A\sqsubseteq_cB$ since $A\! =\! A^e\!\setminus\!\Delta (2^\omega )$ and 
$B\! =\! B^{e'}\!\setminus\!\Delta (2^\omega )$. As $\cal B$ is an antichain, we must have 
$A\! =\! B$. Assume that $e\!\not=\! e'$. As $A^=$ is irreflexive and $A^\square$ is reflexive, 
$e'\! =\sqsubset$.\bigskip

 If $e$ is $=$, then pick $\big(\varepsilon\alpha ,(1\! -\!\varepsilon )\beta\big)\!\in\! A$. As $f$ is injective, $f(\varepsilon\alpha )\!\not=\! f\big( (1\! -\!\varepsilon )\beta\big)$, so that 
$\Big( f(\varepsilon\alpha ),f\big( (1\! -\!\varepsilon )\beta\big)\Big)$ is of the form 
$\big(\varepsilon'\gamma ,(1\! -\!\varepsilon')\delta\big)$. Assume for example that 
$\varepsilon'\! =\! 0$, the other case being similar. Then 
$\big( f(\varepsilon\alpha ),f(\varepsilon\alpha )\big)\!\in\! A^{e'}$, so that 
$(\varepsilon\alpha ,\varepsilon\alpha )\!\in\! A^e\! =\! A$, which is absurd.\bigskip

 If $e$ is not $=$, then it is $\square$. Here again, we pick 
$\big(\varepsilon\alpha ,(1\! -\!\varepsilon )\beta\big)$, and get $\varepsilon'$. Assume for example that $\varepsilon'\! =\! 1$, the other case being similar. Then 
$\big( f(\varepsilon\alpha ),f(\varepsilon\alpha )\big)\!\notin\! A^{e'}$, so that 
$(\varepsilon\alpha ,\varepsilon\alpha )\!\notin\! A^e\! =\! A\cup\Delta (2^\omega )$, which is absurd.\bigskip

\noindent (b) Let $A,B\!\in\! {\cal A}$ and $e,e'\!\in\!\{ =,\square ,\sqsubset ,\sqsupset\}$ such that $A^e\sqsubseteq_cB^{e'}$ with witness $f$. As in (a) we must have $A\! =\! B$, 
$e'\!\in\!\{\sqsubset ,\sqsupset\}$, and $e\!\in\!\{\sqsubset ,\sqsupset\}$ too. Assume that 
$e$ is $\sqsubset$ and $e'$ is $\sqsupset$, the other case being similar. Here again, we pick 
$\big(\varepsilon\alpha ,(1\! -\!\varepsilon )\beta\big)$, and get $\varepsilon'$. Assume for example that $\varepsilon'\! =\! 0$, the other case being similar. Then 
$\big( (1\! -\!\varepsilon )\beta ,(1\! -\!\varepsilon )\beta\big)\!\in\! A^e$, so that $\varepsilon\! =\! 1$. This shows that $\varepsilon\!\not=\!\varepsilon'$. Thus $A\cap (N_0\!\times\! N_1)$ is reducible to $A\cap (N_1\!\times\! N_0)$ with witness $f$, which contradicts our assumption.\bigskip

\noindent (c) Let $A,B\!\in\! {\cal A}\cup {\cal B}$ and 
$e,e'\!\in\!\{ =,\square ,\sqsubset ,\sqsupset\}$ such that $A^e\sqsubseteq_cB^{e'}$ with witness 
$f$. As in (a) we must have $A\! =\! B$. It remains to apply (a) and (b).\hfill{$\square$}

\begin{cor} \label{generalanti} Let $\bf\Gamma$ be a non self-dual Borel class of rank at least two, and $R$ be a true $\check {\bf\Gamma}$ relation on $2^\omega$, contained in 
$N_0\!\times\! N_1$, and such that $\overline{R}\!\setminus\! R$ is dense in $\overline{R}$. Then 
$$\big\{ A^e\mid A\!\in\!\{ R,R\cup\overline{R}^{-1},R\cup (\overline{R}^{-1}\!\setminus\! R^{-1})\}
\ \wedge\ e\!\in\!\{ =,\square ,\sqsubset ,\sqsupset\}\big\}\cup\big\{ s(R)^e\mid 
e\!\in\!\{ =,\square ,\sqsubset\}\big\}$$ 
is an antichain made of $\check {\bf\Gamma}\oplus {\bf\Gamma}$ sets.\end{cor} 

\noindent\bf Proof.\rm ~We set 
${\cal A}\! :=\!\{ R,R\cup\overline{R}^{-1},R\cup (\overline{R}^{-1}\!\setminus\! R^{-1})\}$ and 
${\cal B}\! :=\!\{ s(R)\}$. By Lemma \ref{anti7}.(c), it is enough to check that ${\cal A}\cup {\cal B}$ is an antichain.\bigskip

 Note the elements of $\cal A$ are not reducible to $s(R)$ since they are not symmetric. Similarly, the sets $R\cup\overline{R}^{-1},s(R)$ are not reducible to 
$R,R\cup (\overline{R}^{-1}\!\setminus\! R^{-1})$ since they are not antisymmetric.\bigskip
 
 If $A\! :=\! R\cup (\overline{R}^{-1}\!\setminus\! R^{-1})$ is reducible to $R$ with witness $f$, then $f$ is a homomorphism from $R$ into itself. Thus $f$ is a homomorphism from $\overline{R}$ into itself. Therefore $f$ is a homomorphism from $\overline{R}^{-1}$ into itself, which is absurd.\bigskip 
 
 As $s(R)$ is not closed and $s(R\cup\overline{R}^{-1})\! =\! s(A)\! =\! s(\overline{R})$ is, the sets $R,s(R)$ are not reducible to $R\cup\overline{R}^{-1}$, 
$A$.\bigskip

 If $A$ is reducible to 
$B\! :=\! R\cup\overline{R}^{-1}$ with witness $g$, then $g$ is a homomorphism from 
$\overline{R}\!\setminus\! R$ into itself since 
$\overline{R}\!\setminus\! R\! =\!\overline{B}\!\setminus\! B\!\subseteq\!\overline{A}\!\setminus\! A$. Thus $g$ is a homomorphism from $\overline{R}$ into itself, by our density assumption. Therefore $g$ reduces $R$ and $R^{-1}$ to themselves, which is absurd.\hfill{$\square$}\bigskip

 For ${\bf\Gamma}\! =\!\bormone$, a similar conclusion holds, for slightly different reasons. In this case, we set $R\! :=\!\mathbb{B}_0$, $\mathbb{N}_0\! :=\! R\cup\overline{R}^{-1}$ and 
$\mathbb{M}_0\! :=\! R\cup (\overline{R}^{-1}\!\setminus\! R^{-1})$. 

\begin{prop} \label{smallantipi01} The set  
$\big\{ A^e\mid A\!\in\!\{\mathbb{B}_0,\mathbb{N}_0,\mathbb{M}_0\}
\ \wedge\ e\!\in\!\{ =,\square ,\sqsubset ,\sqsupset\}\big\}\cup\big\{ s(\mathbb{B}_0)^e\mid 
e\!\in\!\{ =,\square ,\sqsubset\}\big\}$ is an antichain made of $D_2(\boraone )$ sets.\end{prop} 

\noindent\bf Proof.\rm ~We set ${\cal A}\! :=\!\{\mathbb{B}_0,\mathbb{N}_0,\mathbb{M}_0\}$ and 
${\cal B}\! :=\!\{ s(\mathbb{B}_0)\}$. By Lemma \ref{anti7}.(c), it is enough to check that 
${\cal A}\cup {\cal B}$ is an antichain. We argue as in the proof of Corollary \ref{generalanti}, except for the following.

\vfill\eject

 If $\mathbb{M}_0$ is reducible to $\mathbb{N}_0$ with witness $g$, then $g$ is a homomorphism from 
$$\overline{\mathbb{B}_0}\!\setminus\!\mathbb{B}_0\! =\!
\{ (0\alpha ,1\alpha )\mid\alpha\!\in\! 2^\omega\}$$ 
into itself again. This gives $k$ injective continuous such that 
$g(\varepsilon\alpha )\! =\!\varepsilon k(\alpha )$ if $\varepsilon\!\in\! 2$ and 
$\alpha\!\in\! 2^\omega$. Therefore $g$ reduces $\mathbb{B}_0$ and $\mathbb{B}_0^{-1}$ to themselves, which is absurd.\hfill{$\square$}\bigskip

 For ${\bf\Gamma}\! =\!\boraone$, we have a smaller antichain. In this case, we set 
$$R\! :=\!\{ (0\alpha ,1\alpha )\mid\alpha\!\in\! 2^\omega\}\! =\! 
R\cup (\overline{R}^{-1}\!\setminus\! R^{-1})\mbox{,}$$ 
so that $R\cup\overline{R}^{-1}\! =\! s(R)$. 

\begin{prop} \label{smallantisigma01} The set  
$$\big\{\{ (0\alpha ,1\alpha )\mid\alpha\!\in\! 2^\omega\}^e\mid 
e\!\in\!\{ =,\square ,\sqsubset ,\sqsupset\}\big\}\cup
\big\{ s(\{ (0\alpha ,1\alpha )\mid\alpha\!\in\! 2^\omega\} )^e\mid 
e\!\in\!\{ =,\square ,\sqsubset\}\big\}$$ 
is an antichain made of non-pot$(\boraone )$ closed sets.\end{prop} 

\noindent\bf Proof.\rm ~The intersection of the elements of our set with $N_0\!\times\! N_1$ is 
$\{ (0\alpha ,1\alpha )\mid\alpha\!\in\! 2^\omega\}$, which is not a countable union of Borel rectangles, and thus is not pot$(\boraone )$. So they are not pot$(\boraone )$. We set 
${\cal A}\! :=\!\big\{\{ (0\alpha ,1\alpha )\mid\alpha\!\in\! 2^\omega\}\big\}$ and 
${\cal B}\! :=\!\big\{ s(\{ (0\alpha ,1\alpha )\mid\alpha\!\in\! 2^\omega\} )\big\}$. By Lemma 
\ref{anti7}.(c), it is enough to check that ${\cal A}\cup {\cal B}$ is an antichain. But 
$\{ (0\alpha ,1\alpha )\mid\alpha\!\in\! 2^\omega\}$ is antisymmetric and 
$s(\{ (0\alpha ,1\alpha )\mid\alpha\!\in\! 2^\omega\} )$ is symmetric.\hfill{$\square$}\bigskip

\noindent $\underline{\mbox{\bf Minimality}}$\bigskip

 We are now interested in the minimality of $R$ and its associated relations among 
non-pot$({\bf\Gamma})$ relations when $R$ is not pot$({\bf\Gamma})$. Indeed, the intersection of the associated relations with $N_0\!\times\! N_1$ is exactly $R$, so that they are not 
pot$({\bf\Gamma})$ in this case. We start with a simple fact.

\begin{prop} \label{minimalitysymmetrization} Let $\bf\Gamma$ be a Borel class, and $R$ be a relation on $2^\omega$, which is $\sqsubseteq_c$-minimal among 
non-$\mbox{pot}({\bf\Gamma})$ relations. Then $s(R)$ is also $\sqsubseteq_c$-minimal among non-$\mbox{pot}({\bf\Gamma})$ relations if it is not $\mbox{pot}({\bf\Gamma})$.\end{prop} 

\noindent\bf Proof.\rm ~Assume that $(X,S)\sqsubseteq_c\big( 2^\omega ,s(R)\big)$ with witness $f$, where $X$ is Polish and $S$ is not $\mbox{pot}({\bf\Gamma})$. We set 
$B\! :=\! (f\!\times\! f)^{-1}(R)$, so that $(X,B)\sqsubseteq_c(2^\omega ,R)$, $S\! =\! B\cup B^{-1}$ and $B\!\notin\!\mbox{pot}({\bf\Gamma})$. By the minimality of $R$, 
$(2^\omega ,R)\sqsubseteq_c(X,B)$, and 
$\big( 2^\omega ,s(R)\big)\sqsubseteq_c\big( X,s(B)\big)\! =\! (X,S)$.\hfill{$\square$}\bigskip

 Similarly, the following holds.
 
\begin{lem} \label{symm} Let $\bf\Gamma$ be a Borel class. Assume that 
$O\!\subseteq\! N_0\!\times\! N_1$ is minimum  among non-pot$({\bf\Gamma})$ Borel subsets of a pot$(F_\sigma )$ Acyclic oriented graph. Then $s(O)$ is minimum among 
non-pot$({\bf\Gamma})$ Borel graphs contained in a pot$(F_\sigma )$ acyclic graph.\end{lem}

\noindent\bf Proof.\rm ~As $O\!\subseteq\! N_0\!\times\! N_1$ is not pot$({\bf\Gamma})$, $s(O)$ is not pot$({\bf\Gamma})$ too. Let $B$ be a non-pot$({\bf\Gamma})$ Borel graph on a Polish space $X$, contained in a pot$(F_\sigma )$ acyclic graph $G$, $C$ be a closed subset of 
$\omega^\omega$ and ${b\! :\! C\!\rightarrow\! X}$ be a continuous bijection.

\vfill\eject

 Note that $B\! =\!\big( B\cap (b\!\times\! b)[\leq_{\mbox{lex}}]\big)\cup
\big( B\cap (b\!\times\! b)[\geq_{\mbox{lex}}]\big)$. Then $B\cap (b\!\times\! b)[\leq_{\mbox{lex}}]$ or $B\cap (b\!\times\! b)[\geq_{\mbox{lex}}]$ is a non-pot$({\bf\Gamma})$ Borel oriented graph, both of them since 
${B\cap (b\!\times\! b)[\geq_{\mbox{lex}}]\! =\!\big( B\cap (b\!\times\! b)[\leq_{\mbox{lex}}]\big)^{-1}}$. Moreover, $B\cap (b\!\times\! b)[\leq_{\mbox{lex}}]$ is contained in the pot$(F_\sigma )$ Acyclic oriented graph $G\cap (b\!\times\! b)[\leq_{\mbox{lex}}]$. Therefore $(2^\omega ,O)$ is reducible to 
$\big( X,B\cap (b\!\times\! b)[\leq_{\mbox{lex}}]\big)$, which implies that 
$\big( 2^\omega ,s(O)\big)\sqsubseteq_c(X,B)$. This finishes the proof.\hfill{$\square$}

\begin{prop} \label{minimalitysquare} Let ${\bf\Gamma}\!\not=\!\boraone$ be a non self-dual Borel class, and $A$ be a digraph on $2^\omega$, $\sqsubseteq_c$-minimal among 
non-$\mbox{pot}({\bf\Gamma})$ relations. Then $A^\square$ is $\sqsubseteq_c$-minimal among 
non-$\mbox{pot}({\bf\Gamma})$ relations if it is not $\mbox{pot}({\bf\Gamma})$.\end{prop} 

\noindent\bf Proof.\rm ~Assume that $(X,S)\sqsubseteq_c(2^\omega ,A^\square )$ with witness 
$f$, where $X$ is Polish and $S$ is not $\mbox{pot}({\bf\Gamma})$. Then $S$ is reflexive and $f$ is also a witness for $\big( X,S\!\setminus\!\Delta (X)\big)\sqsubseteq_c(2^\omega ,A)$. As 
${\bf\Gamma}\!\supseteq\!\bormone$, $S\!\setminus\!\Delta (X)$ is not $\mbox{pot}({\bf\Gamma})$. By the minimality of $A$, $(2^\omega ,A)\sqsubseteq_c\big( X,S\!\setminus\!\Delta (X)\big)$, which implies that $\big( 2^\omega ,A^\square )\sqsubseteq_c(X,S)$. This finishes the proof.
\hfill{$\square$}\bigskip
 
 The reason why we exclude $\sqsupset$ for $s(R)$ is the following.
 
\begin{prop} \label{exclude} Let $R$ be a relation on $2^\omega$ contained in 
$N_0\!\times\! N_1$. We assume that\smallskip

\noindent (1) $(2^\omega ,R)\sqsubseteq_c(2^\omega ,R^{-1})$,\smallskip

\noindent (2) the projections of $\overline{R}$ are $N_0$ and $N_1$.\smallskip

 Then $\big( 2^\omega ,s(R)^\sqsubset\big)\sqsubseteq_c\big( 2^\omega ,s(R)^\sqsupset\big)$.\end{prop} 

\noindent\bf Proof.\rm ~Let $f$ be a witness for $(2^\omega ,R)\sqsubseteq_c(2^\omega ,R^{-1})$. Then $f$ reduces $s(R)$ to itself, and is a homomorphism from $\overline{R}$ into 
$\overline{R}^{-1}$. By (2), $f$ changes the first coordinate. Therefore $f$ reduces 
$s(R)^\sqsubset$ to $s(R)^\sqsupset$.\hfill{$\square$}

\begin{prop} \label{exclude2} Let $\bf\Gamma$ be a non self-dual Borel class, and $R$ be a relation on $2^\omega$ which is minimum among non-$\mbox{pot}({\bf\Gamma})$ Borel subsets of a closed Acyclic oriented graph. Then the reduction
$(2^\omega ,R)\sqsubseteq_c(2^\omega ,R^{-1})$ holds.\end{prop}

\noindent\bf Proof.\rm ~Note that $R^{-1}$ is a non-$\mbox{pot}({\bf\Gamma})$ Borel subset of a closed Acyclic oriented graph, which gives the result.\hfill{$\square$}\bigskip

\noindent $\underline{\mbox{\bf Examples and homomorphisms}}$\bigskip

 For $\bf\Gamma$ of rank at least two, the following is a key tool.

\begin{thm} \label{homo} Let ${\cal I}\!\subseteq\! 2^\omega$ be a vertically invariant set, and $F$ be a $F_\sigma$ relation on $2^\omega$ containing $\lceil T\rceil\cap\mathbb{E}_0$. 
\smallskip

(a) If $F$ is an Acyclic oriented graph, then there is an injective continuous homomorphism 
${f\! :\! 2^\omega\!\rightarrow\! 2^\omega}$ from 
$(\lceil T\rceil ,\neg\lceil T\rceil ,E_{\cal I},\neg E_{\cal I})$ into 
$(\lceil T\rceil ,\neg F,E_{\cal I},\neg E_{\cal I})$.\smallskip

(b) If $F$ is an acyclic graph, then there is an injective continuous homomorphism 
$f\! :\! 2^\omega\!\rightarrow\! 2^\omega$ from 
$(\lceil T\rceil ,\neg s(\lceil T\rceil ),E_{\cal I},\neg E_{\cal I})$ into 
$(\lceil T\rceil ,\neg F,E_{\cal I},\neg E_{\cal I})$ (and thus from 
$\big(\lceil T\rceil\cap E_{\cal I},\lceil T\rceil\!\setminus\! E_{\cal I},\neg s(\lceil T\rceil )\big)$ into 
$(\lceil T\rceil\cap E_{\cal I},\lceil T\rceil\!\setminus\! E_{\cal I},\neg F)$).\end{thm}

\noindent\bf Proof.\rm ~(a) By Lemmas \ref{suffacy} and \ref{meager}, $F$ is meager, which gives a decreasing sequence $(O_m)_{m\in\omega}$ of dense open subsets of $2^\omega$ such that 
$\neg F\! =\!\bigcap_{m\in\omega}~O_m$. We inductively construct $\delta\!\in\!\omega^\omega$, and define a function ${f\! :\! 2^\omega\!\rightarrow\! 2^\omega}$ by 
$f(\alpha )\! :=\!\alpha (0)0^{\delta (0)}\alpha (1)0^{\delta (1)}...$, so that $f$ will be injective continuous. The approximations $f_m\! :\! 2^m\!\rightarrow\! 2^{<\omega}$ of $f$ are defined by 
$f_m(s)\! :=\! s(0)0^{\delta (0)}...s(m\! -\! 1)0^{\delta (m-1)}$. We define $k_m$ by 
$\Sigma_{i<m}~\big( 1\! +\!\delta (i)\big)$, so that $f_m(s)\!\in\! 2^{k_m}$ for each $s\!\in\! 2^m$. We will build $\delta$ satisfying the following properties:
$$\begin{array}{ll}
& (1)~\big( f_m(u_m),f_m(v_m)\big)\!\in\! {\cal F}\mbox{, so that }
(f_m\!\times\! f_m)[T_m]\!\subseteq\! T_{k_m}\cr
& (2)~(k_m)_0\! =\! (m)_0\cr
& (3)~\forall (u,v)\!\in\! (2^m\!\times\! 2^m)\!\setminus\! T_m~~
N_{f_m(u)}\!\times\! N_{f_m(v)}\!\subseteq\! O_m
\end{array}$$
$\bullet$ Assume that this is done. If $(\alpha ,\beta )\!\in\!\lceil T\rceil$, then 
$(\alpha ,\beta )\vert m\!\in\! T_m$ for each $m\!\in\!\omega$, so that 
$$\big( f_m(\alpha\vert m),f_m(\beta\vert m)\big)\! =\!
\big( f(\alpha ),f(\beta )\big)\vert {k_m}\!\in\! T_{k_m}$$ 
for each $m\!\in\!\omega$ and $\big( f(\alpha ),f(\beta )\big)\!\in\!\lceil T\rceil$.\bigskip

 If $(\alpha ,\beta )\!\notin\!\lceil T\rceil$, then there is $m_0\!\in\!\omega$ such that 
$(\alpha ,\beta )\vert m\!\notin\! T_m$ if $m\!\geq\! m_0$. By Condition (4), 
$\big( f_m(\alpha\vert m),f_m(\beta\vert m)\big)\!\subseteq\!\big( f(\alpha ),f(\beta )\big)\!\in\! O_m$ if 
$m\!\geq\! m_0$, so that $\big( f(\alpha ),f(\beta )\big)\!\notin\! F$.\bigskip

 We define $i\! :\!\omega\!\rightarrow\!\omega$ by $i(m)\! :=\! k_m$. Note that $i$ is injective and 
$\big( i(m)\big)_0\! =\! (m)_0$ for each $m\!\in\!\omega$. Fix $\zeta\!\in\! 2^\omega$. We define 
$A\! :=\!\{ m\!\in\!\omega\mid\zeta (m)\! =\! 1\}$. Note that 
$i[A]\! =\!\{ p\!\in\!\omega\mid f(\zeta )(p)\! =\! 1\}$ since $k_m\!\in\! i[A]$ if and only if 
$\zeta (m)\! =\! 1$. As $\cal I$ is vertically invariant, $A\!\in\! {\cal I}$ is equivalent to 
$i[A]\!\in\! {\cal I}$. Thus $\zeta\!\in\! {\cal I}$ is equivalent to $f(\zeta )\!\in\! {\cal I}$. It remains to note that $f(\alpha\Delta\beta )\! =\! f(\alpha )\Delta f(\beta )$, and to apply the previous point to 
$\zeta\! :=\!\alpha\Delta\beta$, to see that $(\alpha ,\beta )\!\in\! E_{\cal I}$ if and only if 
$\big( f(\alpha ),f(\beta )\big)\!\in\! E_{\cal I}$.\bigskip

\noindent $\bullet$ So let us prove that the construction is possible. Note first that 
$$\big( f_0(u_0),f_0(v_0)\big)\! =\!\big( f_0(\emptyset ),f_0(\emptyset )\big)\! =\! 
(\emptyset ,\emptyset )\!\in\! {\cal F}\!\subseteq\! T$$ 
for any $\delta\!\in\!\omega^\omega$. Assume that $\delta (q)$ is constructed for $q\! <\! m$, which is the case for $m\! =\! 0$. If $(u,v)\!\in\! T_{m+1}$, then we can find $q\!\leq\! m$ and 
$w\!\in\! 2^{m-q}$ with $(u,v)\! =\! (u_q0w,v_q1w)$. In particular, 
$\big( f_q(u_q),f_q(v_q)\big)\!\in\! {\cal F}$ and 
$\big( f_{m+1}(u),f_{m+1}(v)\big)\vert (k_m\! +\! 1)$ is equal to 
$$\big(\ f_q(u_q)~0~0^{\delta (q)}w(0)0^{\delta (q+1)}...w(\vert w\vert\! -\! 1)\ ,\ 
f_q(v_q)~1~0^{\delta (q)}w(0)0^{\delta (q+1)}...w(\vert w\vert\! -\! 1)\ \big)\!\in\! T$$
since $q\! +\!\vert w\vert\! =\! m$. This implies that the map 
$\phi\! :\! s\!\mapsto\! f_{m+1}(s)\vert (k_m\! +\! 1)$ is an injective homomorphism of graphs from 
$\big( 2^{m+1},s(T_{m+1})\big)$ into $\big( 2^{<\omega},s(T)\big)$. As 
$\big( 2^{m+1},s(T_{m+1})\big)$ is acyclic connected and $\big( 2^{<\omega},s(T)\big)$ is acyclic, this map is an isomorphism onto its range by Lemma \ref{iso}. In particular, it preserves the lengths of the injective paths. If $(u,v)\!\in\! (2^{m+1}\!\times\! 2^{m+1})\!\setminus\! T_{m+1}$, then there are three cases:\bigskip

\noindent - $(v,u)\!\in\! T_{m+1}$, $\big(\phi(v),\phi(u)\big)\!\in\! T$, 
$\big(\phi(v)0^\infty ,\phi(u)0^\infty\big)\!\in\!\lceil T\rceil\cap\mathbb{E}_0\!\subseteq\! F$, 
$\big(\phi(u)0^\infty ,\phi(v)0^\infty\big)\!\notin\! F$ since $F$ is an oriented graph.\bigskip

\noindent - $u\! =\! v$, and $\big(\phi(u)0^\infty ,\phi(v)0^\infty\big)\!\notin\! F$ since $F$ is irreflexive.\bigskip

\noindent - $(u,v)\!\notin\! s(T_{m+1})\cup\Delta (2^{m+1})$, in which case the injective 
$s(T_{m+1})$-path from $u$ to $v$ has length at least 3. Thus the injective $s(T)$-path from 
$\phi (u)$ to $\phi(v)$ has length at least 3, and the injective $s(\lceil T\rceil\cap\mathbb{E}_0)$-path and the injective $s(F)$-path from $\phi (u)0^\infty$ to $\phi(v)0^\infty$ have length at least 3. Thus $\big(\phi(u)0^\infty ,\phi(v)0^\infty\big)$ is not in $F$ since $F$ is Acyclic.

\vfill\eject

 In every case, $\big(\phi(u)0^\infty ,\phi(v)0^\infty\big)\!\in\! O_{m+1}$. We choose 
$M\!\in\!\omega$ big enough so that 
$$N_{\phi(u)0^M}\!\times\! N_{\phi(v)0^M}\!\subseteq\! O_{m+1}$$ 
for each $(u,v)\!\in\! (2^{m+1}\!\times\! 2^{m+1})\!\setminus\! T_{m+1}$. There is $N\!\in\!\omega$ such that if $\delta (m)\! :=\! M\! +\! N$, then 
$$\big( f_{m+1}(u_{m+1}),f_{m+1}(v_{m+1})\big)\!\in\! {\cal F}$$ 
and $(k_{m+1})_0\! =\! (m\! +\! 1)_0$.\bigskip

\noindent (b) This is a consequence of the proof of (a).\hfill{$\square$}\bigskip

\bf\noindent Remark.\rm\ When $F$ is meager, we can replace the assumption ``$F$ is an Acyclic oriented graph" (resp., ``$F$ is an acyclic graph") with ``$F$ is an oriented graph (resp., a graph) and $F\cap\mathbb{E}_0\!\subseteq\! s(\lceil T\rceil\cap\mathbb{E}_0)$".\bigskip

 The version of Theorem \ref{homo} for ${\bf\Gamma}\! =\!\bormone$ is as follows.
 
\begin{thm} \label{homB0} Let $F$ be a closed relation on $2^\omega$ such that  
$\mathbb{B}_0\!\subseteq\! F\!\subseteq\! (N_0\!\times\! N_1)\cup (N_1\!\times\! N_0)$.\smallskip

(a) If $F$ is an Acyclic oriented graph, then there is an injective continuous homomorphism 
${f\! :\! 2^\omega\!\rightarrow\! 2^\omega}$ from 
$(\mathbb{B}_0,\overline{\mathbb{B}_0}\!\setminus\!\mathbb{B}_0,\neg\overline{\mathbb{B}_0})$ into $(\mathbb{B}_0,\overline{\mathbb{B}_0}\!\setminus\!\mathbb{B}_0,\neg F)$.\smallskip

(b) If $F$ is an acyclic graph, then there is an injective continuous homomorphism 
$f\! :\! 2^\omega\!\rightarrow\! 2^\omega$ from $\big(\mathbb{B}_0,
\overline{\mathbb{B}_0}\!\setminus\!\mathbb{B}_0,\neg s(\overline{\mathbb{B}_0})\big)$ into 
$(\mathbb{B}_0,\overline{\mathbb{B}_0}\!\setminus\!\mathbb{B}_0,\neg F)$.\end{thm}

\noindent\bf Proof.\rm ~(a) The proof is quite similar to that of Theorem \ref{Fsigma}. Note that $\neg F$ is a dense open set, by Lemmas \ref{suffacy} and \ref{meager}. We define $\psi_0\! :=\! {h_0}_{\vert N_0}$ and $\psi_{n+1}\! :\! N_{0s_n0}\!\rightarrow\! N_{1s_n1}$ by 
$\psi_{n+1}(0s_n0\gamma )\! :=\! 1s_n1\gamma$, so that 
$\overline{\mathbb{B}_0}\! =\!\bigcup_{n\in\omega}~\mbox{Gr}(\psi_n)$.\bigskip

\noindent $\bullet$ We construct $\Psi\! :\! 2^{<\omega}\!\rightarrow\! 2^{<\omega}$ and 
$\delta\!\in\!\omega^\omega$ strictly increasing satisfying the following conditions:
$$\begin{array}{ll}
& (1)~\forall s\!\in\! 2^{<\omega}~~\forall\varepsilon\!\in\! 2~~
\Psi (s)\!\subsetneqq\!\Psi (s\varepsilon )\cr
& (2)~\forall l\!\in\!\omega ~~\exists k_l\!\in\!\omega ~~\forall s\!\in\! 2^l~~\vert\Psi (s)\vert\! =\! k_l\cr
& (3)~\delta (0)\! =\! 0\ \wedge\ \forall v\!\in\! 2^{<\omega}~~\exists w\!\in\! 2^{<\omega}~~
\big(\Psi (0v),\Psi (1v)\big)\! =\! (0w,1w)\cr
& (4)~\forall n\!\in\!\omega ~~\forall v\!\in\! 2^{<\omega}~~\exists w\!\in\! 2^{<\omega}~~
\big(\Psi (0s_n0v),\Psi (1s_n1v)\big)\! =\! (0s_{\delta (n+1)-1}0w,1s_{\delta (n+1)-1}1w)\cr
& (5)~\forall (s,t)\!\in\! (2\!\times\! 2)^{<\omega}~~
(N_s\!\times\! N_t)\cap\overline{\mathbb{B}_0}\! =\!\emptyset\Rightarrow 
N_{\Psi (s)}\!\times\! N_{\Psi (t)}\!\subseteq\!\neg F
\end{array}$$
$\bullet$ Assume that this is done. We define $f\! :\! 2^\omega\!\rightarrow\! 2^\omega$ by 
$\{ f(\alpha )\}\! =\!\bigcap_{n\in\omega}~N_{\Psi (\alpha\vert n)}$, and $f$ is continuous. Condition (4) ensures that $\mathbb{B}_0\!\subseteq\! (f\!\times\! f)^{-1}(\mathbb{B}_0)$, and Condition (5) ensures that $\neg\overline{\mathbb{B}_0}\!\subseteq\! (f\!\times\! f)^{-1}(\neg F)$. Note that 
$\overline{\mathbb{B}_0}\!\setminus\!\mathbb{B}_0\! =\!
\{ (0\gamma ,1\gamma )\mid\gamma\!\in\! 2^\omega\}\! =\!\mbox{Gr}(\psi_0)$. Condition (3) ensures that $\overline{\mathbb{B}_0}\!\setminus\!\mathbb{B}_0\!\subseteq\! 
(f\!\times\! f)^{-1}(\overline{\mathbb{B}_0}\!\setminus\!\mathbb{B}_0)$. In order to see that $f$ is injective, it is enough to check that $\Psi (s0)\!\not=\!\Psi (s1)$ if $s\!\in\! 2^{<\omega}$, and we may assume that $s\!\not=\!\emptyset$.\bigskip

 We set, for $l\!\in\!\omega$, $B_l\! :=\!\{ (s,t)\!\in\! 2^l\!\times\! 2^l\mid (N_s\!\times\! N_t)\cap 
s(\overline{\mathbb{B}_0})\!\not=\!\emptyset\}$. Note that $(2^l,B_l)$ is a connected acyclic graph if $l\!\geq\! 1$, by induction on $l$. Indeed, 
$B_1\! =\!\{ (0,1),(1,0)\}$ and 
$$B_{l+1}\! =\!\{ (s\varepsilon ,t\varepsilon )\mid (s,t)\!\in\! B_l\ \wedge\ \varepsilon\!\in\! 2\}\cup
\{ (0s_{l-1}0,1s_{l-1}1),(1s_{l-1}1,0s_{l-1}0)\}$$
if $l\!\geq\! 1$. As $(2^l,B_l)$ is isomorphic to $(2^l\!\times\!\{\varepsilon\} ,B_{l+1})$, 
$(2^{l+1},B_{l+1})$ is a connected acyclic graph.

\vfill\eject

 If $(s,t)\!\in\! (2\!\times\! 2)^{<\omega}$, then $q^{s,t}\! :=\! (v^{s,t}_i)_{i\leq L_{s,t}}$ is the unique injective $B_{\vert s\vert}$-path from $s$ to $t$. Assume that $s\!\in\! 2^l$. We fix, for each 
$i\! <\! L\! :=\! L_{s,0s_{l-1}}$, $n_i\! :=\! n^{s,0s_{l-1}}_i\!\in\!\omega$ and 
$\varepsilon_i\! :=\!\varepsilon^{s,0s_{l-1}}_i\!\in\! 2$ such that 
$v_{i+1}^{s,0s_{l-1}}0^\infty\! =\!\psi_{n_i}^{\varepsilon_i}(v_i^{s,0s_{l-1}}0^\infty )$, so that 
$$\Psi (s1)0^\infty\! =\!\psi_{\delta (n_0)}^{-\varepsilon_0}...
\psi_{\delta (n_{L-1})}^{-\varepsilon_{L-1}}\psi_0\psi^{}_{\delta (l)}
\psi_{\delta (n_{L-1})}^{\varepsilon_{L-1}}...\psi_{\delta (n_0)}^{\varepsilon_0}
\big(\Psi (s0)0^\infty\big) .$$ 
As $k_{l+1}\! >\!\delta (l)\! >\!\mbox{sup}_{i<L}~\delta (n_i)$, $\Psi (s0)\!\not=\!\Psi (s1)$.\bigskip

\noindent $\bullet$ It remains to prove that the construction is possible. We first set 
$\Psi (\emptyset )\! :=\!\emptyset$. As $F$ is a closed oriented graph and 
$(0^\infty ,10^\infty )\!\in\!\overline{\mathbb{B}_0}$, $(10^\infty ,0^\infty )\!\notin\! F$. This gives 
$N\!\in\!\omega$ such that $N_{10^N}\!\times\! N_{0^{N+1}}\!\subseteq\!\neg F$, and we set 
$\Psi (\varepsilon )\! :=\!\varepsilon 0^N$. Assume that $\Psi [2^{\leq l}]$ and 
$\big(\delta (j)\big)_{j<l}$ satisfying (1)-(5) have been constructed, which is the case for 
$l\!\leq\! 1$. Let $l\!\geq\! 1$. Note that $\Psi_{\vert 2^l}$ is an injective homomorphism from 
$s(B_l)$ into $s(B_{k_l})$, and therefore an isomorphism of graphs onto its range by Lemma 
\ref{iso}. Moreover, $\delta (n\! +\! 1)\! <\! k_l$ if $n\! <\! l\! -\! 1$. Let 
$\delta (l)\! >\!\mbox{sup}_{n<l-1}~\delta (n\! +\! 1)$ such that 
$\Psi (0s_{l-1})\! -\! 0\!\subseteq\! s_{\delta (l)-1}$. We define temporary versions 
$\tilde\Psi (u\varepsilon )$ of the $\Psi (u\varepsilon )$'s by $\tilde\Psi (u\varepsilon )\! :=\!
\Psi (u)\big( s_{\delta (l)-1}\varepsilon\! -\! s_{\delta (l)-1}\vert (k_l\! -\! 1)\big)$, ensuring Conditions (1)-(4).\bigskip

 For Condition (5), if $s(0)\! =\! t(0)$, then $N_{\Psi (s)}\!\times\! N_{\Psi (t)}$ will be a subset of 
$N_0^2\cup N_1^2\!\subseteq\!\neg F$. If 
$(N_t\!\times\! N_s)\cap\overline{\mathbb{B}_0}\!\not=\!\emptyset$, then 
$\big(\tilde\Psi (t)0^\infty ,\tilde\Psi (s)0^\infty\big)\!\in\!\overline{\mathbb{B}_0}\!\subseteq\! F$ and 
$\big(\tilde\Psi (s)0^\infty ,\tilde\Psi (t)0^\infty\big)\!\notin\! F$. This gives $M\!\in\!\omega$ such that $N_{\tilde\Psi (s)0^M}\!\times\! N_{\tilde\Psi (t)0^M}\!\subseteq\!\neg F$, and we set 
$\Psi'(u\varepsilon )\! :=\!\tilde\Psi (u\varepsilon )0^M$.\bigskip

 So we may assume that $L\! :=\! L^{s,t}\!\geq\! 2$. Here again, $\tilde\Psi_{\vert 2^{l+1}}$ is an isomorphism of graphs onto its range. This implies that $\big(\tilde\Psi (v^{s,t}_i)\big)_{i\leq L}$ is the injective $s(B_{\vert\tilde\Psi (s)\vert})$-path from $\tilde\Psi (s)$ to $\tilde\Psi (t)$. Thus 
$\big(\tilde\Psi (v^{s,t}_i)0^\infty\big)_{i\leq L}$ is the injective $s(\overline{\mathbb{B}_0})$-path (and also $s(F)$-path) from $\tilde\Psi (s)0^\infty$ to $\tilde\Psi (t)0^\infty$. Thus 
$\big(\tilde\Psi (s)0^\infty ,\tilde\Psi (t)0^\infty\big)\!\notin\! F$ since $L\!\geq\! 2$. We conclude as in the previous case.\bigskip

\noindent (b) This is a consequence of the proof of (a) (here, $\Psi (\varepsilon )\! :=\!\varepsilon$).
\hfill{$\square$}\bigskip

\noindent\bf Remark.\rm ~This proof shows that we can replace the assumption ``$F$ is closed" with ``$F$ is $F_\sigma$ and the disjoint union $s(F)\cup\mbox{Gr}(h_0)$ is acyclic". In the proof, we write $\neg F\! =\!\bigcap_{l\in\omega}~O_l$, where $O_l$ is dense open, and replace 
$\neg F$ with $O_{\vert s\vert}$ in (5).\bigskip

 For ${\bf\Gamma}\! =\!\bormone$, the following holds.

\begin{lem} \label{propertiespi01} The set $R\! :=\!\mathbb{B}_0$ is a $D_2(\boraone )$ relation on $2^\omega$, contained in $N_0\!\times\! N_1$, satisfying the following properties.\smallskip

\noindent (1) For each $s\!\in\! 2^{<\omega}$, and for each dense $G_\delta$ subset $C$ of 
$2^\omega$, $R\cap\big( 2\!\times\! (N_s\cap C)\big)^2$ is not $\mbox{pot}(\bormone )$.\smallskip

\noindent (2) $\overline{R}$ is Acyclic.\smallskip

\noindent (3) The projections of $\overline{R}$ are $N_0$ and $N_1$.\end{lem}

\noindent\bf Proof.\rm ~(1) As the maps $f\! :\!\alpha\!\mapsto\! 0\alpha$ and $g\! :\!\beta\!\mapsto\! 1\beta$ satisfy 
$$\mathbb{G}_0\cap (N_s\cap C)^2\! =\! (f\!\times\! g)^{-1}\Big( R\cap\big( 2\!\times\! (N_s\cap C)\big)^2\Big)\mbox{,}$$ 
it is enough to see that $\mathbb{G}_0\cap (N_s\cap C)^2\!\notin\!\mbox{pot}(\bormone )$. We argue by contradiction, which gives a countable partition of $N_s\cap C$ into Borel sets whose square does not meet $\mathbb{G}_0$. One of these Borel sets has to be non-meager, which is absurd, as in the proof of Proposition 6.2 in [K-S-T].

\vfill\eject

\noindent (2) The map $\varepsilon\alpha\!\mapsto\! (\varepsilon ,\alpha )$ is an isomorphism from 
$s(\overline{\mathbb{B}_0})$ onto $s(G_{\overline{\mathbb{G}_0}})$, which is acyclic by Proposition \ref{conn} and Lemma \ref{suffacy}.\bigskip

\noindent (3) Note that $\{ (0\alpha ,1\alpha )\mid\alpha\!\in\! 2^\omega\}\!\subseteq\!\overline{\mathbb{B}_0}\!\subseteq\! N_0\!\times\! N_1$.\hfill{$\square$}\bigskip

\noindent $\underline{\mbox{\bf Basis}}$\bigskip

 We first introduce a definition generalizing the conclusion of Corollary \ref{corDL}. In order to make it work for the first Borel classes, we add an acyclicity assumption.

\begin{defi} Let ${\cal I}\!\subseteq\! 2^\omega$, and ${\bf\Gamma},{\bf\Gamma}'$ be classes of Borel sets closed under continuous pre-images. We say that $S_{\cal I}$ has the \bf 
$({\bf\Gamma},{\bf\Gamma}')$-basis property\ \it if for each Polish space $X$, and for each pair $A,B$ of disjoint analytic relations on $X$ such that $A$ is contained in a 
$\mbox{pot}({\bf\Gamma}')$ symmetric acyclic relation, exactly one of the following holds:\smallskip

(a) the set $A$ is separable from $B$ by a $\mbox{pot}({\bf\Gamma})$ set,\smallskip

(b) there is $g\! :\! 2^\omega\!\rightarrow\! X$ injective continuous such that 
$S_{\cal I}\!\subseteq\! (g\!\times\! g)^{-1}(A)$ and 
$\lceil T\rceil\!\setminus\! S_{\cal I}\!\subseteq\! (g\!\times\! g)^{-1}(B)$.\end{defi}

 Corollary \ref{corDL} says that if $\bf\Gamma$ is a non self-dual Borel class of rank at least three and 
$\cal I$ is a vertically and $h_0$-invariant true $\check {\bf\Gamma}$ set, then $S_{\cal I}$ has the 
$({\bf\Gamma},{\bf\Gamma}')$-basis property for each class of Borel sets ${\bf\Gamma}'$ closed under continuous pre-images.

\begin{thm} \label{fourpi0>1} Let $\bf\Gamma$ be a non self-dual Borel class of rank at least two, ${\cal I}\!\subseteq\! 2^\omega$ be a vertically and $\mathbb{E}_0$-invariant true 
$\check {\bf\Gamma}$ set such that $R\! :=\! S_{\cal I}$ is dense in $\lceil T\rceil$ 
(${\cal I}\! =\!\mbox{FIN}$ if ${\bf\Gamma}\! =\!\bormtwo$), and 
$\bormone\!\subseteq\! {\bf\Gamma}'\!\subseteq\! F_\sigma$ be a class of Borel sets closed under continuous pre-images. We assume that $R$ has the 
$({\bf\Gamma},{\bf\Gamma}')$-basis property. Then 
$\big\{ R,R\cup\overline{R}^{-1},R\cup (\overline{R}^{-1}\!\setminus\! R^{-1}),s(R)\big\}$ is a basis for the class of non-$\mbox{pot}({\bf\Gamma})$ Borel subsets of a $\mbox{pot}({\bf\Gamma}')$ acyclic graph.\end{thm}

\noindent\bf Proof.\rm\ By Theorem \ref{compl}, all the examples are in the context of the theorem. So let $B$ be a non-pot$({\bf\Gamma})$ Borel relation on a Polish space $X$, contained in a 
pot$({\bf\Gamma}')$ acyclic graph $H$. We can change the Polish topology and assume that $H$ is in ${\bf\Gamma}'$. We set $G\! :=\! B\cap B^{-1}$.\bigskip

\noindent\bf Case 1\rm\ $G$ is $\mbox{pot}({\bf\Gamma})$.\bigskip

 Assume first that ${\bf\Gamma}\!\not=\!\bormtwo$. Note that $B\!\setminus\! G$ is not separable from $H\!\setminus\! B$ by a $\mbox{pot}({\bf\Gamma})$ set $P$, since otherwise 
$B\! =\! (P\cap H)\cup G\!\in\!\mbox{pot}({\bf\Gamma})$. As $R$ has the 
$({\bf\Gamma},{\bf\Gamma}')$-basis property, there is $g\! :\! 2^\omega\!\rightarrow\! X$ injective continuous such that $R\!\subseteq\! (g\!\times\! g)^{-1}(B\!\setminus\! G)$ and 
$\overline{R}\!\setminus\! R\!\subseteq\! (g\!\times\! g)^{-1}(H\!\setminus\! B)$. Theorem \ref{homo} gives an injective continuous homomorphism $h\! :\! 2^\omega\!\rightarrow\! 2^\omega$ from 
$(\lceil T\rceil ,\neg s(\lceil T\rceil ),E_{\cal I},\neg E_{\cal I})$ into 
$(\lceil T\rceil ,\neg (g\!\times\! g)^{-1}(H),E_{\cal I},\neg E_{\cal I})$. We set $k\! :=\! g\circ h$ and $B'\! :=\! (k\!\times\! k)^{-1}(B)$, so that $(2^\omega ,B')\sqsubseteq_c(X,B)$ and 
$R\!\subseteq\! B'\!\subseteq\! R\cup (\overline{R}^{-1}\!\setminus\! R^{-1})$. Indeed, $h$ is a homomorphism from $R^{-1}$ into itself, and $g$ is a homomorphism from $R^{-1}$ into $\neg B$, since otherwise there is $(\alpha ,\beta )\!\in\! R^{-1}$ with $\big( g(\alpha ),g(\beta )\big)\!\in\! B$, and $\big( g(\beta ),g(\alpha )\big)\!\in\! G\!\setminus\! G$. If ${\bf\Gamma}\! =\!\bormtwo$, then we argue similarly: $B\!\setminus\! G$ is not separable from $\neg B$ by a $\mbox{pot}({\bf\Gamma})$ set, and we can apply Theorem \ref{homo} since $(g\!\times\! g)^{-1}(H)$ contains 
$R\! =\!\lceil T\rceil\cap\mathbb{E}_0$. So we may assume that 
$R\!\subseteq\! B\!\subseteq\! R\cup (\overline{R}^{-1}\!\setminus\! R^{-1})$ and $X\! =\! 2^\omega$. We write $B\! =\! R\cup S$, where $S$ is a Borel subset of $\overline{R}^{-1}\!\setminus\! R^{-1}$.

\vfill\eject

\noindent\bf Case 1.1\rm\ $R$ is not separable from $S^{-1}$ by a 
$\mbox{pot}({\bf\Gamma})$ set.\bigskip

 As $R$ has the $({\bf\Gamma},{\bf\Gamma}')$-basis property, we can find  
$g'\! :\! 2^\omega\!\rightarrow\! 2^\omega$ injective continuous such that 
$R\!\subseteq\! (g'\!\times\! g')^{-1}(R)$ and 
$\overline{R}\!\setminus\! R\!\subseteq\! (g'\!\times\! g')^{-1}(S^{-1})\!\subseteq\! 
(g'\!\times\! g')^{-1}(\neg B)$. Note that $(g'\!\times\! g')^{-1}\big( s(\overline{R})\big)$ is a closed acyclic graph containing $\lceil T\rceil$. Theorem \ref{homo} gives an injective continuous homomorphism $h'\! :\! 2^\omega\!\rightarrow\! 2^\omega$ from 
$(\lceil T\rceil ,\neg s(\lceil T\rceil ),E_{\cal I},\neg E_{\cal I})$ into 
$(\lceil T\rceil ,\neg (g'\!\times\! g')^{-1}\big( s(\overline{R})\big) ,E_{\cal I},\neg E_{\cal I})$. The map $k'\! :=\! g'\circ h'$ reduces $R\cup (\overline{R}^{-1}\!\setminus\! R^{-1})$ to $B$.\bigskip

\noindent\bf Case 1.2\rm\ $R$ is separable from $S^{-1}$ by a 
$\mbox{pot}({\bf\Gamma})$ set.\bigskip
 
 Let $Q\!\subseteq\!\lceil T\rceil\!\subseteq\! N_0\!\times\! N_1$ be such a set. Note that 
$R$ is not separable from $Q\!\setminus\! R$ by a $\mbox{pot}({\bf\Gamma})$ set, by Theorem 
\ref{compl}. As $R$ has the $({\bf\Gamma},{\bf\Gamma}')$-basis property, there is  
$g''\! :\! 2^\omega\!\rightarrow\! 2^\omega$ injective continuous with 
$R\!\subseteq\! (g''\!\times\! g'')^{-1}(R)$ and 
$\overline{R}\!\setminus\! R\!\subseteq\! (g''\!\times\! g'')^{-1}(Q\!\setminus\! R)$. Note that $g''$ reduces $R$ to $B$ on $s(\overline{R})$. Theorem \ref{homo} gives an injective continuous homomorphism $l''\! :\! 2^\omega\!\rightarrow\! 2^\omega$ from 
$(\lceil T\rceil ,\neg s(\lceil T\rceil ),E_{\cal I},\neg E_{\cal I})$ into 
$(\lceil T\rceil ,\neg (g''\!\times\! g'')^{-1}\big( s(\overline{R})\big) ,E_{\cal I},\neg E_{\cal I})$. Note that $g''\circ l''$ reduces $R$ to $B$.\bigskip

\noindent\bf Case 2\rm\ $G$ is not $\mbox{pot}({\bf\Gamma})$.\bigskip

 Assume first that ${\bf\Gamma}\!\not=\!\bormtwo$. Note that $G$ is not separable from $H\!\setminus\! B$ by a $\mbox{pot}({\bf\Gamma})$ set $P$, since otherwise 
${G\! =\! (P\cap H)\cap (P\cap H)^{-1}}$ 
would be $\mbox{pot}({\bf\Gamma})$. As in Case 1 we get $g\! :\! 2^\omega\!\rightarrow\! X$ injective continuous such that $R\!\subseteq\! (g\!\times\! g)^{-1}(G)$ and 
$\overline{R}\!\setminus\! R\!\subseteq\! (g\!\times\! g)^{-1}(H\!\setminus\! B)$. Theorem \ref{homo} gives an injective continuous homomorphism ${h\! :\! 2^\omega\!\rightarrow\! 2^\omega}$ from 
$(\lceil T\rceil ,\neg s(\lceil T\rceil ),E_{\cal I},\neg E_{\cal I})$ 
into $(\lceil T\rceil ,\neg (g\!\times\! g)^{-1}(H),E_{\cal I},\neg E_{\cal I})$. We set 
$k\! :=\! g\circ h$ and also ${B'\! :=\! (k\!\times\! k)^{-1}(B)}$, so that 
$s(R)\!\subseteq\! B'\!\subseteq\! R\cup\overline{R}^{-1}$ 
and $(2^\omega ,B')\sqsubseteq_c(X,B)$. Indeed, 
$h$ is a homomorphism from $\overline{R}^{-1}\cap E_{\cal I}$ into itself, and $g$ is a homomorphism from $\overline{R}^{-1}\cap E_{\cal I}$ into $B$, since $G$ is symmetric. If 
${\bf\Gamma}\! =\!\bormtwo$, then we argue similarly: $G$ is not separable from 
$\neg B$ by a $\mbox{pot}({\bf\Gamma})$ set, and we can apply Theorem \ref{homo} since 
$(g\!\times\! g)^{-1}(H)$ contains $R\! =\!\lceil T\rceil\cap\mathbb{E}_0$. So we may assume that 
$X\! =\! 2^\omega$ and $s(R)\!\subseteq\! B\!\subseteq\! R\cup\overline{R}^{-1}$. We write 
$B\! =\! s(R)\cup S$, where $S$ is a Borel subset of $\overline{R}^{-1}\!\setminus\! R^{-1}$.\bigskip

\noindent\bf Case 2.1\rm\ $R$ is not separable from $S^{-1}$ by a 
$\mbox{pot}({\bf\Gamma})$ set.\bigskip

 We argue as in Case 1.1 to see that $\big( 2^\omega ,R\cup\overline{R}^{-1}\big)\sqsubseteq_c(2^\omega ,B)$.\bigskip

\noindent\bf Case 2.2\rm\ $R$ is separable from $S^{-1}$ by a 
$\mbox{pot}({\bf\Gamma})$ set.\bigskip

 We argue as in Case 1.2 to see that $\big( 2^\omega ,s(R)\big)\sqsubseteq_c(2^\omega ,B)$.
\hfill{$\square$}\bigskip

\noindent\bf Remark.\rm ~This shows that, under the same assumptions, 
$\big\{ R,R\cup\overline{R}^{-1},s(R)\big\}$ is a basis for the class of non-$\mbox{pot}(\bormtwo )$ pot$(F_\sigma )$ Acyclic digraphs. Indeed, $R\cup (\overline{R}^{-1}\!\setminus\! R^{-1})$ is not pot$(F_\sigma )$. 

\begin{thm} \label{baspi0>1} Let $\bf\Gamma$ be a non self-dual Borel class of rank at least two, 
${\cal I}\!\subseteq\! 2^\omega$ given by Lemma \ref{exist}, and 
$\bortwo\!\subseteq\! {\bf\Gamma}'\!\subseteq\! F_\sigma$ be Borel class. We assume that 
$R\! :=\! S_{\cal I}$ has the $({\bf\Gamma},{\bf\Gamma}')$-basis property. Then $\cal{A}$ is a basis for the class of non-$\mbox{pot}({\bf\Gamma})$ Borel subsets of a $\mbox{pot}({\bf\Gamma}')$ symmetric acyclic relation.\end{thm}

\noindent\bf Proof.\rm ~By Theorem \ref{compl}, all the examples are in the context of the theorem. So let $B$ be a non-pot$({\bf\Gamma})$ Borel relation on a Polish space $X$, contained in a 
pot$({\bf\Gamma}')$ symmetric acyclic relation. Note that $B\!\setminus\!\Delta (X)$ is a 
non-pot$({\bf\Gamma})$ Borel relation on $X$, contained in a pot$({\bf\Gamma}')$ acyclic graph.

\vfill\eject

 Theorem \ref{fourpi0>1} gives $A$ in 
$\{ R,R\cup\overline{R}^{-1},R\cup (\overline{R}^{-1}\!\setminus\! R^{-1}),s(R)\}$ reducible to 
$B\!\setminus\!\Delta (X)$ with witness $f$. We set 
$B'\! :=\! (f\!\times\! f)^{-1}(B)$, so that $(2^\omega ,B')\sqsubseteq_c(X,B)$ and 
$A\!\subseteq\! B'\!\subseteq\! A\cup\Delta (2^\omega )$. This means that we may assume that 
$X\! =\! 2^\omega$ and there is a Borel subset $J$ of $2^\omega$ such that 
$B\! =\! A\cup\Delta (J)$. We set, for $\varepsilon\!\in\! 2$, 
$S_\varepsilon\! :=\!\{\alpha\!\in\! 2^\omega\mid\varepsilon\alpha\!\in\! J\}$. This defines a partition 
$\{ S_0\cap S_1,S_0\!\setminus\! S_1,S_1\!\setminus\! S_0,(\neg S_0)\cap (\neg S_1)\}$ of 
$2^\omega$ into Borel sets. By Baire's theorem, one of these sets is not meager.\bigskip

\noindent\bf Claim\it\ Let $s\!\in\! 2^{<\omega}$, $C$ be a dense $G_\delta$ subset of 
$2^\omega$, and $e\!\in\!\{\sqsubset ,\sqsupset\}$. Then\smallskip 

(a) $A\cap\big( 2\!\times\! (N_s\cap C)\big)^2$ is not $\mbox{pot}({\bf\Gamma})$,\smallskip 

(b) $(2^\omega ,s(R)^e)\sqsubseteq_c\big( 2\!\times\! N_s,s(R)^e\cap (2\!\times\! N_s)^2\big)$.\bigskip\rm

\noindent (a) Indeed, $A\cap\big( 2\!\times\! (N_s\cap C)\big)^2\cap (N_0\!\times\! N_1)\! =\! 
R\cap\big( 2\!\times\! (N_s\cap C)\big)^2$. It remains to apply Theorem \ref{compl}.\bigskip

\noindent (b) By (a), $R\cap (2\!\times\! N_s)^2$ is not $\mbox{pot}({\bf\Gamma})$, and it is reducible to $R$. By Corollary \ref{generalanti} and Theorem \ref{fourpi0>1}, $R$ is minimal among non-$\mbox{pot}({\bf\Gamma})$ sets, so that 
$(2^\omega ,R)\sqsubseteq_c\big( 2\!\times\! N_s,R\cap (2\!\times\! N_s)^2\big)$ with witness $f$. Note that $f$ is a homomorphism from $\lceil T\rceil$ into itself, by density. In particular, $f$ sends $N_\varepsilon$ into itself for each $\varepsilon\!\in\! 2$. This shows that $f$ reduces 
$s(R)^e$ to $s(R)^e\cap (2\!\times\! N_s)^2$.\hfill{$\diamond$}\bigskip

\noindent\bf Case 1\rm\ $S_0\cap S_1$ is not meager.\bigskip

 Let $s\!\in\! 2^{<\omega}$ and $C$ be a dense $G_\delta$ subset of $2^\omega$ such that 
$N_s\cap C\!\subseteq\! S_0\cap S_1$. We set 
$$A'\! :=\! A\cap\big( 2\!\times\! (N_s\cap C)\big)^2\mbox{,}$$ 
so that $\big( 2\!\times\! (N_s\cap C),A'\big)\sqsubseteq_c(2^\omega ,A)$. The claim implies that $A'$ is not $\mbox{pot}({\bf\Gamma})$. Corollary \ref{generalanti} and Theorem \ref{fourpi0>1} show that $A$ is minimal among non-$\mbox{pot}({\bf\Gamma})$ sets, so that 
$(2^\omega ,A)\sqsubseteq_c\big( 2\!\times\! (N_s\cap C),A'\big)$ with witness $f'$. The map $f'$ is also a witness for $(2^\omega ,A^\square )\sqsubseteq_c
\Big( 2\!\times\! (N_s\cap C),A'\cup\Delta\big( 2\!\times\! (N_s\cap C)\big)\Big)$. Now 
$(2^\omega ,A^\square )\sqsubseteq_c(2^\omega ,B)$ since 
$B\cap\big( 2\!\times\! (N_s\cap C)\big)^2\! =\! A'\cup\Delta\big( 2\!\times\! (N_s\cap C)\big)$.\bigskip

\noindent\bf Case 2\rm\ $S_0\!\setminus\! S_1$ is not meager.\bigskip

 As in Case 1 we get $s$, $C$ with $N_s\cap C\!\subseteq\! S_0\!\setminus\! S_1$, $A'$, $f'$. The map $f'$ is also a witness for $(2^\omega ,A^\sqsubset )\sqsubseteq_c
\Big( 2\!\times\! (N_s\cap C),A'\cup\Delta\big(\{ 0\}\!\times\! (N_s\cap C)\big)\Big)$ if 
$A\!\not=\! s(R)$, for topological complexity reasons. If $A\! =\! s(R)$, then we can find 
$t\!\in\! 2^{<\omega}$ and $e\!\in\!\{\sqsubset ,\sqsupset\}$ such that 
$$\big( 2\!\times\! N_t,A^e\cap (2\!\times\! N_t)^2\big)\sqsubseteq_c
\Big( 2\!\times\! (N_s\cap C),A'\cup\Delta\big(\{ 0\}\!\times\! (N_s\cap C)\big)\Big) .$$ 
Now note that 
$B\cap\big( 2\!\times\! (N_s\cap C)\big)^2\! =\! A'\cup\Delta\big(\{ 0\}\!\times\! (N_s\cap C)\big)$, so that $(2^\omega ,A^\sqsubset )\sqsubseteq_c(2^\omega ,B)$. Indeed, by Proposition 
\ref{exclude2} and Theorem \ref{fourpi0>1}, $R$ is reducible to $R^{-1}$ since $R$ is contained in a closed Acyclic oriented graph, which is not the case of 
$R\cup (\overline{R}^{-1}\!\setminus\! R^{-1})$. This implies that $s(R)^\sqsubset$ is reducible to $s(R)^\sqsupset$. It remains to note that 
$(2^\omega ,s(R)^e)\sqsubseteq_c\big( 2\!\times\! N_t,s(R)^e\cap (2\!\times\! N_t)^2\big)$, by the claim.\bigskip

\noindent\bf Case 3\rm\ $S_1\!\setminus\! S_0$ is not meager.\bigskip

 We argue as in Case 2 to see that $(2^\omega ,A^\sqsupset )\sqsubseteq_c(2^\omega ,B)$.
 
\vfill\eject

\noindent\bf Case 4\rm\ $(\neg S_0)\cap (\neg S_1)$ is not meager.\bigskip

 As in Case 1 we get $s$, $C$ with $N_s\cap C\!\subseteq\! (\neg S_0)\cap (\neg S_1)$, $A'$. Now note that $B\cap\big( 2\!\times\! (N_s\cap C)\big)^2\! =\! A'$, so that 
$(2^\omega ,A^=)\sqsubseteq_c(2^\omega ,B)$.\hfill{$\square$}\bigskip

\noindent\bf Remark.\rm ~This shows that, under the same assumptions, 
$$\big\{ A^e\mid A\!\in\!\{ R,R\cup\overline{R}^{-1}\}\ \wedge\ 
e\!\in\!\{ =,\square ,\sqsubset ,\sqsupset\}\big\}\cup\big\{ s(R)^e\mid 
e\!\in\!\{ =,\square ,\sqsubset\}\big\}$$ 
is a basis for the class of non-$\mbox{pot}(\bormtwo )$ pot$(F_\sigma )$ Acyclic relations.\bigskip

\noindent $\underline{\mbox{\bf Conditions implying Theorem \ref{rect}}}$

\begin{lem} \label{gorect} Let $\bf\Gamma$ be a Borel class. Assume that\smallskip

\noindent (1) $O$ is a $\check {\bf\Gamma}$ relation on $2^\omega$,\smallskip

\noindent (2) $O$ is contained in a closed Acyclic oriented graph 
$H\!\subseteq\! N_0\!\times\! N_1$,\smallskip

\noindent (3) $O$ is minimum  among non-pot$({\bf\Gamma})$ Borel subsets of a 
pot$(F_\sigma )$ Acyclic oriented graph,\smallskip

\noindent (4) $N_\varepsilon\!\subseteq\!\overline{\Pi_\varepsilon[O]}$.\smallskip

 Then $S\! :=\!\{ (\alpha ,\beta )\!\in\! 2^\omega\!\times\! 2^\omega\mid (0\alpha ,1\beta )\!\in\! O\}$ satisfies the conclusion of Theorem \ref{rect}.\end{lem}

\noindent\bf Proof.\rm ~We set $O'\! :=\! S$. As $O\!\in\!\check {\bf\Gamma}$, 
$S\!\in\!\check {\bf\Gamma}$. As $O$ is contained in $H$, $O'$ is contained in the closed set 
$C\! :=\! H'$. As $H\!\subseteq\! N_0\!\times\! N_1$, the map 
$\varepsilon z\!\mapsto\! (\varepsilon ,z)$ is an isomorphism from $H$ onto $G_{H'}$. Thus 
$s(G_C)$ is acyclic since it is isomorphic to $s(H)$. The shift maps 
$s_\varepsilon\! :\!\varepsilon z\!\mapsto\! z$ defined on $N_\varepsilon$ satisfy 
$O\! =\! (s_0\!\times\! s_1)^{-1}(O')$, which shows that $O'$ is not pot$({\bf\Gamma})$. This shows that (a) and (b) cannot hold simultaneously.\bigskip

 Note that $G_B$ is a Borel oriented graph on the Polish space $X\!\oplus\! Y$ contained in 
$G_F$, which is a pot$(F_\sigma )$ Acyclic oriented graph since the map 
$(\varepsilon ,z)\!\mapsto\! z$ reduces $G_F$ to $F$ on 
$(\{ 0\}\!\times\! X)\!\times\! (\{ 1\}\!\times\! X)$. Assume that $B$ is not pot$({\bf\Gamma})$. Then 
$G_B$ is not pot$({\bf\Gamma})$ since the maps $z\!\mapsto\! (\varepsilon ,z)$ reduce $B$ to $G_B$. As $O$ is minimum, we get $i\! :\! 2^\omega\!\mapsto\! X\!\oplus\! Y$ injective continuous such that $O\! =\! (i\!\times\! i)^{-1}(G_B)$. It remains to set $f(\alpha )\! :=\! i_1(0\alpha )$ and 
$g(\beta )\! :=\! i_1(1\beta )$. Indeed, if $\alpha\!\in\! N_\varepsilon$, then $\alpha$ is the limit of points of $\Pi_\varepsilon [O]$, so that $i_0(\alpha )\! =\!\varepsilon$.\hfill{$\square$}

\section{$\!\!\!\!\!\!$ Study when the rank of $\bf\Gamma$ is at least three} \label{cor>=3}

\noindent\bf Theorem 5\it\ Let $\bf\Gamma$ be a non self-dual Borel class of rank at least three, 
${\cal I}\!\subseteq\! 2^\omega$ given by Lemma \ref{exist}, and $R\! :=\! S_{\cal I}$.\smallskip

(a) the set $\cal A$ defined in Theorem \ref{generalpos} is a basis for the class of 
non-$\mbox{pot}({\bf\Gamma})$ Borel subsets of a $\mbox{pot}(F_\sigma )$ Acyclic relation.\smallskip

(b) $R$ is minimum among non-$\mbox{pot}({\bf\Gamma})$ Borel subsets of a 
$\mbox{pot}(F_\sigma )$ Acyclic oriented graph.\smallskip

(c) $s(R)$ is minimum among non-$\mbox{pot}({\bf\Gamma})$ Borel graphs contained in a 
$\mbox{pot}(F_\sigma )$ acyclic graph.\smallskip

(d) $R\cup\Delta (2^\omega )$ is minimum among non-$\mbox{pot}({\bf\Gamma})$ Borel quasi-orders (or partial orders) contained in a $\mbox{pot}(F_\sigma )$ Acyclic relation.

\vfill\eject

\noindent\bf Proof.\rm ~(a) We apply Theorem \ref{baspi0>1} to ${\bf\Gamma}'\! :=\! F_\sigma$. This is possible, by the remark before Theorem \ref{fourpi0>1}.\bigskip

\noindent (b) Assume that $B$ is a non-$\mbox{pot}({\bf\Gamma})$ Borel subset of a 
$\mbox{pot}(F_\sigma )$ Acyclic oriented graph. By (a), $R$ or 
$R\cup (\overline{R}^{-1}\!\setminus\! R^{-1})$ is reducible to $B$ since $B$ is an oriented graph. It cannot be $R\cup (\overline{R}^{-1}\!\setminus\! R^{-1})$, which is not contained in a 
$\mbox{pot}(F_\sigma )$ Acyclic oriented graph since $R$ is not $\mbox{pot}(F_\sigma )$.\bigskip

\noindent (c) We apply Lemma \ref{symm} and (b).\bigskip

\noindent (d) As $R\!\subseteq\! N_0\!\times\! N_1$, $R\cup\Delta (2^\omega )$ is a Borel quasi-order. By (a), $R\cup\Delta (2^\omega )$ is not $\mbox{pot}({\bf\Gamma})$ and is contained in a 
$\mbox{pot}(F_\sigma )$ Acyclic relation. Assume that $Q$ is a non-$\mbox{pot}({\bf\Gamma})$ Borel quasi-order on a Polish space $X$, contained in a $\mbox{pot}(F_\sigma )$ Acyclic relation. (a) gives $A\!\in\! {\cal A}$ with $(2^\omega ,A)\sqsubseteq_c(X,Q)$. As $Q$ is reflexive, $A$ has to be reflexive too, so that $e\! =\!\square$. We saw that $\cal I$ can be a free ideal if the rank of 
$\bf\Gamma$ is infinite or if 
$\bf\Gamma\!\in\!\{\bormtwo ,\borathree ,{\bf\Pi}^0_4,{\bf\Sigma}^0_5,...\}$. Note that 
$S_{\cal I}\! =\!\{ (\alpha ,\beta )\!\in\!\lceil T\rceil\mid\alpha\Delta\beta\!\in\! {\cal I}\}$, 
$\overline{R}\! =\!\lceil T\rceil$ since $S_{\cal I}$ is dense in $\lceil T\rceil$.\bigskip

 If the rank of $\bf\Gamma$ is infinite or if 
${\bf\Gamma}\!\in\!\{\bormtwo ,\borathree ,{\bf\Pi}^0_4,{\bf\Sigma}^0_5,...\}$, then 
$(0^\infty ,1^20^\infty ),(010^\infty ,1^20^\infty )\!\in\! R$, but 
$(0^\infty ,010^\infty )\!\notin\! s(R)\cup\Delta (2^\omega )$, so that $s(R)\cup\Delta (2^\omega )$ is not transitive. Pick $(0\alpha ,1\beta )\!\in\! S_{\neg {\cal I}}$, which is dense in $\lceil T\rceil$. Then 
$(0\beta ,1\beta ),(1\beta ,0\alpha )\!\in\! R\cup (\overline{R}^{-1}\!\setminus\! R^{-1})$, and 
$(0\beta ,0\alpha )\!\notin\! R\cup\overline{R}^{-1}\cup\Delta (2^\omega )$ since $\beta\!\not=\!\alpha$, so that $R\cup\overline{R}^{-1}\cup\Delta (2^\omega )$ and 
$R\cup (\overline{R}^{-1}\!\setminus\! R^{-1})\cup\Delta (2^\omega )$ are not transitive. This shows that 
$A\! =\! R\cup\Delta (2^\omega )$.\bigskip

 If ${\bf\Gamma}\!\in\!\{\boratwo ,\bormthree ,{\bf\Sigma}^0_4,{\bf\Pi}^0_5,...\}$, then $\cal I$ can be the complement of the set $\cal I$ previously considered. As $R$ is not 
$\mbox{pot}({\bf\Gamma})$, there are $\alpha ,\beta ,\gamma$ with $\beta\!\not=\!\gamma$ and 
$(0\alpha ,1\beta ),(0\alpha ,1\gamma )\!\in\! R$. Then $(1\gamma ,0\alpha )\!\in\! s(R)$ and 
$(1\gamma ,1\beta )\!\notin\! s(R)\cup\Delta (2^\omega )$, so that $s(R)\cup\Delta (2^\omega )$ is not transitive. Pick $(0\alpha ,1\beta )\!\in\! S_{\cal I}$, which is dense in $\lceil T\rceil$. Then 
$(0\alpha ,1\beta ),(1\beta ,0\beta )\!\in\! R\cup (\overline{R}^{-1}\!\setminus\! R^{-1})$, and 
$(0\alpha ,0\beta )\!\notin\! R\cup\overline{R}^{-1}\cup\Delta (2^\omega )$ since $\beta\!\not=\!\alpha$, so that $R\cup\overline{R}^{-1}\cup\Delta (2^\omega )$ and 
$R\cup (\overline{R}^{-1}\!\setminus\! R^{-1})\cup\Delta (2^\omega )$ are not transitive. This shows that 
$A\! =\! R\cup\Delta (2^\omega )$ again.\hfill{$\square$}\bigskip

\noindent\bf Proof of Theorem \ref{generalpos} when the rank is at least three.\rm ~Fix $\cal I$ given by Lemma \ref{exist}. We set $R\! :=\! S_{\cal I}$.\bigskip

\noindent (1) We apply Theorem \ref{compl}.\bigskip

\noindent (2) We apply Theorem \ref{compl}, Corollary \ref{generalanti}, and the beginning of Section 3 (which ensures that $s(\lceil T\rceil )$ is acyclic).\bigskip

\noindent (3) We apply Theorem \ref{cor>=3}.\hfill{$\square$}\bigskip

\noindent\bf Proof of Theorem \ref{generalneg} when the rank is at least three.\rm ~(2) We argue by contradiction, which gives $O$. Lemma \ref{exist} gives $\cal I$. By Theorem \ref{cor>=3}, $O$ is reducible to $S_{\cal I}$, $O$ is contained in a closed Acyclic oriented graph, and $S_{\cal I}$ is reducible to $O$. By Theorem \ref{compl}, 
$U_{\cal I}\! :=\! S_{\cal I}\cup (\lceil T\rceil^{-1}\!\setminus\! E_{\cal I})$ is a 
$({\bf\Gamma}\oplus\check {\bf\Gamma})\!\setminus\!\mbox{pot}({\bf\Gamma})$ Acyclic oriented graph. Thus $S_{\cal I}$ is reducible to $U_{\cal I}$, which contradicts Corollary 
\ref{generalanti}.\bigskip

\noindent (1) We argue by contradiction, which gives $O'$. Lemma \ref{exist} gives $\cal I$. By Theorem \ref{cor>=3}, $O'$ is reducible to $S_{\cal I}$, and $O'$ is contained in a closed Acyclic oriented graph. Thus $O'$ is minimum among Borel relations, contained in a 
$\mbox{pot}({\bf\Gamma}\oplus\check {\bf\Gamma})$ Acyclic oriented graph, which are not 
$\mbox{pot}({\bf\Gamma})$. We just saw that this cannot happen.\hfill{$\square$}

\vfill\eject

\noindent\bf Proof of Theorem \ref{rect} when the rank is at least three.\rm ~We apply Theorem 
\ref{cor>=3} and Lemmas \ref{exist}, \ref{gorect}. Lemma \ref{gorect} is applied to 
$O\! :=\! S_{\cal I}$ and $H\! :=\!\lceil T\rceil$. If $\alpha\!\in\! N_\varepsilon$ and $s$ is the shift map, then $\big( 0s(\alpha ),1s(\alpha )\big)$ is in $\lceil T\rceil$ and is the limit of points of $O$.\hfill{$\square$}

\section{$\!\!\!\!\!\!$ Study when the rank of $\bf\Gamma$ is two}\indent 

 We start with a consequence of Corollary 6.4 in [L-Z].

\begin{cor} \label{cordelta02} Let ${\bf\Gamma}'$ be a class of Borel sets closed under continuous pre-images. Then $\lceil T\rceil\cap\mathbb{E}_0$ has the $(\bormtwo ,{\bf\Gamma}')$-basis property.\end{cor}

\noindent\bf Proof.\rm ~By Theorem \ref{compl}, $\lceil T\rceil\cap\mathbb{E}_0$ is not 
$\mbox{pot}(\bormtwo )$, so that (a) and (b) cannot hold simultaneously. It remains to apply Corollary 6.4 in [L-Z].\hfill{$\square$}

\begin{thm} \label{cor=2pi} Let $R\! :=\!\lceil T\rceil\cap\mathbb{E}_0$.\smallskip

(a) the set $\cal A$ defined in Theorem \ref{generalpos} is a basis for the class of 
non-$\mbox{pot}(\bormtwo )$ Borel subsets of a $\mbox{pot}(F_\sigma )$ Acyclic relation.\smallskip

(b) $R$ is minimum among non-$\mbox{pot}(\bormtwo )$ Borel subsets of a 
$\mbox{pot}(F_\sigma )$ Acyclic oriented graph.\smallskip

(c) $s(R)$ is minimum among non-$\mbox{pot}(\bormtwo )$ Borel graphs contained in a 
$\mbox{pot}(F_\sigma )$ acyclic graph.\smallskip

(d) $R\cup\Delta (2^\omega )$ is minimum among non-$\mbox{pot}(\bormtwo )$ Borel quasi-orders (or partial orders) contained in a $\mbox{pot}(F_\sigma )$ Acyclic relation.\end{thm}

\noindent\bf Proof.\rm ~(a) We apply Theorem \ref{baspi0>1} to ${\bf\Gamma}'\! :=\! F_\sigma$. This is possible, by Corollary \ref{cordelta02}.\bigskip

\noindent (b) Assume that $B$ is a non-$\mbox{pot}(\bormtwo )$ Borel subset of a 
$\mbox{pot}(F_\sigma )$ Acyclic oriented graph. By (a), $R$ or 
$R\cup (\overline{R}^{-1}\!\setminus\! R^{-1})$ is reducible to $B$ since $B$ is an oriented graph. It cannot be $R\cup (\overline{R}^{-1}\!\setminus\! R^{-1})$, which is not contained in a 
$\mbox{pot}(F_\sigma )$ Acyclic oriented graph since $R$ is not $\mbox{pot}(\bortwo )$.\bigskip

\noindent (c) We apply Lemma \ref{symm} and (b).\bigskip

\noindent (d) We argue as in the proof of Theorem 5.\hfill{$\square$}\bigskip

\noindent\bf Proof of Theorem \ref{generalpos}.(1)-(3) when ${\bf\Gamma}\! =\!\bormtwo$.\rm ~We set 
$R\! :=\!\lceil T\rceil\cap\mathbb{E}_0$, and argue as when the rank of ${\bf\Gamma}$ is at least three (we just have to replace Theorem \ref{cor>=3} with Theorem \ref{cor=2pi}).\hfill{$\square$}\bigskip

\noindent\bf Proof of Theorems \ref{generalneg} and \ref{rect} when ${\bf\Gamma}\! =\!\bormtwo$.\rm ~We argue as when the rank of ${\bf\Gamma}$ is at least three (we just have to replace Theorem \ref{cor>=3} with Theorem \ref{cor=2pi}).\hfill{$\square$}

\begin{thm} \label{cor=2sigma} Let $R\! :=\!\lceil T\rceil\!\setminus\!\mathbb{E}_0$.\smallskip

(a) the set $\cal A$ defined in Theorem \ref{generalpos} is a basis for the class of 
non-$\mbox{pot}(\boratwo )$ Borel subsets of a $\mbox{pot}(F_\sigma )$ Acyclic relation.\smallskip

(b) $R$ is minimum among non-$\mbox{pot}(\boratwo )$ Borel subsets of a 
$\mbox{pot}(F_\sigma )$ Acyclic oriented graph.\smallskip

(c) $s(R)$ is minimum among non-$\mbox{pot}(\boratwo )$ Borel graphs contained in a 
$\mbox{pot}(F_\sigma )$ acyclic graph.\smallskip

(d) $R\cup\Delta (2^\omega )$ is minimum among non-$\mbox{pot}(\boratwo )$ Borel quasi-orders (or partial orders) contained in a $\mbox{pot}(F_\sigma )$ Acyclic relation.\end{thm}

\noindent\bf Proof.\rm ~(a) Let us check that $R$ has the $(\boratwo ,F_\sigma )$-basis property. Let $X$ be a Polish space, and $A,B$ be disjoint analytic relations on $X$ such that 
$A$ is contained in a $\mbox{pot}(F_\sigma )$ symmetric acyclic relation $F$. Note first that $R$ is not separable from $\overline{R}\!\setminus\! R$ by a $\mbox{pot}(\boratwo )$ set, by Theorem 
\ref{compl}. So assume that $A$ is not separable from $B$ by a $\mbox{pot}(\boratwo )$ set. Note that $A$ is not separable from $B\cap F$ by a $\mbox{pot}(\boratwo )$ set. Corollary 
\ref{cordelta02} gives $g\! :\! 2^\omega\!\rightarrow\! X$ injective continuous such that 
$\lceil T\rceil\cap\mathbb{E}_0\!\subseteq\! (g\!\times\! g)^{-1}(B\cap F)$ 
and $\lceil T\rceil\!\setminus\!\mathbb{E}_0\!\subseteq\! (g\!\times\! g)^{-1}(A)$, and we are done.\bigskip

We can now apply Theorem \ref{baspi0>1} to ${\bf\Gamma}'\! :=\! F_\sigma$.\bigskip

\noindent (b) Assume that $B$ is a non-$\mbox{pot}(\boratwo )$ Borel subset of a 
$\mbox{pot}(F_\sigma )$ Acyclic oriented graph. By (a), $R$ or 
$R\cup (\overline{R}^{-1}\!\setminus\! R^{-1})$ is reducible to $B$ since $B$ is an oriented graph. It cannot be $R\cup (\overline{R}^{-1}\!\setminus\! R^{-1})$, which is not contained in a 
$\mbox{pot}(F_\sigma )$ Acyclic oriented graph since $R$ is not $\mbox{pot}(F_\sigma )$.\bigskip

\noindent (c) We apply Lemma \ref{symm} and (b).\bigskip

\noindent (d) We argue as in the proof of Theorem 5.\hfill{$\square$}\bigskip

\noindent\bf Proof of Theorem \ref{generalpos}.(1)-(3) when ${\bf\Gamma}\! =\!\boratwo$.\rm ~We set $R\! :=\!\lceil T\rceil\!\setminus\!\mathbb{E}_0$, and argue as when the rank of 
${\bf\Gamma}$ is at least 3 (we just have to replace Theorem \ref{cor>=3} with Theorem 
\ref{cor=2sigma}).\hfill{$\square$}\bigskip

\noindent\bf Proof of Theorems \ref{generalneg} and \ref{rect} when ${\bf\Gamma}\! =\!\boratwo$.\rm ~We argue as when the rank of ${\bf\Gamma}$ is at least three (we just have to replace Theorem \ref{cor>=3} with Theorem \ref{cor=2sigma}).\hfill{$\square$}\bigskip

  If we add an acyclicity assumption to Corollary 6.5 in [L-Z], then we get a reduction on the whole product, namely Theorem \ref{generalpos}.(4). We can prove it using Corollary 6.5 in [L-Z], but in fact it is just a corollary of Theorem \ref{generalpos}.(3).\bigskip

\noindent\bf Proof of Theorem \ref{generalpos}.(4).\rm ~We apply the fact, noted in the introduction, that a Borel locally countable relation is pot$(F_\sigma )$, and Theorem \ref{generalpos}.(3). We use the fact that $R\cup\overline{R}^{-1}$ and $R\cup (\overline{R}^{-1}\!\setminus\! R^{-1})$ are not localy countable.\hfill{$\square$}
   
\section{$\!\!\!\!\!\!$ Study when the rank of $\bf\Gamma$ is one}\indent 

 We first study the case ${\bf\Gamma}\! =\!\boraone$.\bigskip

\noindent\bf Proof of Theorem \ref{generalpos}.(6).\rm ~As the $\mbox{pot}(\boraone )$ sets are exactly the countable unions of Borel rectangles, 
$\Delta (2^\omega ),\mbox{Gr}({h_0}_{\vert N_0}),\mbox{Gr}(h_0)$ are not 
$\mbox{pot}(\boraone )$. Note that these relations are closed and Acyclic since 
$\mbox{Gr}(h_0)$ is acyclic. Considerations about reflexivity and Proposition 
\ref{smallantisigma01} show that these relations form a $\sqsubseteq_c$-antichain. So assume that $B$ is a non-$\mbox{pot}(\boraone )$ Borel Acyclic relation, so that $B$ is not a countable union of Borel rectangles.\bigskip

 If $\{ x\!\in\! X\mid (x,x)\!\in\! B\}$ is uncountable, then it contains a Cantor set $C$. Lemmas 
\ref{suffacy} and \ref{meager} show that $B\cap C^2$ is meager in $C^2$. Mycielski's theorem gives a Cantor subset $K$ of $C$ such that $K^2\cap B\! =\!\Delta (K)$ (see 19.1 in [K]). This implies that $\big( 2^\omega ,\Delta (2^\omega )\big)\sqsubseteq_c(X,B)$.

\vfill\eject

 So we may assume that $\{ x\!\in\! X\mid (x,x)\!\in\! B\}$ is countable, and in fact that $B$ is irreflexive. As $B$ is not a countable union of Borel rectangles, we can find Cantor subsets 
$C$, $D$ of $X$ and a homeomorphism $\varphi\! :\! C\!\rightarrow\! D$ whose graph is contained in $B$ (see [P]). As $B$ is irreflexive, $\varphi$ is fixed point free and we may assume that $C$ and $D$ are disjoint. Let $\Psi_0\! :\! N_0\!\rightarrow\! C$ be a homeomorphism, and 
$\Psi_1\! :=\!\varphi\circ\Psi_0\circ {h_0}_{\vert N_1}$, so that 
$\Psi_1\! :\! N_1\!\rightarrow\! D$ is a homeomorphism too. We set 
$\Psi (\alpha )\! :=\!\Psi_\varepsilon (\alpha )$ if $\alpha\!\in\! N_\varepsilon$, so that 
$\Psi\! :\! 2^\omega\!\rightarrow\! X$ is a continuous injection.\bigskip

 We also set $B'\! :=\! (\Psi\!\times\!\Psi )^{-1}(B)$, so that $B'$ is a relation on $2^\omega$ containing $\mbox{Gr}({h_0}_{\vert N_0})$ and satisfying the same properties as $B$. By Lemmas 
\ref{suffacy} and \ref{meager}, $B'$ is meager. Let 
$\varphi_\varepsilon\! :\! 2^\omega\!\rightarrow\! N_\varepsilon$ be the homeomorphism defined by $\varphi_\varepsilon (\alpha )\! :=\!\varepsilon\alpha$, and 
$B''\! :=\!\bigcup_{\varepsilon ,\varepsilon'\in 2}~
(\varphi_\varepsilon\!\times\!\varphi_{\varepsilon'})^{-1}(B')$, so that $B''$ is a reflexive meager relation on $2^\omega$. Mycielski's theorem gives a Cantor subset $K$ of $2^\omega$ such that $K^2\cap B''\! =\!\Delta (K)$. Let $h\! :\! 2^\omega\!\rightarrow\! K$ be a homeomorphism, and $g(\varepsilon\alpha )\! :=\!\varphi_\varepsilon\big( h(\alpha )\big)$. Then 
$g$ is injective continuous. We set $B'''\! :=\! (g\!\times\! g)^{-1}(B')$, so that 
$\mbox{Gr}({h_0}_{\vert N_0})\!\subseteq\! B'''\!\subseteq\!\mbox{Gr}(h_0)$. We then set 
$S\! :=\!\{\alpha\!\in\! 2^\omega\mid (1\alpha ,0\alpha )\!\in\! B'''\}$.\bigskip

 If $S$ is meager, then let $P$ be a Cantor subset disjoint from $S$. Then 
$$B'''\cap (2\!\times\! P)^2\! =\!\mbox{Gr}({h_0}_{\vert N_0})\cap (2\!\times\! P)^2$$ 
is a non-pot$(\boraone )$ Acyclic oriented graph on $2\!\times\! P$, and, repeating the previous discussion, we see that 
$$\big( 2^\omega ,\mbox{Gr}({h_0}_{\vert N_0})\big)\sqsubseteq_c\big( 2\!\times\! P^2,
B'''\cap (2\!\times\! P)^2\big)\sqsubseteq_c(2^\omega ,B''')\sqsubseteq_c(X,B).$$ 
Similarly, if $S$ is not meager, then let $Q$ be a Cantor subset of $S$. Then 
$$B'''\cap (2\!\times\! Q)^2\! =\!\mbox{Gr}(h_0)\cap (2\!\times\! Q)^2$$ 
is a non-pot$(\boraone )$ acyclic graph on $2\!\times\! Q$, and, repeating the previous discussion, we see that 
$$\big( 2^\omega ,\mbox{Gr}(h_0)\big)\sqsubseteq_c\big( 2\!\times\! Q^2,
B'''\cap (2\!\times\! Q)^2\big)\sqsubseteq_c(2^\omega ,B''')\sqsubseteq_c(X,B).$$
For the last assertion, let $Q$ be a non-$\mbox{pot}({\bf\Gamma})$ Borel Acyclic quasi-order on a Polish space $X$. Theorem \ref{generalpos} gives $A\!\in\!\{ \Delta (2^\omega ),R,s(R)\}$ with 
$(2^\omega ,A)\sqsubseteq_c(X,Q)$. As $R$ and $s(R)$ are not reflexive, $A$ has to be 
$\Delta (2^\omega )$.\hfill{$\square$}\bigskip

 If we apply Theorem \ref{generalpos}.(6) and Lemma \ref{gorect}, then we get a version of Theorem \ref{rect} for ${\bf\Gamma}\! =\!\boraone$. Let us mention a corollary in the style of Corollary 6.4 in [L-Z]. 
 
\begin{cor} \label{partialopen} Let $X$ be a Polish space, and $B$ be a Borel Acyclic relation on $X$. Then exactly one of the following holds:\smallskip  

(a) the set $B$ is $\mbox{pot}(\boraone )$,\smallskip  

(b) there are $f,g\! :\! 2^\omega\!\rightarrow\! X$ injective continuous with 
$\Delta (2^\omega )\! =\! (f\!\times\! g)^{-1}(B)$.\end{cor}
 
 We now study the case ${\bf\Gamma}\! =\!\bormone$. We will apply several times Corollary 3.10 in [L-Z] and use the following lemma.

\begin{lem} \label{values} Let $X$ be a Polish space, $B$ be a relation on $X$, $C,D$ be closed subsets of $X$, and $f,g\! :\! 2^\omega\!\rightarrow\! X$ be continuous maps such that 
$\mathbb{G}_0\!\subseteq\! (f\!\times\! g)^{-1}\big( B\cap (C\!\times\! D)\big)$. Then $f$ (resp., $g$) takes values in $C$ (resp., $D$).\end{lem}

\noindent\bf Proof.\rm ~The first projection of $\mathbb{G}_0$ is comeager, so that 
$f(\alpha )\!\in\! C$ for almost all $\alpha$, and all $\alpha$ by continuity. Similarly, 
$g(\beta )\!\in\! D$ for all $\beta$.\hfill{$\square$}

\vfill\eject

 In our results about potentially closed sets, the assumption of being $\mbox{pot}\big(\check D_2(\boraone )\big)$ is equivalent to being 
$\mbox{pot}(\bormone )$, in the acyclic context. We indicate the class $\check D_2(\boraone )$ for optimality reasons.

\begin{prop} \label{simpler} Any $\mbox{pot}\big(\check D_2(\boraone )\big)$ Acyclic relation is 
$\mbox{pot}(\bormone )$.\end{prop}

\noindent\bf Proof.\rm ~Let $G$ be a $\mbox{pot}\big(\check D_2(\boraone )\big)$ Acyclic relation. We can write $G\! =\! O\cup C$, with $O\!\in\!\mbox{pot}(\boraone )$ and 
$C\!\in\!\mbox{pot}(\bormone )$. As $O\!\setminus\!\Delta (X)$ is $\mbox{pot}(\boraone )$, irreflexive and Acyclic, $\big( O\!\setminus\!\Delta (X)\big)\cap (C\!\times\! D)$ is meager in 
$C\!\times\! D$ if $C,D$ are Cantor subsets of $X$ by Lemmas \ref{suffacy} and \ref{meager}, so that we can write $O\!\setminus\!\Delta (X)\! =\!\bigcup_{n\in\omega}~A_n\!\times\! B_n$, with 
$A_n$ or $B_n$ countable for each $n$. In particular, $O\!\setminus\!\Delta (X)$ is 
$\mbox{pot}(\borone )$ by Remark 2.1 in [L1]. Note that $O\!\cap\!\Delta (X)$ is a Borel set with closed vertical sections and is therefore $\mbox{pot}(\bormone )$ (see [Lo1]). Thus 
$O\! =\!\big( O\!\setminus\!\Delta (X)\big)\cup\big( O\cap\Delta (X)\big)$ and $G$ are 
$\mbox{pot}(\bormone )$.\hfill{$\square$}\bigskip

\noindent\bf Proof of Theorem \ref{generalpos}.(1)-(2) and (5).(i) when 
${\bf\Gamma}\! =\!\bormone$.\rm ~(1) By Lemma \ref{propertiespi01}, $R$ is $D_2(\boraone )$, not $\mbox{pot}(\bormone )$, and is Acyclic. By Proposition \ref{simpler}, $R$ is not 
$\mbox{pot}\big(\check D_2(\boraone )\big)$.\bigskip

\noindent (5).(i) Note first that Lemma \ref{suffacy} implies that $\mathbb{B}_0$ is in the context of Theorem \ref{generalpos}.(5).(i), in the sense that it is a Borel subset of the closed Acyclic oriented graph 
$\overline{\mathbb{B}_0}\! =\!\mathbb{B}_0\cup\{ (0\alpha ,1\alpha )\mid\alpha\!\in\! 2^\omega\}$. Assume that $B$ is a non-pot$(\bormone )$ Borel subset of a pot$(\bormone )$ Acyclic oriented graph. Note that there is a Borel countable coloring of $(X,B)$. Indeed, we argue by contradiction. Theorem \ref{G0} gives $f\! :\! 2^\omega\!\rightarrow\! X$ injective continuous such that 
$\mathbb{G}_0\! =\! (f\!\times\! f)^{-1}(B)$. This shows the existence of a $\mbox{pot}(\bormone )$ oriented graph separating $\mathbb{G}_0$ from $\Delta (2^\omega)$. This gives a Borel countable coloring of $(2^\omega ,\mathbb{G}_0)$, which is absurd.\bigskip
 
 This shows the existence of a Borel partition $(B_n)_{n\in\omega}$ of $X$ into $B$-discrete sets. This gives $m\!\not=\! n$ such that 
$B\cap (B_m\!\times\! B_n)$ is not $\mbox{pot}(\bormone )$. We can change the Polish topology, so that we can assume that the $B_n$'s are clopen and $B$ is contained in a closed Acyclic oriented graph $F$. Note that 
$\Big( B_m\cup B_n,\big( B\cap (B_m\!\times\! B_n)\big)\cup\big( B\cap (B_n\!\times\! B_m)\big)\Big)\sqsubseteq_c(X,B)$, and that 
$$F'\! :=\!\big( F\cap (B_m\!\times\! B_n)\big)\cup\big( F\cap (B_n\!\times\! B_m)\big)$$ 
is a closed Acyclic oriented graph on $B_m\cup B_n$ containing ${\big( B\cap (B_m\!\times\! B_n)\big)\cup\big( B\cap (B_n\!\times\! B_m)\big)}$. Corollary 3.10 in [L-Z] gives $f',g'\! :\! 2^\omega\!\rightarrow\! B_m\cup B_n$ injective continuous with 
$$\mathbb{G}_0\!\subseteq\! (f'\!\times\! g')^{-1}\big( B\cap (B_m\!\times\! B_n)\big)$$ 
and $\Delta (2^\omega )\!\subseteq\!\neg (f'\!\times\! g')^{-1}\big( B\cap (B_m\!\times\! B_n)\big)$. By Lemma \ref{values}, $f'(\alpha )\!\in\! B_m$ for all 
$\alpha$, and $g'(\beta )\!\in\! B_n$ for all $\beta$. Thus $\Delta (2^\omega )\!\subseteq\! (f'\!\times\! g')^{-1}(\neg B)$. The shift maps 
$s_\varepsilon\! :\! N_\varepsilon\!\rightarrow\! 2^\omega$, for $\varepsilon\!\in\! 2$, are continuous injections and 
$\mathbb{B}_0\! =\!\overline{\mathbb{B}_0}\cap (s_0\!\times\! s_1)^{-1}(\mathbb{G}_0)$. The map $f''\! :\! N_0\!\rightarrow\! B_m$ (resp., $g''\! :\! N_1\!\rightarrow\! B_n$) defined by 
$f''\! :=\! f'\circ s_0$ (resp., $g''\! :=\! g'\circ s_1$) is injective continuous, 
${\mathbb{B}_0\!\subseteq\! (f''\!\times\! g'')^{-1}(B)}$ and 
${\overline{\mathbb{B}_0}\!\setminus\!\mathbb{B}_0\!\subseteq\! (f''\!\times\! g'')^{-1}(\neg B)}$. We set ${h(\alpha )\! :=\! f''(\alpha )}$ if $\alpha (0)\! =\! 0$, 
${h(\alpha )\! :=\! g''(\alpha )}$ otherwise. Note that ${h\! :\! 2^\omega\!\rightarrow\! B_m\cup B_n}$ is injective continuous, 
$\mathbb{B}_0\!\subseteq\! (h\!\times\! h)^{-1}(B)$ and $\overline{\mathbb{B}_0}\!\setminus\!\mathbb{B}_0\!\subseteq\! (h\!\times\! h)^{-1}(\neg B)$. Moreover, $F''\! :=\! (h\!\times\! h)^{-1}(F')$ is a closed Acyclic oriented graph on $2^\omega$ containing $\mathbb{B}_0$, and contained in 
$(N_0\!\times\! N_1)\cup (N_1\!\times\! N_0)$. Theorem \ref{homB0} gives 
$i\! :\! 2^\omega\!\rightarrow\! 2^\omega$ injective continuous with 
$\mathbb{B}_0\!\subseteq\! (i\!\times\! i)^{-1}(\mathbb{B}_0)$, $\overline{\mathbb{B}_0}\!\setminus\!\mathbb{B}_0\!\subseteq\! (i\!\times\! i)^{-1}(\overline{\mathbb{B}_0}\!\setminus\!\mathbb{B}_0)$, and $\neg\overline{\mathbb{B}_0}\!\subseteq\! (i\!\times\! i)^{-1}(\neg F'')$. Then $f\! :=\! h\circ i$ is an injective continuous reduction of $\mathbb{B}_0$ to $B$.\bigskip

 For $s(\mathbb{B}_0)$, we apply the proof of Lemma \ref{symm} and the previous argument. For the last assertion, we argue as in the proof of Theorem 5 (assuming that Theorem 
\ref{generalpos}.(5).(ii) is proved, which will be done later).\bigskip
 
\noindent (2) We apply Proposition \ref{smallantipi01} and Lemma \ref{propertiespi01}.
\hfill{$\square$}

\vfill\eject

 The proof of Lemma \ref{gorect} and Theorem \ref{generalpos} give the version of Theorem 
\ref{rect} for ${\bf\Gamma}\! =\!\bormone$ announced in the introduction.

\begin{prop} \label{not} $\mathbb{B}_0\not\sqsubseteq_cG_{s(\mathbb{G}_0)}$.\end{prop}

\noindent\bf Proof.\rm ~Assume that $f\! :\! 2^\omega\!\rightarrow\! 2^\omega$ is injective continuous and $\mathbb{B}_0\! =\! (f\!\times\! f)^{-1}(G_{s(\mathbb{G}_0)})$. Let 
$S\! :\! 2^\omega\!\rightarrow\! 2^\omega$ be the shift map defined by 
$S(\varepsilon\alpha )\! :=\!\alpha$. Then the maps $\alpha\!\mapsto\! S\big( f(0\alpha )\big)$ and 
$\beta\!\mapsto\! S\big( f(1\beta )\big)$ define a rectangular continuous reduction of 
$\mathbb{G}_0$ to $s(\mathbb{G}_0)$. Indeed, it is clearly a homomorphism. The first projection of $\mathbb{G}_0$ is comeager, so that $0\!\subseteq\! f(0\alpha )$ for almost all $\alpha$, and all $\alpha$ by continuity. Similarly, $1\!\subseteq\! f(1\beta )$ for all $\beta$, which gives a rectangular reduction. As $\overline{\mathbb{G}_0}\!\setminus\!\mathbb{G}_0\! =\!
\Delta (2^\omega )\! =\!\overline{s(\mathbb{G}_0)}\!\setminus\! s(\mathbb{G}_0)$, we have in fact a square rectangular continuous reduction, which is not possible since $\mathbb{G}_0$ is antisymmetric and $s(\mathbb{G}_0)$ is symmetric.\hfill{$\square$}\bigskip

\noindent\bf Remarks.\rm ~(a) The assumptions ``$F$ is closed" and ``$F$ is Acyclic" in Theorem 
\ref{homB0} are useful. Indeed, for the first one, assume that $F$ is $G_{s(\mathbb{G}_0)}$. Then $F$ satisfies the assumptions of Theorem \ref{homB0}, except that it is not $\bormone$. If the conclusion was true, then we would have $\mathbb{B}_0\sqsubseteq_c G_{s(\mathbb{G}_0)}$, which is absurd by Proposition \ref{not}.\bigskip

 For the second one, assume that $F$ is $\overline{G_{s(\mathbb{G}_0)}}$. Then $F$ satisfies the assumptions of Theorem \ref{homB0}, except that it is not Acyclic. If the conclusion was true, then we would have 
$\overline{\mathbb{B}_0}\sqsubseteq_c\overline{G_{s(\mathbb{G}_0)}}$. As in the proof of Proposition \ref{not}, this would give a rectangular continuous reduction of 
$\overline{\mathbb{G}_0}$ to $\overline{s(\mathbb{G}_0)}$, with witnesses $f',g'$. As in the proof of Proposition \ref{not}, we cannot have $f'\! =\! g'$. The proof of Proposition \ref{not} shows that $f',g'$ are injective. Let $\alpha\!\in\! 2^\omega$ with $f'(\alpha )\!\not=\! g'(\alpha )$. Then for example there is $n\!\in\!\omega$ such that $g'(\alpha )\! =\!\varphi_n\big( f'(\alpha )\big)$ (we use the notation in the proof of Theorem \ref{Fsigma}). In particular, there are clopen sets $U,V$ such that $\big( f'(\alpha ),g'(\alpha )\big)\!\in\! U\!\times\! V$ and 
${\overline{s(\mathbb{G}_0)}\cap (U\!\times\! V)\! =\!\mbox{Gr}(\varphi_n)\cap (U\!\times\! V)}$. We set ${W\! :=\! {f'}^{-1}(U)\cap {g'}^{-1}(V)}$, which is a clopen neighborhood of $\alpha$ such that 
${\overline{\mathbb{G}_0}\cap W^2\! =\! 
(f'\!\times\! g')^{-1}\big(\mbox{Gr}(\varphi_n)\big)\cap W^2}$. Pick $p\!\in\!\omega$, $\beta\!\in\! W$ with $\varphi_p(\beta )\!\in\! W$. Then 
$g'(\beta )\! =\!\varphi_n\big( f'(\beta )\big)\! =\! g'\big(\varphi_p(\beta )\big)$, which contradicts the injectivity of $g'$.\bigskip

\noindent (b) We cannot replace the class $\check D_2(\boraone )$ with $D_2(\bormone )$ in the version of Theorem \ref{rect} for ${\bf\Gamma}\! =\!\bormone$. Indeed, take 
$B\! :=\! s(\mathbb{G}_0)$. Note that 
$B\! =\!\overline{B}\!\setminus\!\Delta (2^\omega )\!\in\! D_2(\bormone )$. Moreover, $B$ is irreflexive, symmetric and acyclic (see Proposition \ref{conn}). Thus $s(G_B)$ is acyclic by Lemma \ref{suffacy}. Theorem \ref{potclo} shows that $B\!\notin\!\mbox{pot}(\bormone )$. The proof of Proposition \ref{not} shows that we cannot find $f,g\! :\! 2^\omega\!\rightarrow\! 2^\omega$ injective continuous with $\mathbb{G}_0\! =\! (f\!\times\! g)^{-1}(B)$.\bigskip

\noindent\bf Proof of Theorem \ref{generalneg}.(1) and (3) when ${\bf\Gamma}\! =\!\bormone$.\rm ~(3) We argue by contradiction, which gives $O$. As 
$\overline{\mathbb{B}_0}$ is a locally countable Acyclic oriented graph, $O$ is reductible to 
$\mathbb{B}_0$, $O$ is contained in a closed Acyclic oriented graph, and $\mathbb{B}_0$ is reducible to $O$. Note that $G_{s(\mathbb{G}_0)}$ is a locally countable $D_2(\boraone )$ 
non-$\mbox{pot}(\bormone )$ Acyclic oriented graph. Indeed, by Lemma \ref{suffacy} and the remark after it, $s(\mathbb{T}_0)$ is acyclic. Thus $\mathbb{B}_0$ is reducible to 
$G_{s(\mathbb{G}_0)}$, which contradicts Proposition \ref{not}.\bigskip

\noindent (1) We argue as when the rank of $\bf\Gamma$ is at least three.\hfill{$\square$}\bigskip

\noindent $\underline{\mbox{\bf An antichain made non-pot$(\bormone )$ relations}}$

\begin{prop} \label{biganti10} $\{\mathbb{G}_0,\mathbb{B}_0,\mathbb{N}_0,\mathbb{M}_0,
G_{s(\mathbb{G}_0)},\mathbb{T}_0,\mathbb{U}_0,s(\mathbb{G}_0),s(\mathbb{B}_0),
s(\mathbb{T}_0)\}$ is a $\sqsubseteq_c$-antichain made of $D_2(\boraone )$ Acyclic digraphs, with locally countable closure, which are not $\mbox{pot}(\bormone )$.\end{prop}

\vfill\eject

\noindent\bf Proof.\rm ~By Theorem \ref{potclo}, $\mathbb{G}_0$ and $s(\mathbb{G}_0)$ are not 
$\mbox{pot}(\bormone )$. As there is a rectangular continuous reduction of $\mathbb{G}_0$ or 
$s(\mathbb{G}_0)$ to the intersection of any of the other examples with $N_0\!\times\! N_1$, they are not $\mbox{pot}(\bormone )$. All the examples are $D_2(\boraone )$. They are clearly irreflexive, and have 
locally countable closure, like $\mathbb{G}_0$. We saw the acyclicity of 
$s(\overline{\mathbb{G}_0})$ in Proposition \ref{conn}, that of $s(\overline{\mathbb{B}_0})$ in Lemma \ref{propertiespi01}, and that of $s(\mathbb{T}_0)$ in the proof of Theorem \ref{generalneg}. The symmetrization of any of the ten sets is a subset of one of these three symmetrizations, and thus is acyclic.\bigskip

\noindent $\bullet$ By Proposition \ref{smallantipi01}, $\{\mathbb{B}_0,\mathbb{N}_0,\mathbb{M}_0,
s(\mathbb{B}_0)\}$ is an antichain.\bigskip

\noindent $\bullet$ As $\mathbb{U}_0$ is neither an oriented graph, nor a graph, it is not reducible to the other examples, except maybe $\mathbb{N}_0$. The set $\mathbb{U}_0$ is not reducible to $\mathbb{N}_0$ since $s(\mathbb{N}_0)$ is closed and $s(\mathbb{U}_0)$ is not.\bigskip

\noindent $\bullet$ As $\mathbb{N}_0$ is neither an oriented graph, nor a graph, it is not reducible to any of the other examples, except maybe $\mathbb{U}_0$. As its symmetrization is closed, the other examples different from $\mathbb{M}_0$ are not reducible to it.\bigskip

\noindent $\bullet$ Assume, towards a contradiction, that $\mathbb{N}_0$ is reducible to 
$\mathbb{U}_0$, with witness $f$. Then 
$$\big( f(0^\infty ),f(10^\infty )\big)\!\in\! (f\!\times\! f)[\overline{\mathbb{N}_0}
\!\setminus\!\mathbb{N}_0]\!\subseteq\!\overline{\mathbb{U}_0}\!\setminus\!\mathbb{U}_0
\mbox{,}$$ 
which gives $\varepsilon\!\in\! 2$ and $\beta\!\in\! 2^\omega$ such that 
$\big( f(0^\infty ),f(10^\infty )\big)\! =\! (\varepsilon\beta ,(1\! -\!\varepsilon )\beta )$. Thus 
$\big( (1\! -\!\varepsilon )\beta ,\varepsilon\beta )\!\in\!\mathbb{U}_0$, which is absurd. Note that this argument also shows that $\mathbb{M}_0$ is not reducible to $\mathbb{B}_0$ and $\mathbb{U}_0$.\bigskip

\noindent $\bullet$ As $\mathbb{G}_0,\mathbb{B}_0,\mathbb{M}_0,G_{s(\mathbb{G}_0)},
\mathbb{T}_0$ are oriented graphs and 
$s(\mathbb{G}_0),s(\mathbb{B}_0),s(\mathbb{T}_0)$ are graphs, the elements of the first set are incomparable with the elements of the second one. So we can consider these two sets separately.\bigskip

\noindent $\bullet$ Let us consider the first one. Note that 
$\mathbb{G}_0\not\sqsubseteq_c\mathbb{B}_0$ and 
$\mathbb{B}_0\not\sqsubseteq_c\mathbb{G}_0$. Indeed, for the first claim, there is a Borel countable coloring of $\mathbb{B}_0$. For the second one, we argue by contradiction, which gives $f$ continuous. As 
$(0^\infty ,10^\infty )\!\in\!\overline{\mathbb{B}_0}\!\setminus\!\mathbb{B}_0$, 
$\big( f(0^\infty ),f(10^\infty )\big)\!\in\!\overline{\mathbb{G}_0}\!\setminus\!\mathbb{G}_0\! =\!
\Delta (2^\omega )$, so that $f$ is not injective. Moreover, 
$\mathbb{B}_0\not\sqsubseteq_cG_{s(\mathbb{G}_0)}$, by Proposition \ref{not}.\bigskip

 Using the same arguments as in these proofs, we see that 
$\{\mathbb{G}_0,\mathbb{B}_0,G_{s(\mathbb{G}_0)}\}$ is an antichain, that $\mathbb{G}_0$ is incomparable with the other examples, and that $\mathbb{T}_0,\mathbb{M}_0$ are not reducible to $G_{s(\mathbb{G}_0)}$. The symmetrization of $\mathbb{M}_0$ is closed, which is not the case of the other symmetrizations, so that $\mathbb{M}_0$ cannot $\sqsubseteq_c$-reduce another one. Thus $\{\mathbb{G}_0,\mathbb{B}_0,\mathbb{M}_0,G_{s(\mathbb{G}_0)}\}$ is an antichain.\bigskip

 The set $\mathbb{T}_0$ is not reducible to $\mathbb{B}_0$. Indeed, we argue by contradiction, so that $\mathbb{T}_0$ is a subset of a $\mbox{pot}(\bormone )$ Acyclic oriented graph, by Theorem \ref{generalpos}. Thus $s(\mathbb{T}_0)$ is a subset of a $\mbox{pot}(\bormone )$ acyclic graph $G$, and $G_{s(\mathbb{G}_0)}$ ``is"  a subset of the $\mbox{pot}(\bormone )$ Acyclic oriented graph $G\cap (N_0\!\times\! N_1)$. By Theorem \ref{generalpos} again, 
$\mathbb{B}_0$ is reducible to $G_{s(\mathbb{G}_0)}$, which is absurd.\bigskip

 It remains to see that $\mathbb{B}_0,\mathbb{M}_0,G_{s(\mathbb{G}_0)}$ are not reducible to 
$\mathbb{T}_0$. If $\mathbb{B}_0$ is reductible to $\mathbb{T}_0$ with witness $f$, then $f$ is also a reduction of $s(\mathbb{B}_0)$ to $s(\mathbb{T}_0)$. As 
$(0^\infty ,10^\infty )\!\in\!\overline{\mathbb{B}_0}\!\setminus\!\mathbb{B}_0$, 
$$\big( f(0^\infty ),f(10^\infty )\big)\!\in\!\overline{\mathbb{T}_0}\!\setminus\!\mathbb{T}_0\! =\!
\{ (\varepsilon\gamma ,(1\! -\!\varepsilon )\gamma )\mid\varepsilon\!\in\! 2\ \wedge\ 
\gamma\!\in\! 2^\omega\} .$$ 

 This gives $\varepsilon\!\in\! 2$, $N\!\in\!\omega$ such that 
$f[N_{0^{N+1}}]\!\subseteq\! N_\varepsilon$ and $f[N_{10^N}]\!\subseteq\! N_{1-\varepsilon}$. Therefore the maps 
$\alpha\!\mapsto\! S\big( f(0\alpha )\big)$ and $\alpha\!\mapsto\! S\big( f(1\alpha )\big)$ define a rectangular continuous reduction of $\mathbb{G}_0\cap N_{0^N}^2$ to $s(\mathbb{G}_0)$. By Theorem \ref{G0}, $\mathbb{G}_0\sqsubseteq_c\mathbb{G}_0\cap N_{0^N}^2$. This gives a rectangular continuous reduction of $\mathbb{G}_0$ to $s(\mathbb{G}_0)$, which is absurd.\bigskip

 Similarly, $G_{s(\mathbb{G}_0)}$ is not reductible to $\mathbb{T}_0$. If 
$\mathbb{M}_0$ is reductible to $\mathbb{T}_0$ with witness $f$, then 
$$\big( f(0^\infty ),f(10^\infty )\big)\!\in\! (f\!\times\! f)[\overline{\mathbb{M}_0}\!\setminus\!
\mathbb{M}_0]\!\subseteq\!\overline{\mathbb{T}_0}\!\setminus\!\mathbb{T}_0\mbox{,}$$ 
which gives $\varepsilon\!\in\! 2$ and $\beta\!\in\! 2^\omega$ such that 
$\big( f(0^\infty ),f(10^\infty )\big)\! =\!\big(\varepsilon\beta ,(1\! -\!\varepsilon)\beta\big)$. Thus 
$\big( (1\! -\!\varepsilon)\beta ,\varepsilon\beta\big)\!\in\!\mathbb{T}_0$, which is absurd.\bigskip

\noindent $\bullet$ As $s(\mathbb{B}_0)$ is not reducible to $s(\mathbb{T}_0)$, 
$\mathbb{B}_0$ and $s(\mathbb{B}_0)$ are not reducible to $\mathbb{U}_0$. Let us prove that 
$G_{s(\mathbb{G}_0)},\mathbb{T}_0,s(\mathbb{T}_0)$ are not reducible to 
$\mathbb{U}_0$. Let us do it for $\mathbb{T}_0$, the other cases being similar.\bigskip

 We argue by contradiction, which gives $f$. Note that 
$\big( f(0^\infty ),f(10^\infty )\big)\!\in\!\overline{\mathbb{U}_0}\!\setminus\!\mathbb{U}_0$, which gives $\varepsilon\!\in\! 2$ and $N\!\in\!\omega$ such that 
$f[N_{0^{N+1}}]\!\subseteq\! N_\varepsilon$ and $f[N_{10^N}]\!\subseteq\! N_{1-\varepsilon}$. We can write $f(0^{N+1}\alpha )\! =\!\varepsilon g(\alpha )$, where $g$ is injective continuous. Similarly, $f(10^N\beta )\! =\! (1\! -\!\varepsilon )h(\beta )$, where $h$ is injective continuous. As 
$\big( f(0^{N+1}\alpha ),f(10^N\alpha )\big)\!\in\!\overline{\mathbb{U}_0}\!\setminus\!
\mathbb{U}_0$, $h\! =\! g$. Now ${(0^{N+1}\alpha ,10^N\beta )\!\in\!\mathbb{T}_0\Leftrightarrow
\big(\varepsilon g(\alpha ),(1\! -\!\varepsilon )g(\beta )\big)\!\in\!\mathbb{U}_0}$. Moreover, this implies that $(0^{N+1}\beta ,10^N\alpha )\!\notin\!\mathbb{T}_0$ and 
$\big(\varepsilon g(\beta ),(1\! -\!\varepsilon )g(\alpha )\big)\!\notin\!\mathbb{U}_0$, so that 
${\varepsilon\! =\! 1}$. Now $(10^N\alpha ,0^{N+1}\beta )\!\in\!\mathbb{T}_0$ is equivalent to 
$\big( 0g(\alpha ),1g(\beta )\big)\!\in\!\mathbb{U}_0$, 
$\big( 0g(\beta ),1g(\alpha )\big)\!\in\!\mathbb{U}_0$, and 
$$(10^N\beta ,0^{N+1}\alpha )\!\in\!\mathbb{T}_0\mbox{,}$$ 
which is absurd.\bigskip

\noindent $\bullet$ Let us consider the second one. As in the previous point, $s(\mathbb{G}_0)$ is not comparable with $s(\mathbb{B}_0)$ and $s(\mathbb{T}_0)$. We saw that $s(\mathbb{B}_0)$ is not reducible to $s(\mathbb{T}_0)$. If $s(\mathbb{T}_0)$ is reducible to $s(\mathbb{B}_0)$, then it is a subset of a $\mbox{pot}(\bormone )$ acyclic graph, which is absurd as before.
\hfill{$\square$}\bigskip

\noindent $\underline{\mbox{\bf A basis result}}$\bigskip

 We will see that the elements of this antichain are minimal. In fact, we prove more.\bigskip

\noindent\bf Proof of Theorem \ref{generalpos}.(5).(ii).\rm ~We set 
${\cal A}''\! :=\!\{\mathbb{B}_0,\mathbb{N}_0,\mathbb{M}_0\}$ and 
${\cal B}''\! :=\!\{s(\mathbb{B}_0)\}$. By Lemma \ref{propertiespi01} and Proposition \ref{biganti10}, ${\cal A}''\cup {\cal B}''$ is a $\sqsubseteq_c$-antichain made of $D_2(\boraone )$ relations, whose closure is Acyclic and is contained in $(N_0\!\times\! N_1)\cup (N_0\!\times\! N_1)$, which are not pot$(\bormone )$. By Lemma \ref{anti7}, $\cal A$ is also a $\sqsubseteq_c$-antichain. The proof of Proposition \ref{biganti10} shows that $ {\cal A}\cup\{\mathbb{G}_0,s(\mathbb{G}_0)\}$ is a 
$\sqsubseteq_c$-antichain, which is made of relations in the context of the theorem.\bigskip

\noindent $\bullet$ We first consider the case of digraphs. So assume that $B$ is a 
non-$\mbox{pot}(\bormone )$ Borel digraph contained in a $\mbox{pot}(\bormone )$ symmetric acyclic relation $F$. By Theorem \ref{G0}, we may assume that there is a Borel countable coloring 
$(B_n)_{n\in\omega}$ of $B$. We can change the Polish topology, so that we may assume that the $B_n$'s are clopen and $F$ is closed. Let $m\!\not=\! n$ such that 
$B\cap (B_m\!\times\! B_n)$ is not pot$(\bormone )$. Note that 
$F'\! :=\! F\cap\big( (B_m\!\times\! B_n)\cup (B_n\!\times\! B_m)\big)$ is a closed acyclic graph.\bigskip

 Corollary 3.10 in [L-Z] gives $f',g'$ injective continuous with $\mathbb{G}_0\! =\!
\overline{\mathbb{G}_0}\cap (f'\!\times\! g')^{-1}\big( B\cap (B_m\!\times\! B_n)\big)$. Lemma 
\ref{values} shows that $f'$ takes values in $B_m$ and $g'$ takes values in $B_n$, so that 
$$\mathbb{G}_0\! =\!\overline{\mathbb{G}_0}\cap (f'\!\times\! g')^{-1}(B).$$ 
The proof of Theorem \ref{generalpos}.(5).(i) gives $h\! :\! 2^\omega\!\rightarrow\! B_m\cup B_n$ injective continuous such that 
$${\mathbb{B}_0\! =\!\overline{\mathbb{B}_0}\cap (h\!\times\! h)^{-1}(B)}\mbox{,}$$ 
and $F''\! :=\! (h\!\times\! h)^{-1}(F')$ is a closed acyclic graph on $2^\omega$ contained in 
$(N_0\!\times\! N_1)\cup (N_1\!\times\! N_0)$ and containing $\mathbb{B}_0$. Theorem 
\ref{homB0} gives $i\! :\! 2^\omega\!\rightarrow\! 2^\omega$ injective continuous such that 
$\mathbb{B}_0\!\subseteq\! (i\!\times\! i)^{-1}(\mathbb{B}_0)$, 
$\overline{\mathbb{B}_0}\!\setminus\!\mathbb{B}_0\!\subseteq\! 
(i\!\times\! i)^{-1}(\overline{\mathbb{B}_0}\!\setminus\!\mathbb{B}_0)$ 
and $\neg s(\overline{\mathbb{B}_0})\!\subseteq\! (i\!\times\! i)^{-1}(\neg F'')$. We set 
$\tilde f\! :=\! h\circ i$, so that $\tilde f$ is injective continuous, 
$\mathbb{B}_0\!\subseteq\! B'\! :=\! (\tilde f\!\times\!\tilde f)^{-1}(B)$,  
$\overline{\mathbb{B}_0}\!\setminus\!\mathbb{B}_0\!\subseteq\!\neg B'$,  
$\neg s(\overline{\mathbb{B}_0})\!\subseteq\! (\tilde f\!\times\!\tilde f)^{-1}(\neg F')$, and thus 
$\neg s(\overline{\mathbb{B}_0})\!\subseteq\!\neg B'$. We proved that 
$\mathbb{B}_0\!\subseteq\! B'\!\subseteq\!\mathbb{N}_0$.\bigskip

\noindent\bf Case 1\rm\ $S\! :=\!\{\alpha\!\in\! 2^\omega\mid (1\alpha ,0\alpha )\!\in\! B'\}$ is meager.\bigskip

 Then $(2^\omega ,A)\sqsubseteq_c(X,B)$ for some $A\!\in\!\{\mathbb{B}_0,s(\mathbb{B}_0)\}$. Indeed, let 
$\tilde G$ be a dense $G_\delta$ subset of $2^\omega$ disjoint from $S$, and $G\! :=\! 2\!\times\!\tilde G$. Then 
$$\mathbb{B}_0\cap G^2\!\subseteq\! B'\cap G^2\!\subseteq\! s(\mathbb{B}_0)\cap G^2.$$
We set $B''\! :=\!\{ (\alpha ,\beta )\!\in\! {\tilde G}^2\mid (1\alpha ,0\beta )\!\in\! B'\}$. Note that $B''$ is a Borel oriented graph on $\tilde G$ contained in the $F_\sigma$ acyclic graph 
$s(\mathbb{G}_0)\cap {\tilde G}^2$. By Theorem \ref{G0}, either $B''$ has a Borel countable coloring, or $(2^\omega ,\mathbb{G}_0^{-1})\sqsubseteq_c(\tilde G,B'')$ with witness $g$.\bigskip

\noindent - In the first case, we find a non meager $G_\delta$ subset $G'$ of $2^\omega$ contained in 
$\tilde G$ which is $B''$-discrete. Note that 
$B'\cap (2\!\times\! G')^2\! =\!\mathbb{B}_0\cap (2\!\times\! G')^2$ and 
$(2^\omega ,\mathbb{B}_0)\sqsubseteq_c
\big( 2\!\times\! G',\mathbb{B}_0\cap (2\!\times\! G')^2\big)\sqsubseteq_c(X,B)$, by Theorem 
\ref{generalpos} and Lemma \ref{propertiespi01}.\bigskip

\noindent - In the second case, note that $\mathbb{G}_0\!\subseteq\! (g\!\times\! g)^{-1}({B''}^{-1})
\!\subseteq\! (g\!\times\! g)^{-1}(\mathbb{G}_0)$. Theorem \ref{Fsigma} gives 
${g''\! :\! 2^\omega\!\rightarrow\! 2^\omega}$ injective continuous such that 
$$\mathbb{G}_0\!\subseteq\! (g''\!\times\! g'')^{-1}(\mathbb{G}_0)\!\subseteq\! 
(g''\!\times\! g'')^{-1}\big( (g\!\times\! g)^{-1}(\mathbb{G}_0)\big)\!\subseteq\! s(\mathbb{G}_0).$$ 
As $(g\!\times\! g)^{-1}(\mathbb{G}_0)$ is an oriented graph, 
$\mathbb{G}_0\! =\! (g''\!\times\! g'')^{-1}\big( (g\!\times\! g)^{-1}(\mathbb{G}_0)\big)$. We set 
$$f''(\varepsilon\alpha )\! :=\!\varepsilon g\big( g''(\alpha )\big)\mbox{,}$$ 
so that $f''$ is injective continuous. If ${(0\alpha ,1\beta )\!\in\! s(\mathbb{B}_0)}$, then 
$\big( g''(\alpha ),g''(\beta )\big)\!\in\!\mathbb{G}_0$, 
$$\big( (gg'')(\alpha ),(gg'')(\beta )\big)\!\in\! {B''}^{-1}$$ 
and $\big( 1(gg'')(\beta ),0(gg'')(\alpha )\big)\!\in\! B'\cap G^2$. Thus 
$\big( 1(gg'')(\beta ),0(gg'')(\alpha )\big)\!\in\! s(\mathbb{B}_0)$, 
$${\big( 1(gg'')(\beta ),0(gg'')(\alpha )\big)\!\in\!\mathbb{B}_0^{-1}}\mbox{,}$$ 
${\big( 0(gg'')(\alpha ),1(gg'')(\beta )\big)\!\in\!\mathbb{B}_0\cap G^2}$ and 
$\big( f''(0\alpha ),f''(1\beta )\big)\!\in\! B'$. In particular, if 
${(1\alpha ,0\beta )\!\in\! s(\mathbb{B}_0)}$, then $\big( f''(1\alpha ),f''(0\beta )\big)\!\in\! B'$.

\vfill\eject

 Conversely, if $\big( f''(0\alpha ),f''(1\beta )\big)\!\in\! B'$, then 
$\big( f''(0\alpha ),f''(1\beta )\big)\!\in\!\mathbb{B}_0$, 
${\big( (gg'')(\alpha ),(gg'')(\beta )\big)\!\in\!\mathbb{G}_0}$, $(\alpha ,\beta )\!\in\!\mathbb{G}_0$ and 
$(0\alpha ,1\beta )\!\in\! s(\mathbb{B}_0)$. If $\big( f''(1\alpha ),f''(0\beta )\big)\!\in\! B'$, then 
$\big( (gg'')(\alpha ),(gg'')(\beta )\big)\!\in\! B''$, 
$\big( g''(\beta ),g''(\alpha )\big)\!\in\!\mathbb{G}_0$, $(0\beta ,1\alpha )\!\in\! s(\mathbb{B}_0)$ and $(1\alpha ,0\beta )\!\in\! s(\mathbb{B}_0)$. Thus $f''$ is a witness for 
$$\big( 2^\omega ,s(\mathbb{B}_0)\big)\sqsubseteq_c(2^\omega ,B')\sqsubseteq_c(X,B).$$
\bf Case 2\rm\ $S$ is not meager.\bigskip

 Then let us show that $(2^\omega ,A)\sqsubseteq_c(X,B)$ for some 
$A\!\in\!\{\mathbb{N}_0,\mathbb{M}_0\}$. Indeed, let $\tilde H$ be a non-meager $G_\delta$ subset of 
$2^\omega$ contained in $S$, and $H\! :=\! 2\!\times\!\tilde H$. Then 
$\mathbb{M}_0\cap H^2\!\subseteq\! B'\cap H^2\!\subseteq\!\mathbb{N}_0\cap H^2$. We set 
$B''\! :=\!\{ (\alpha ,\beta )\!\in\! {\tilde H}^2\mid\alpha\!\not=\!\beta\ \wedge\ (1\alpha ,0\beta )\!\in\! B'\}$. Note that $B''\!\subseteq\!\mathbb{G}_0^{-1}$ is Borel. By Theorem \ref{G0}, either there is a Borel countable coloring of $B''$, or $(2^\omega ,\mathbb{G}_0^{-1})\sqsubseteq_c(\tilde H,B'')$ with witness $g$.\bigskip

\noindent - In the first case, there is a non meager $G_\delta$ subset $H'$ of $2^\omega$ contained in 
$\tilde H$ which is $B''$-discrete. Note that 
$B'\cap (2\!\times\! H')^2\! =\!\mathbb{M}_0\cap (2\!\times\! H')^2$ and 
$\big(2\!\times\! H',\mathbb{M}_0\cap (2\!\times\! H')^2\big)\sqsubseteq_c(2^\omega ,B')\sqsubseteq_c(X,B)$. So we are done if we prove that $(2^\omega ,\mathbb{M}_0)\sqsubseteq_c\big(2\!\times\! H',\mathbb{M}_0\cap (2\!\times\! H')^2\big)$. Note that $\mathbb{G}_0\cap {H'}^2$ is a $F_\sigma$ Acyclic oriented graph on $H'$ without Borel countable coloring. Theorem \ref{G0} gives $\tilde g\! :\! 2^\omega\!\rightarrow\! H'$ injective continuous such that 
$\mathbb{G}_0\! =\! (\tilde g\!\times\!\tilde g)^{-1}(\mathbb{G}_0)$. It remains to consider 
$\varepsilon\alpha\!\mapsto\!\varepsilon\tilde g(\alpha )$ to get our reduction.\bigskip

\noindent - In the second case, note that 
$$\mathbb{G}_0\!\subseteq\!
\{ (\alpha ,\beta )\!\in\! 2^\omega\!\times\! 2^\omega\mid\alpha\!\not=\!\beta\ \wedge\ 
\big( 1g(\beta ),0g(\alpha )\big)\!\in\! B'\}\!\subseteq\! (g\!\times\! g)^{-1}(\mathbb{G}_0).$$ 
Theorem \ref{Fsigma} gives $g''\! :\! 2^\omega\!\rightarrow\! 2^\omega$ injective continuous such that 
$$\mathbb{G}_0\!\subseteq\! (g''\!\times\! g'')^{-1}(\mathbb{G}_0)\!\subseteq\! 
(g''\!\times\! g'')^{-1}\big( (g\!\times\! g)^{-1}(\mathbb{G}_0)\big)\!\subseteq\! 
s(\mathbb{G}_0)\mbox{,}$$ 
and $\mathbb{G}_0\! =\! (g'\!\times\! g')^{-1}\big( (g\!\times\! g)^{-1}(\mathbb{G}_0)\big)$ as in the previous point. We define $f''$ as in the previous point, and here again $f''$ is a witness for 
$(2^\omega ,\mathbb{N}_0)\sqsubseteq_c(2^\omega ,B')\sqsubseteq_c(X,B)$.\bigskip

\noindent $\bullet$ We now consider the general case of non necessarily irreflexive relations. Let $B$ be a non-pot$(\bormone )$ Borel subset of a pot$(\bormone )$ symmetric acyclic relation. We may assume that $B$ is contained in a closed symmetric acyclic relation $F$.\bigskip

\noindent\bf Case 1\rm\ If $C,D$ are disjoint Borel subsets of $X$, then $B\cap (C\!\times\! D)$ is 
pot$(\bormone )$.\bigskip

 We set $N\! :=\!\{ x\!\in\! X\mid (x,x)\!\notin\! B\}$. Note that $N$ is Borel, so that we may assume that $N$ is clopen, and $B\cap (N\!\times\!\neg N)$ and $B\cap (\neg N\!\times\! N)$ are 
pot$(\bormone )$. This is also the case of $B\cap (\neg N)^2$, which is a reflexive relation on 
$\neg N$. Indeed, we may assume that $\neg N$ is uncountable, which gives a Borel isomorphism 
$\Psi\! :\! 2^\omega\!\rightarrow\!\neg N$. Note that $(\neg N)^2\!\setminus\!\Delta (\neg N)\! =\!
\bigcup_{s\in 2^{<\omega },\varepsilon\in 2}~
\Psi [N_{s\varepsilon}]\!\times\!\Psi [N_{s(1-\varepsilon )}]$, so that 
$$B\cap (\neg N)^2\! =\!\Delta (\neg N)\cup\bigcup_{s\in 2^{<\omega },\varepsilon\in 2}~
B\cap (\Psi [N_{s\varepsilon}]\!\times\!\Psi [N_{s(1-\varepsilon )}])$$ 
and $(\Psi\!\times\!\Psi)^{-1}\big( B\cap (\neg N)^2\big)\! =\!\Delta (2^\omega )\cup
\bigcup_{s\in 2^{<\omega },\varepsilon\in 2}~
(\Psi\!\times\!\Psi)^{-1}(B)\cap (N_{s\varepsilon}\!\times\! N_{s(1-\varepsilon )})$. By our assumption, the $(\Psi\!\times\!\Psi)^{-1}(B)\cap (N_{s\varepsilon}\!\times\! N_{s(1-\varepsilon )})$'s are pot$(\bormone )$. We are done since they can accumulate only on the diagonal. This shows that $B\cap N^2$ is a non-pot$(\bormone )$ Borel digraph on $N$. By our assumption, it has no Borel countable coloring. Theorem \ref{G0} gives $A\!\in\!\{\mathbb{G}_0,s(\mathbb{G}_0)\}$ such that 
$$(2^\omega ,A)\sqsubseteq_c(N,B\cap N^2)\sqsubseteq_c(X,B).$$ 

\vfill\eject

\noindent\bf Case 2\rm\ There are disjoint Borel subsets $C,D$ of $X$ such that $B\cap (C\!\times\! D)$ is not 
pot$(\bormone )$.\bigskip

 Note that we may assume that $C,D$ are clopen. The case of digraphs gives 
$$A\!\in\!\{\mathbb{B}_0,\mathbb{N}_0,\mathbb{M}_0,s(\mathbb{B}_0)\}$$ 
such that 
$(2^\omega ,A)\sqsubseteq_c\Big( C\cup D,B\cap\big( (C\!\times\! D)\cup (D\!\times\! C)\big)\Big)$ with witness $g$, for coloring reasons. Note that 
$\big( g(0^\infty ),g(10^\infty )\big)\!\in\! (g\!\times\! g)[\overline{A}]$, so that 
$\big( g(0^\infty ),g(10^\infty )\big)\!\in\! C\!\times\! D$, for example. The continuity of $g$ gives 
$N\!\in\!\omega$ such that $g[N_{0^{N+1}}]\!\subseteq\! C$ and $g[N_{10^N}]\!\subseteq\! D$. Note that 
$$(g\!\times\! g)^{-1}(B)\cap\big( (N_{0^{N+1}}\!\times\! N_{10^N})\cup (N_{10^N}\!\times\! N_{0^{N+1}})\big)\! =\! A\cap (N_{0^{N+1}}\cup N_{10^N})^2$$
is not pot$(\bormone )$. This implies that we may assume that $X\! =\! N_{0^{N+1}}\cup N_{10^N}$ and 
$$B\cap\big( (N_{0^{N+1}}\!\times\! N_{10^N})\cup (N_{10^N}\!\times\! N_{0^{N+1}})\big)\! =\! A\cap (N_{0^{N+1}}\cup N_{10^N})^2.$$ 
We set $F'\! :=\!\{ (\alpha ,\beta )\!\in\! N_{0^N}^2\mid\exists\varepsilon ,\varepsilon'\!\in\! 2~~(\varepsilon\alpha ,\varepsilon'\beta )\!\in\! F\}$. Note that $F'$ is a closed symmetric relation on $N_{0^N}$ containing $\mathbb{G}_0\cap N_{0^N}^2$. Moreover, $F'$ is acyclic. Indeed, we argue by contradiction to see this, which gives $n\!\geq\! 2$ and $(\gamma_i)_{i\leq n}$ injective with $(\gamma_i,\gamma_{i+1})\!\in\! F'$ for each $i\! <\! n$ and 
$(\gamma_0,\gamma_n)\!\in\! F'$. This provides $(\varepsilon_j)_{j\leq 2n+1}\!\in\! 2^{2n+2}$ such that $(\varepsilon_{2i}\gamma_i,\varepsilon_{2i+1}\gamma_{i+1})\!\in\! F$ if $i\! <\! n$ and 
$(\varepsilon_{2n}\gamma_0,\varepsilon_{2n+1}\gamma_n)$ is in $F$. If 
$\varepsilon_1\!\not=\!\varepsilon_2$, then 
$(\varepsilon_1\gamma_1,\varepsilon_2\gamma_1)\!\in\! s(\overline{B})\!\subseteq\! F$. This gives an injective $F$-path with at least $n\! +\! 1$ elements contradicting the acyclicity of $F$.\bigskip

 Corollary \ref{corFsigma} gives $h\! :\! 2^\omega\!\rightarrow\! N_{0^N}$ injective continuous such that 
$$\mathbb{G}_0\!\subseteq\! (h\!\times\! h)^{-1}\big(\mathbb{G}_0\cap N_{0^N}^2\big)
\!\subseteq\! (h\!\times\! h)^{-1}\big( F'\!\setminus\!\Delta (X)\big)\!\subseteq\! s(\mathbb{G}_0).$$ Symmetry considerations show that in fact 
$$\mathbb{G}_0\! =\! (h\!\times\! h)^{-1}\big(\mathbb{G}_0\cap N_{0^N}^2\big)\!\subseteq\! 
(h\!\times\! h)^{-1}\big( F'\!\setminus\!\Delta (X)\big)\! =\! s(\mathbb{G}_0).$$ 
We set $k(\varepsilon\alpha )\! :=\!\varepsilon h(\alpha )$, which defines 
$k\! :\! 2^\omega\!\rightarrow\! N_{0^{N+1}}\cup N_{10^N}$ injective continuous with 
$$\mathbb{B}_0\!\subseteq\! (k\!\times\! k)^{-1}(B)\!\subseteq\!\{ (\varepsilon\alpha ,\varepsilon'\beta )\!\in\! 2^\omega\!\times\! 2^\omega\mid 
(\alpha ,\beta )\!\in\!\overline{s(\mathbb{G}_0)}\}\!\setminus\!\{ (0\gamma ,1\gamma )\mid\gamma\!\in\! 2^\omega\} .$$ 
This means that we may assume that $X\! =\! 2^\omega$ and 
$$\mathbb{B}_0\!\subseteq\! B\!\subseteq\!
\{ (\varepsilon\alpha ,\varepsilon'\beta )\!\in\! 2^\omega\!\times\! 2^\omega\mid (\alpha ,\beta )\!\in\!\overline{s(\mathbb{G}_0)}\}
\!\setminus\!\{ (0\gamma ,1\gamma )\mid\gamma\!\in\! 2^\omega\} .$$ 
This proof also shows that we may assume that 
$B\cap\big( (N_0\!\times\! N_1)\cup (N_0\!\times\! N_1)\big)\! =\! A$. It remains to study 
$B\cap (N_0^2\cup N_1^2)$. Assume that $(0\alpha ,0\beta )\!\in\! B$ and $\alpha\!\not=\!\beta$. Then we can find $n\!\in\!\omega$, $\varepsilon\!\in\! 2$ and 
$\gamma\!\in\! 2^\omega$ such that 
$(\alpha ,\beta )\!\in\! (s_n\varepsilon\gamma ,s_n(1\! -\!\varepsilon )\gamma )$. Then 
$(0s_n1\gamma ,0s_n0\gamma ,1s_n1\gamma )$ is an injective $F$-path contradicting the acyclicity of $F$ since 
$(0s_n1\gamma ,1s_n1\gamma )\!\in\!\overline{\mathbb{B}_0}\!\subseteq\! F$. Similarly, $(1\alpha ,1\beta )$ cannot be in $B$ if $\alpha\!\not=\!\beta$. This proves that we may assume that $A\!\subseteq\! B\!\subseteq\! A\cup\Delta (2^\omega )$. This means that we may assume that 
$X\! =\! 2^\omega$ and there is a Borel subset $I$ of $2^\omega$ such that $B\! =\! A\cup\Delta (I)$. Then we argue as in the proof of Theorem 
\ref{baspi0>1}. The (a) part of the claim comes from Lemma \ref{propertiespi01}.(1). For the (b) part of the claim, the minimality of $\mathbb{B}_0$ comes from the case of digraphs. The witness $f$ is a homomorphism from $\overline{\mathbb{B}_0}$ into itself and sends $N_\varepsilon$ into itself. For  Case 2, $\mathbb{B}_0$ is minimum among non-pot$(\bormone )$ subsets of a closed Acyclic oriented graph, by Theorem \ref{generalpos}.(5).(i), so that we can apply Proposition 
\ref{exclude2}. We conclude as in the proof of Theorem \ref{baspi0>1}.\hfill{$\square$}

\vfill\eject
 
 Note that this implies that $\mathbb{G}_0,\mathbb{B}_0,\mathbb{N}_0,\mathbb{M}_0,
s(\mathbb{G}_0),s(\mathbb{B}_0)$ are $\sqsubseteq_c$-minimal among non-pot$(\bormone )$ relations. Theorem \ref{generalpos}.(5).(ii) is optimal in terms of potential complexity, because of 
$G_{s(\mathbb{G}_0)}$, $\mathbb{T}_0$ and $s(\mathbb{T}_0)$, by Proposition \ref{biganti10}. We now give a consequence of our results of injective reduction on a closed set.
 
\begin{cor} \label{situation} Let $X$ be a Polish space, and $B$ be a Borel Acyclic digraph on 
$X$ contained in a $\mbox{pot}\big(\check D_2(\boraone )\big)$ locally countable relation. Then exactly one of the following holds:\smallskip  

(a) the set $B$ is $\mbox{pot}(\bormone )$,\smallskip  

(b) $(2^\omega ,\mathbb{G}_0)\sqsubseteq_c(X,B)$ or 
$\big( 2^\omega ,s(\mathbb{G}_0)\big)\sqsubseteq_c(X,B)$ or there is a $F_\sigma$ Acyclic digraph $B'$ on $2^\omega$ with locally countable closure contained in 
$(N_0\!\times\! N_1)\cup (N_1\!\times\! N_0)$ such that 
$\mathbb{B}_0\! =\!\overline{\mathbb{B}_0}\cap B'$ and $(2^\omega ,B')\sqsubseteq_c(X,B)$.
\end{cor}

\noindent\bf Proof.\rm ~Let $F$ be a $\mbox{pot}\big(\check D_2(\boraone )\big)$ locally countable relation containing $B$. Then $F$ is in fact 
$\mbox{pot}(\bormone )$. Assume that (a) does not hold. By Theorem \ref{G0}, we may assume that there is a Borel coloring $(B_n)_{n\in\omega}$ of $B$. As $B$ is locally countable and we can change the Polish topology, we may assume that $B$ is $F_\sigma$, $F$ is closed, and the $B_n$'s are clopen. Let 
$m\!\not=\! n$ such that $B\cap (B_m\!\times\! B_n)$ is not $\mbox{pot}(\bormone )$. Corollary 3.10 in [L-Z] gives $f_0\! :\! 2^\omega\!\rightarrow\! B_m$ and 
$f_1\! :\! 2^\omega\!\rightarrow\! B_n$ injective continuous with 
$\mathbb{G}_0\! =\!\overline{\mathbb{G}_0}\cap (f_0\!\times\! f_1)^{-1}\big( B\cap (B_m\!\times\! B_n)\big)$. We set 
$h(\varepsilon\alpha )\! :=\! f_\varepsilon (\alpha )$, so that $h$ is injective continuous and 
$\mathbb{B}_0\! =\!\overline{\mathbb{B}_0}\cap (h\!\times\! h)^{-1}(B)$. It remains to set $B'\! :=\! (h\!\times\! h)^{-1}(B)$.\hfill{$\square$}\bigskip

\noindent $\underline{\mbox{\bf Minimality}}$
 
\begin{thm} \label{minGsG0} The sets $G_{s(\mathbb{G}_0)},s(\mathbb{T}_0)$ are 
$\sqsubseteq_c$-minimal among non-pot$(\bormone )$ relations.\end{thm}

\noindent\bf Proof.\rm ~By Lemma \ref{minimalitysymmetrization}, it is enough to prove that 
$G_{s(\mathbb{G}_0)}$ is $\sqsubseteq_c$-minimal among non-pot$(\bormone )$ relations. So assume that $A\!\subseteq\! X^2$ is not pot$(\bormone )$ and $(X,A)\sqsubseteq_c(2^\omega ,G_{s(\mathbb{G}_0)})$ with witness $g$. Then $A$ is a $D_2(\boraone )$ Acyclic oriented graph with locally countable closure. By Corollary \ref{situation}, 
$(2^\omega ,\mathbb{G}_0)\sqsubseteq_c(X,A)$ or there is a $D_2(\boraone )$ Acyclic oriented graph $B$ on $2^\omega$ with locally countable closure contained in 
$(N_0\!\times\! N_1)\cup (N_1\!\times\! N_0)$ such that 
$\mathbb{B}_0\! =\!\overline{\mathbb{B}_0}\cap B$ and $(2^\omega ,B)\sqsubseteq_c(X,A)$. So we may assume that $X$ is compact and $A\!\in\! K_\sigma$. We set $R\! :=\! (g\!\times\! g)[A]$, so that $R\!\subseteq\! G_{s(\mathbb{G}_0)}\cap g[X]^2$ is $K_\sigma$, 
$(X,A)\sqsubseteq_c(g[X],R)$ and $(g[X],R)\sqsubseteq_c(X,A)$. In particular, $R$ is not 
pot$(\bormone )$. We set 
$$B\! :=\!\{ (\alpha ,\beta )\!\in\! 2^\omega\!\times\! 2^\omega\mid (0\alpha ,1\beta )\!\in\! R\} .$$ 
Note that $B\!\subseteq\! s(\mathbb{G}_0)$ and the shift map is a rectangular reduction of 
$R\!\subseteq\! N_0\!\times\! N_1$ to $B$. Thus $B$ is a non-pot$(\bormone )$ subset of 
$s(\mathbb{G}_0)$. This implies that $B$ has no Borel countable coloring. Theorem \ref{G0} implies that $\mathbb{G}_0\sqsubseteq_cB$ or $s(\mathbb{G}_0)\sqsubseteq_cB$, with witness $h$. We set 
$f(\varepsilon\alpha )\! :=\!\varepsilon h(\alpha )$, so that $f$ is injective continuous. As 
$R\!\subseteq\! N_0\!\times\! N_1$, we get $(2^\omega ,\mathbb{B}_0)\sqsubseteq_c(g[X],R)$ or 
$(2^\omega ,G_{s(\mathbb{G}_0)})\sqsubseteq_c(g[X],R)$. The first possibility cannot occur by Proposition \ref{biganti10} and we are done.\hfill{$\square$}\bigskip

 We need several results to prepare the proof of the minimality of $\mathbb{T}_0$ and 
$\mathbb{U}_0$.

\begin{thm} \label{graphsT0} Let $B\!\subseteq\! (N_0\!\times\! N_1)\cup (N_1\!\times\! N_0)$ be a $F_\sigma$ acyclic graph on $2^\omega$ such that 
$\mathbb{B}_0\! =\!\overline{\mathbb{B}_0}\cap B$. We assume that $(0\alpha ,1\beta )\!\in\! B\Leftrightarrow (0\beta ,1\alpha )\!\in\! B$. Then there is $f\! :\! 2^\omega\!\rightarrow\! 2^\omega$ injective continuous such that $s(\mathbb{T}_0)\! =\! (f\!\times\! f)^{-1}(B)$ and 
$\mathbb{B}_0\!\subseteq\! (f\!\times\! f)^{-1}(\mathbb{B}_0)$.\end{thm}

\noindent\bf Proof.\rm ~We set $B'\! :=\!\{ (\alpha ,\beta )\!\in\! 2^\omega\!\times\! 2^\omega\mid (0\alpha ,1\beta )\!\in\! B\}$. Note that $B'$ is a $F_\sigma$ acyclic graph on $2^\omega$.\bigskip

 Indeed, assume that $n\!\geq\! 2$ and $(x_i)_{i\leq n}$ is an injective $B'$-path with $(x_0,x_n)\!\in\! B'$. If 
$n$ is odd, then $(0x_0,1x_1,0x_2,1x_3,...,1x_n)$ is an injective $B$-path contradicting the acyclicity of $B$. If $n$ is even, then $(0x_0,1x_1,0x_2,1x_3,...,0x_n,1x_0,0x_1,1x_2,0x_3,...,1x_n)$ is also an injective $B$-path contradicting the acyclicity of $B$.\bigskip

 Theorem \ref{Fsigma} gives $g\! :\! 2^\omega\!\rightarrow\! 2^\omega$ injective continuous satisfying the inclusions $s(\mathbb{G}_0)\! =\! (g\!\times\! g)^{-1}(B')$ and 
$\mathbb{G}_0\!\subseteq\! (g\!\times\! g)^{-1}(\mathbb{G}_0)$ since $\mathbb{G}_0\!\subseteq\! B'$. We set $f(\varepsilon\alpha )\! :=\!\varepsilon g(\alpha )$, so that $f$ is injective continuous. Note that 
$\mathbb{B}_0\!\subseteq\! (f\!\times\! f)^{-1}(\mathbb{B}_0)$ and 
$$(0\alpha ,1\beta )\!\in\! G_{s(\mathbb{G}_0)}\!\Leftrightarrow\! (\alpha ,\beta )\!\in\! s(\mathbb{G}_0)\!\Leftrightarrow\!\big( g(\alpha ),g(\beta )\big)\!\in\! B'
\!\Leftrightarrow\!\big( f(0\alpha ),f(1\beta )\big)\! =\!\big( 0g(\alpha ),1g(\beta )\big)\!\in\! B$$  
and we are done since $B\!\subseteq\! (N_0\!\times\! N_1)\cup (N_1\!\times\! N_0)$ is symmetric.\hfill{$\square$}

\begin{thm} \label{ogasGsG0} Let $B\!\subseteq\! (N_0\!\times\! N_1)\cup (N_1\!\times\! N_0)$ be a $F_\sigma$ Acyclic oriented graph on $2^\omega$ such that 
$\mathbb{B}_0\! =\!\overline{\mathbb{B}_0}\cap B$. We assume that 
$\big(\varepsilon\alpha ,(1\! -\!\varepsilon )\beta\big)\!\in\! B\Leftrightarrow\big(\varepsilon\beta ,(1\! -\!\varepsilon )\alpha\big)\!\in\! B$. Then 
$(2^\omega ,G_{s(\mathbb{G}_0)})\sqsubseteq_c(2^\omega ,B)$.\end{thm}

\noindent\bf Proof.\rm ~We set $B'\! :=\!\{ (\alpha ,\beta )\!\in\! 2^\omega\!\times\! 2^\omega\mid 
(0\alpha ,1\beta )\!\in\! B\}$. Note that $B'$ is a $F_\sigma$ acyclic graph on $2^\omega$, as in the proof of Theorem \ref{graphsT0}. We set $M\! :=\!\{\alpha\!\in\! 2^\omega\mid (1\alpha ,0\alpha )\!\in\! B\}$. Note that $M$ is meager. Indeed, we argue by contradiction. As $M$ is $F_\sigma$, this gives 
$(\alpha ,\beta )\!\in\!\mathbb{G}_0\cap M^2$. Then by assumption $(0\alpha ,1\beta ,0\beta ,1\alpha )$ is an injective $s(B)$-path contradicting the Acyclicity of $B$. So let $G$ be a dense $G_\delta$ subset of 
$2^\omega$ disjoint from $M$. Note that $\mathbb{G}_0\!\subseteq\! B'$, so that there is no Borel countable coloring of $B'\cap G^2$. Theorem \ref{G0} gives $g\! :\! 2^\omega\!\rightarrow\! G$  injective continuous such that $s(\mathbb{G}_0)\! =\! (g\!\times\! g)^{-1}(B')$. We set 
$f(\varepsilon\alpha )\! :=\!\varepsilon g(\alpha )$, so that $f$ is injective continuous. Note that 
$(0\alpha ,1\beta )\!\in\! G_{s(\mathbb{G}_0)}\!\Leftrightarrow\!\big( f(0\alpha ),f(1\beta )\big)\!\in\! B$, as in the proof of Theorem \ref{graphsT0}. We set 
$B''\! :=\! (f\!\times\! f)^{-1}(B)$. Then $B''$ satisfies the same assumptions as $B$, $B''\sqsubseteq_cB$ and $(1\alpha ,0\alpha )\!\notin\! B''$ for each 
$\alpha\!\in\! 2^\omega$.\bigskip

 We now set $B'''\! :=\!\{ (\alpha ,\beta )\!\in\! 2^\omega\!\times\! 2^\omega\mid 
(0\alpha ,1\beta )\!\in\! B''\ \vee\ (1\alpha ,0\beta )\!\in\! B''\}$. As for $B'$, Theorem \ref{G0} gives 
$h\! :\! 2^\omega\!\rightarrow\! 2^\omega$ injective continuous such that 
$s(\mathbb{G}_0)\! =\! (h\!\times\! h)^{-1}(B''')$. We set $l(\varepsilon\alpha )\! :=\!\varepsilon h(\alpha )$, so that $l$ is injective continuous. As in the previous paragraph, 
$(0\alpha ,1\beta )\!\in\! G_{s(\mathbb{G}_0)}\!\Leftrightarrow\!\big( l(0\alpha ),l(1\beta )\big)\!\in\! B''$. It remains to see that $\big( l(1\alpha ),l(0\beta )\big)\!\notin\! B''$. We argue by contradiction, so that 
$\big( 1h(\alpha ),0h(\beta )\big)\!\in\! B''$, $\big( h(\alpha ),h(\beta )\big)\!\in\! B'''$, 
$(\alpha ,\beta )\!\in\! s(\mathbb{G}_0)$, $(\beta ,\alpha )\!\in\! s(\mathbb{G}_0)$, 
$(0\beta ,1\alpha )\!\in\! G_{s(\mathbb{G}_0)}$ and $\big( l(0\beta ),l(1\alpha )\big)\!\in\! B''$. As $B''$ is antisymmetric, we get $l(0\beta )\! =\! l(1\alpha )$, which is absurd.\hfill{$\square$}

\begin{thm} \label{insidesT0} Let $B$ be a Borel relation on $2^\omega$ such that $s(B)\! =\! s(\mathbb{T}_0)$. Then $(2^\omega ,A)\sqsubseteq_c(2^\omega ,B)$ for some $A$ in $\{ G_{s(\mathbb{G}_0)},\mathbb{T}_0,\mathbb{U}_0,s(\mathbb{T}_0)\}$.\end{thm}

\noindent\bf Proof.\rm ~We set $B'\! :=\!\{ (\alpha ,\beta )\!\in\! 2^\omega\!\times\! 2^\omega\mid 
(0\alpha ,1\beta )\!\in\! B\}$. Note that $B'\!\subseteq\! s(\mathbb{G}_0)$.\bigskip

\noindent $\bullet$ Let us prove that we may assume that $B'$ is not pot$(\bormone )$. We argue by contradiction, which gives a non-meager $G_\delta$ subset $G$ of $2^\omega$ which is $B'$-discrete. Note that $B\cap (2\!\times\! G)^2\! =\! G_{s(\mathbb{G}_0)}^{-1}\cap (2\!\times\! G)^2$ since 
$B\!\subseteq\! s(B)\! =\! s(\mathbb{T}_0)$. As $G$ is not meager, there is no Borel countable coloring of 
$\big( G,s(\mathbb{G}_0)\cap G^2\big)$. Theorem \ref{G0} shows that 
$\big( 2^\omega ,s(\mathbb{G}_0)\big)\sqsubseteq_c\big( G,s(\mathbb{G}_0)\cap G^2\big)$ with witness $g$. The map $\varepsilon\alpha\!\mapsto\!\varepsilon g(\alpha )$ is a witness for 
$(2^\omega ,G_{s(\mathbb{G}_0)}^{-1})\sqsubseteq_c
\big( 2\!\times\! G,G_{s(\mathbb{G}_0)}^{-1}\cap (2\!\times\! G)^2\big)$. As 
$(2^\omega ,G_{s(\mathbb{G}_0)})\sqsubseteq_c(2^\omega ,G_{s(\mathbb{G}_0)}^{-1})$ with witness 
$\varepsilon\alpha\!\mapsto\! (1\! -\!\varepsilon )\alpha$, 
$(2^\omega ,G_{s(\mathbb{G}_0)})\sqsubseteq_c(2^\omega ,B)$.\bigskip

\noindent $\bullet$ Theorem \ref{G0} implies that 
$( 2^\omega ,\mathbb{G}_0)\sqsubseteq_c( 2^\omega ,B')$ or 
$\big( 2^\omega ,s(\mathbb{G}_0)\big)\sqsubseteq_c( 2^\omega ,B')$ with witness $h$.\bigskip

\noindent\bf Case 1.1\rm\ $( 2^\omega ,\mathbb{G}_0)\sqsubseteq_c( 2^\omega ,B')$\bigskip

 In this case, 
$(0\alpha ,1\beta )\!\in\!\mathbb{T}_0\Leftrightarrow\big( 0h(\alpha ),1h(\beta )\big)\!\in\! B$. As 
${B\!\subseteq\! s(B)\! =\! s(\mathbb{T}_0)}$, 
$$\big( 0h(\alpha ),1h(\beta )\big)\!\in\! B
\Rightarrow\big( 0h(\beta ),1h(\alpha )\big)\!\in\! s(B)\!\setminus\! B\! =\! B^{-1}
\Rightarrow\big( 1h(\alpha ),0h(\beta )\big)\!\in\! B.$$ 
\bf Case 1.2\rm\ $\big( 2^\omega ,s(\mathbb{G}_0)\big)\sqsubseteq_c( 2^\omega ,B')$\bigskip

In this case, 
$(0\alpha ,1\beta )\!\in\! G_{s(\mathbb{G}_0)}\Leftrightarrow\big( 0h(\alpha ),1h(\beta )\big)\!\in\! B$, which implies that 
$$\big( 0h(\alpha ),1h(\beta )\big)\!\in\! B\Leftrightarrow\big( 0h(\beta ),1h(\alpha )\big)\!\in\! B.$$
In both cases, $B\cap\{\big( 0h(\alpha ),1h(\beta )\big)\mid\alpha ,\beta\!\in\! 2^\omega\}$ is 
$D_2(\boraone )\!\subseteq\! F_\sigma$.\bigskip

\noindent $\bullet$ We set 
$B''\! :=\!\{ (\alpha ,\beta )\!\in\! h[2^\omega ]\!\times\! h[2^\omega ]\mid (1\alpha ,0\beta )\!\in\! B\}$. Here again, $B''\!\subseteq\! s(\mathbb{G}_0)$.\bigskip

\noindent $\bullet$ Let us prove that we may assume that $B''$ is not pot$(\bormone )$. We argue by contradiction, which gives a non-meager $G_\delta$ subset $G'$ of $2^\omega$ such that 
$h[G']$ is $B''$-discrete. Note that 
$$B\cap (2\!\times\! h[G'])^2\! =\! G_{s(\mathbb{G}_0)}\cap (2\!\times\! h[G'])^2.$$
As $G'$ is not meager, there is no Borel countable coloring of $(G',\mathbb{G}_0\cap {G'}^2)$. Theorem \ref{KSTinjdigr} gives 
$f\! :\! 2^\omega\!\rightarrow\! G'$ injective continuous such that $\mathbb{G}_0\!\subseteq\! (f\!\times\! f)^{-1}(\mathbb{G}_0\cap {G'}^2)$. The map $h\circ f$ is a witness for the fact that there is no Borel countable coloring of $\big( h[G'],s(\mathbb{G}_0)\cap h[G']^2\big)$. Theorem \ref{G0} shows that 
$\big( 2^\omega ,s(\mathbb{G}_0)\big)\sqsubseteq_c\big( h[G'],s(\mathbb{G}_0)\cap h[G']^2\big)$ with witness $g'$. The map 
$\varepsilon\alpha\!\mapsto\!\varepsilon g'(\alpha )$ is a witness for 
$(2^\omega ,G_{s(\mathbb{G}_0)})\sqsubseteq_c\big( 2\!\times\! h[G'],G_{s(\mathbb{G}_0)}\cap (2\!\times\! h[G'])^2\big)$. Thus 
$(2^\omega ,G_{s(\mathbb{G}_0)})\sqsubseteq_c(2^\omega ,B)$.\bigskip

\noindent $\bullet$ By Theorem \ref{G0} again, 
$( 2^\omega ,\mathbb{G}_0)\sqsubseteq_c( h[2^\omega ],B'')$ or 
$\big( 2^\omega ,s(\mathbb{G}_0)\big)\sqsubseteq_c( h[2^\omega ],B'')$ with witness $h'$.\bigskip

\noindent\bf Case 2.1\rm\ $( 2^\omega ,\mathbb{G}_0)\sqsubseteq_c( h[2^\omega ],B'')$\bigskip

 In this case, $(1\alpha ,0\beta )\!\in\!\mathbb{T}_0\Leftrightarrow
\big( 1h'(\alpha ),0h'(\beta )\big)\!\in\! B$, and 
$$\big( 1h'(\alpha ),0h'(\beta )\big)\!\in\! B
\Rightarrow\big( 1h'(\beta ),0h'(\alpha )\big)\!\in\! s(B)\!\setminus\! B\! =\! B^{-1}
\Rightarrow\big( 0h'(\alpha ),1h'(\beta )\big)\!\in\! B.$$ 
\bf Case 2.2\rm\ $\big( 2^\omega ,s(\mathbb{G}_0)\big)\sqsubseteq_c( h[2^\omega ],B'')$\bigskip

 In this case, $(1\alpha ,0\beta )\!\in\! G^{-1}_{s(\mathbb{G}_0)}\Leftrightarrow
\big( 1h'(\alpha ),0h'(\beta )\big)\!\in\! B$, which implies that
$$\big( 1h'(\alpha ),0h'(\beta )\big)\!\in\! B\Leftrightarrow\big( 1h'(\beta ),0h'(\alpha )\big)\!\in\! B.$$
In both cases, $B\cap\{\big( 1h'(\alpha ),0h'(\beta )\big)\mid\alpha ,\beta\!\in\! 2^\omega\}$ is 
$D_2(\boraone )\!\subseteq\! F_\sigma$. As $h'[2^\omega ]\!\subseteq\! h[2^\omega ]$, the set 
$B\cap (2\!\times\! h'[2^\omega ])^2$ is $F_\sigma$.

\vfill\eject

\noindent $\bullet$ Now four new cases are possible.\bigskip
 
\noindent\bf Cases 1.1 and 2.1 hold\rm\bigskip

 As $h'[2^\omega ]\!\subseteq\! h[2^\omega ]$, 
$\big( 0h'(\alpha ),1h'(\beta )\big)\!\in\! B\Rightarrow\big( 1h'(\alpha ),0h'(\beta )\big)\!\in\! B$, so that $$\big( 0h'(\alpha ),1h'(\beta )\big)\!\in\! B\Leftrightarrow\big( 1h'(\alpha ),0h'(\beta )\big)\!\in\! B.$$ 
Moreover, this is equivalent to $(1\alpha ,0\beta )\!\in\!\mathbb{T}_0$, so that 
$(2^\omega ,\mathbb{T}_0)\sqsubseteq_c(X,B)$ with witness 
$\varepsilon\alpha\!\mapsto\!\varepsilon h'(\alpha )$.\bigskip
 
\noindent\bf Cases 1.1 and 2.2 hold\rm\bigskip

 Here again, 
$\big( 0h'(\alpha ),1h'(\beta )\big)\!\in\! B\Rightarrow\big( 1h'(\alpha ),0h'(\beta )\big)\!\in\! B$, and 
$(1\alpha ,0\beta )\!\in\! G^{-1}_{s(\mathbb{G}_0)}$ is equivalent to 
$\big( 1h'(\alpha ),0h'(\beta )\big)\!\in\! B\Leftrightarrow\big( 1h'(\beta ),0h'(\alpha )\big)\!\in\! B$. We set $f'(\varepsilon\alpha )\! :=\!\varepsilon h'(\alpha )$ and $B_0\! :=\! (f'\!\times\! f')^{-1}(B)$. Note that $(2^\omega ,B_0)\sqsubseteq_c(2^\omega ,B)$, 
$B_0\!\subseteq\! (N_0\!\times\! N_1)\cup (N_1\!\times\! N_0)$ is a $F_\sigma$ relation on 
$2^\omega$, contained in the $F_\sigma$ acyclic graph 
$(f'\!\times\! f')^{-1}\big( s(\mathbb{T}_0)\big)$, 
$B_0\cap (N_1\!\times\! N_0)\! =\! G^{-1}_{s(\mathbb{G}_0)}$, and 
$(0\alpha ,1\beta )\!\in\! B_0$ implies that $(1\alpha ,0\beta )\!\in\! B_0$.\bigskip

 We set $B'_0\! :=\!\{ (\alpha ,\beta )\!\in\! 2^\omega\!\times\! 2^\omega\mid (0\alpha ,1\beta )\!\in\! B_0\}$. Note that $B'_0$ is $F_\sigma$ and contained in the set $(h'\!\times\! h')^{-1}\big( s(\mathbb{G}_0)\big)$. By Theorem 1.1, there is a Borel countable coloring of $B'_0$, or 
$$(2^\omega ,\mathbb{G}_0)\sqsubseteq_c(2^\omega ,B'_0)\mbox{,}$$ 
or $\big( 2^\omega ,s(\mathbb{G}_0)\big)\sqsubseteq_c(2^\omega ,B'_0)$ with witness $h_0$.\bigskip

\noindent - In the first case, we get a non-meager $B'_0$-discrete $G_\delta$ subset $G_0$ of $2^\omega$. Note that 
$$B_0\cap (2\!\times\! G_0)^2\! =\! G^{-1}_{s(\mathbb{G}_0)}\cap (2\!\times\! G_0)^2.$$ 
As $(2^\omega ,G^{-1}_{s(\mathbb{G}_0)})\sqsubseteq_c
\big( 2\!\times\! G_0,G^{-1}_{s(\mathbb{G}_0)}\cap (2\!\times\! G_0)^2\big)$, 
$(2^\omega ,G_{s(\mathbb{G}_0)})\sqsubseteq_c(2^\omega ,B_0)\sqsubseteq_c(2^\omega ,B)$.\bigskip

\noindent - In the second case, we set $f_0(\varepsilon\alpha )\! :=\!\varepsilon h_0(\alpha )$ and 
$B_1\! :=\! (f_0\!\times\! f_0)^{-1}(B_0)$. Note that 
$$(2^\omega ,B_1)\sqsubseteq_c(2^\omega ,B_0)\mbox{,}$$ 
$B_1\!\subseteq\! (N_0\!\times\! N_1)\cup (N_1\!\times\! N_0)$ is a $F_\sigma$ relation on 
$2^\omega$, $B_1\cap (N_0\!\times\! N_1)\! =\! G_{\mathbb{G}_0}$, 
$(0\alpha ,1\beta )\!\in\! B_1$ implies that $(1\alpha ,0\beta )\!\in\! B_1$, and 
$(1\alpha ,0\beta )\!\in\! B_1\Leftrightarrow (1\beta ,0\alpha )\!\in\! B_1$.\bigskip

 Note that 
$$\mathbb{G}_0\! =\!\{ (\alpha ,\beta )\!\in\! 2^\omega\!\times\! 2^\omega\mid 
(0\alpha ,1\beta )\!\in\! B_1\}\!\subseteq\! B'_1\! :=\!\{ (\alpha ,\beta )\mid (1\alpha ,0\beta )\!\in\! B_1\}\!\subseteq\! (h_0\!\times\! h_0)^{-1}\big( s(\mathbb{G}_0)\big) .$$ 
Corollary \ref{corFsigma} gives $h_1\! :\! 2^\omega\!\rightarrow\! 2^\omega$ injective continuous such that 
$$\mathbb{G}_0\!\subseteq\! (h_1\!\times\! h_1)^{-1}(\mathbb{G}_0)\!\subseteq\! 
(h_1\!\times\! h_1)^{-1}(B'_1)\!\subseteq\! s(\mathbb{G}_0).$$ 
 By symmetry considerations, we see that 
$\mathbb{G}_0\! =\! (h_1\!\times\! h_1)^{-1}(\mathbb{G}_0)$ and 
$(h_1\!\times\! h_1)^{-1}(B'_1)\! =\! s(\mathbb{G}_0)$. This shows that the map 
$\varepsilon\alpha\!\mapsto\!\varepsilon h_1(\alpha )$ is a witness for 
$(2^\omega ,\mathbb{T}_0\cup G^{-1}_{s(\mathbb{G}_0)})\sqsubseteq_c(2^\omega ,B_1)$. Now the map $\varepsilon\alpha\!\mapsto\! (1\! -\!\varepsilon )\alpha$ is a witness for the fact that 
$(2^\omega ,\mathbb{U}_0)\sqsubseteq_c(2^\omega ,\mathbb{T}_0\cup 
G^{-1}_{s(\mathbb{G}_0)})$.\bigskip

\noindent - The third case is similar to and simpler than the second one. We get 
$\big( 2^\omega ,s(\mathbb{T}_0)\big)\sqsubseteq_c(2^\omega ,B_1)$.

\vfill\eject
 
\noindent\bf Cases 1.2 and 2.1 hold\rm\bigskip

 Here, 
$\big( 0h'(\alpha ),1h'(\beta )\big)\!\in\! B\Leftrightarrow\big( 0h'(\beta ),1h'(\alpha )\big)\!\in\! B$ and $(1\alpha ,0\beta )\!\in\!\mathbb{T}_0$ is equivalent to 
$\big( 1h'(\alpha ),0h'(\beta )\big)\!\in\! B$, which implies that 
$\big( 0h'(\alpha ),1h'(\beta )\big)\!\in\! B$. We set $f'(\varepsilon\alpha )\! :=\!\varepsilon h'(\alpha )$ and $B_0\! :=\! (f'\!\times\! f')^{-1}(B)$. Note that $(2^\omega ,B_0)\sqsubseteq_c(2^\omega ,B)$, 
$B_0\!\subseteq\! (N_0\!\times\! N_1)\cup (N_1\!\times\! N_0)$ is a $F_\sigma$ relation on 
$2^\omega$, contained in the $F_\sigma$ acyclic graph 
$(f'\!\times\! f')^{-1}\big( s(\mathbb{T}_0)\big)$, 
$B_0\cap (N_1\!\times\! N_0)\! =\!\mathbb{U}_0\cap (N_1\!\times\! N_0)$,  
$(0\alpha ,1\beta )\!\in\! B_0\Leftrightarrow (0\beta ,1\alpha )\!\in\! B_0$, and 
$(1\alpha ,0\beta )\!\in\! B_0\Rightarrow (0\alpha ,1\beta )\!\in\! B_0$. Note that 
$$\mathbb{G}_0\! =\!\{ (\alpha ,\beta )\!\in\! 2^\omega\!\times\! 2^\omega\mid 
(1\alpha ,0\beta )\!\in\! B_0\}\!\subseteq\! B'_0\! :=\!\{ (\alpha ,\beta )\mid (0\alpha ,1\beta )\!\in\! B_0\}\!\subseteq\! (h'\!\times\! h')^{-1}\big( s(\mathbb{G}_0)\big) .$$ 
Corollary \ref{corFsigma} gives $h_0\! :\! 2^\omega\!\rightarrow\! 2^\omega$ injective continuous such that 
$$\mathbb{G}_0\!\subseteq\! (h_0\!\times\! h_0)^{-1}(\mathbb{G}_0)\!\subseteq\! 
(h_0\!\times\! h_0)^{-1}(B'_0)\!\subseteq\! s(\mathbb{G}_0).$$ 
 By symmetry considerations, we see that 
$\mathbb{G}_0\! =\! (h_0\!\times\! h_0)^{-1}(\mathbb{G}_0)$ and 
$(h_0\!\times\! h_0)^{-1}(B'_0)\! =\! s(\mathbb{G}_0)$. This shows that the map 
$\varepsilon\alpha\!\mapsto\!\varepsilon h_0(\alpha )$ is a witness for 
$(2^\omega ,\mathbb{U}_0)\sqsubseteq_c(2^\omega ,B_0)$.\bigskip
 
\noindent\bf Cases 1.2 and 2.2 hold\rm\bigskip

 Here again, 
$\big( 0h'(\alpha ),1h'(\beta )\big)\!\in\! B\Leftrightarrow\big( 0h'(\beta ),1h'(\alpha )\big)\!\in\! B$ and 
$(1\alpha ,0\beta )\!\in\! G^{-1}_{s(\mathbb{G}_0)}$ is equivalent to 
$\big( 1h'(\alpha ),0h'(\beta )\big)\!\in\! B\Leftrightarrow\big( 1h'(\beta ),0h'(\alpha )\big)\!\in\! B$. We set $f'(\varepsilon\alpha )\! :=\!\varepsilon h'(\alpha )$ and $B_0\! :=\! (f'\!\times\! f')^{-1}(B)$. Note that $(2^\omega ,B_0)\sqsubseteq_c(2^\omega ,B)$, 
$B_0\!\subseteq\! (N_0\!\times\! N_1)\cup (N_1\!\times\! N_0)$ is a $F_\sigma$ relation on 
$2^\omega$, contained in the $F_\sigma$ acyclic graph 
$(f'\!\times\! f')^{-1}\big( s(\mathbb{T}_0)\big)$, 
$B_0\cap (N_1\!\times\! N_0)\! =\! G^{-1}_{s(\mathbb{G}_0)}$, and 
$(0\alpha ,1\beta )\!\in\! B_0\Leftrightarrow (0\beta ,1\alpha )\!\in\! B_0$. We set 
$B'_0\! :=\!\{ (\alpha ,\beta )\!\in\! 2^\omega\!\times\! 2^\omega\mid (0\alpha ,1\beta )\!\in\! B_0\}$. 
Note that $B'_0$ is a $F_\sigma$ graph on $2^\omega$ contained in the acyclic graph 
$(h'\!\times\! h')^{-1}\big( s(\mathbb{G}_0)\big)$. By Theorem \ref{G0}, either there is a Borel countable coloring of $B'_0$, or $\big( 2^\omega ,s(\mathbb{G}_0)\big)\sqsubseteq_c(2^\omega ,B'_0)$ with witness $h_0$.\bigskip

 In the first case, $(2^\omega ,G_{s(\mathbb{G}_0)})\sqsubseteq_c(2^\omega ,B_0)$, as when 1.1 and 2.2 hold. In the second case, we set $f_0(\varepsilon\alpha )\! :=\!\varepsilon h_0(\alpha )$ and $B_1\! :=\! (f_0\!\times\! f_0)^{-1}(B_0)$. Note that 
$(2^\omega ,B_1)\sqsubseteq_c(2^\omega ,B_0)$, 
$$B_1\!\subseteq\! (N_0\!\times\! N_1)\cup (N_1\!\times\! N_0)$$ 
is a $F_\sigma$ relation on $2^\omega$, $B_1\cap (N_0\!\times\! N_1)\! =\! G_{s(\mathbb{G}_0)}$, and $(1\alpha ,0\beta )\!\in\! B_1$ is equivalent to $(1\beta ,0\alpha )\!\in\! B_1$.\bigskip
 
 We set $S\! :=\!\big\{ (\alpha ,\beta )\!\in\! 2^\omega\!\times\! 2^\omega\mid 
\big( 0h_0(\alpha ),1h_0(\beta )\big)\!\in\! B_0\ \wedge\ 
\big( 1h_0(\alpha ),0h_0(\beta )\big)\!\in\! B_0\big\}$. Note that $S$ is a graph on $2^\omega$ contained in $s(\mathbb{G}_0)$. By Corollary \ref{corFsigma}, either there is a Borel countable coloring of 
$S$, or there is $g_0\! :\! 2^\omega\!\rightarrow\! 2^\omega$ injective continuous such that 
$\mathbb{G}_0\!\subseteq\! (g_0\!\times\! g_0)^{-1}(S)\!\subseteq\! 
(g_0\!\times\! g_0)^{-1}\big( s(\mathbb{G}_0)\big)\!\subseteq\! s(\mathbb{G}_0)$.\bigskip

\noindent - In the first subcase, we get a non-meager $S$-discrete $G_\delta$ subset $G_1$ of $2^\omega$. Note that $B_1\cap (2\!\times\! G_1)^2$ is a $F_\sigma$ Acyclic oriented graph on $2\!\times\! G_1$. Theorem \ref{G0} shows that 
$\big( 2^\omega ,s(\mathbb{G}_0)\big)\sqsubseteq_c\big( G_1,s(\mathbb{G}_0)\cap G_1^2\big)$ with witness $g_1$. The map $f_1\! :\!\varepsilon\alpha\!\mapsto\!\varepsilon g_1(\alpha )$ is a witness for 
$$(2^\omega ,G_{s(\mathbb{G}_0)})\sqsubseteq_c
\big( 2\!\times\! G_1,G_{s(\mathbb{G}_0)}\cap (2\!\times\! G_1)^2\big) .$$ 
We set $B_2\! :=\! (f_1\!\times\! f_1)^{-1}(B_1)$. Note that 
$(2^\omega ,B_2)\sqsubseteq_c(2^\omega ,B_1)$, 
$B_2\!\subseteq\! (N_0\!\times\! N_1)\cup (N_1\!\times\! N_0)$ is a $F_\sigma$ Acyclic oriented graph on $2^\omega$, $B_2\cap (N_0\!\times\! N_1)\! =\! G_{s(\mathbb{G}_0)}$, and 
$(1\alpha ,0\beta )\!\in\! B_2$ is equivalent to $(1\beta ,0\alpha )\!\in\! B_2$. By Theorem 
\ref{ogasGsG0}, 
$(2^\omega ,G_{s(\mathbb{G}_0)})\sqsubseteq_c(2^\omega ,B_2)\sqsubseteq_c(2^\omega ,B)$.

\vfill\eject

\noindent - In the second subcase, $(g_0\!\times\! g_0)^{-1}(S)\! =\! 
(g_0\!\times\! g_0)^{-1}\big( s(\mathbb{G}_0)\big)\! =\! s(\mathbb{G}_0)$ since $S$ is a graph. We set $f_2(\varepsilon\alpha )\! :=\!\varepsilon g_0(\alpha )$ and 
$B_3\! :=\! (f_2\!\times\! f_2)^{-1}(B_1)$. Note that 
$(2^\omega ,B_3)\sqsubseteq_c(2^\omega ,B_1)$, 
$$B_3\!\subseteq\! (N_0\!\times\! N_1)\cup (N_1\!\times\! N_0)$$ 
is a $F_\sigma$ relation on $2^\omega$, $B_3\cap (N_0\!\times\! N_1)\! =\! G_{s(\mathbb{G}_0)}$, and $(1\alpha ,0\beta )\!\in\! B_3$ is equivalent to $(1\beta ,0\alpha )\!\in\! B_3$. Moreover, 
$(0\alpha ,1\beta )\!\in\! B_3$ implies that $(1\alpha ,0\beta )\!\in\! B_3$.\bigskip

 We set 
$B'_3\! :=\!\{ (\alpha ,\beta )\!\in\! 2^\omega\!\times\! 2^\omega\mid (1\alpha ,0\beta )\!\in\! B_3\}$. We repeat the previous argument, which gives a relation $B_4$ on $2^\omega$ with 
$(2^\omega ,B_4)\sqsubseteq_c(2^\omega ,B_3)$, 
$B_4\cap (N_1\!\times\! N_0)\! =\! G^{-1}_{s(\mathbb{G}_0)}$, 
$(0\alpha ,1\beta )\!\in\! B_4$ is equivalent to $(0\beta ,1\alpha )\!\in\! B_4$, and 
$(1\alpha ,0\beta )\!\in\! B_4$ is equivalent to $(0\alpha ,1\beta )\!\in\! B_4$. This means that 
$B_4\! =\! s(\mathbb{T}_0)$.\hfill{$\square$}

\begin{thm} \label{minT0} The set $\mathbb{T}_0$ is $\sqsubseteq_c$-minimal among non-pot$(\bormone )$ relations.\end{thm}

\noindent\bf Proof.\rm ~Assume that $A\!\subseteq\! X^2$ is not pot$(\bormone )$ and 
$(X,A)\sqsubseteq_c(2^\omega ,\mathbb{T}_0)$ with witness $g$. Then $A$ is a 
$D_2(\boraone )$ Acyclic oriented graph with locally countable closure. By Corollary 
\ref{situation}, ${(2^\omega ,\mathbb{G}_0)\sqsubseteq_c(X,A)}$ or there is a $D_2(\boraone )$ Acyclic oriented graph $B\!\subseteq\! (N_0\!\times\! N_1)\cup (N_1\!\times\! N_0)$ on $2^\omega$ with locally countable closure such that $\mathbb{B}_0\! =\!\overline{\mathbb{B}_0}\cap B$ and 
$(2^\omega ,B)\sqsubseteq_c(X,A)$. In particular, $(2^\omega ,B)\sqsubseteq_c(2^\omega ,\mathbb{T}_0)$ with witness $h$. As $(0\alpha ,1\alpha )\!\in\!\overline{B}\!\setminus\! B$, 
$$\big( h(0\alpha ),h(1\alpha )\big)\!\in\!\overline{\mathbb{T}_0}\!\setminus\!\mathbb{T}_0\! =\!
\{ (\varepsilon\gamma ,(1\! -\!\varepsilon )\gamma )\mid\varepsilon\!\in\! 2\ \wedge\ 
\gamma\!\in\! 2^\omega\} .$$ 
If $(1\alpha ,0\alpha )\!\in\! B$, then $\big( h(1\alpha ),h(0\alpha )\big)\!\in\!\mathbb{T}_0\cap
\{ (\varepsilon\gamma ,(1\! -\!\varepsilon )\gamma )\mid\varepsilon\!\in\! 2\ \wedge\ 
\gamma\!\in\! 2^\omega\}$, which is absurd. This implies that 
$(\varepsilon\gamma ,(1\! -\!\varepsilon )\gamma )\!\notin\! B$ if $\varepsilon\!\in\! 2$ and 
$\gamma\!\in\! 2^\omega$. Thus $\mathbb{B}_0\! =\!\overline{\mathbb{B}_0}\cap s(B)$ and $s(B)$ is not 
pot$(\bormone )$. Note that $h$ is a witness for 
$\big( 2^\omega ,s(B)\big)\sqsubseteq_c\big( 2^\omega ,s(\mathbb{T}_0)\big)$. The minimality of 
$s(\mathbb{T}_0)$ implies that 
${\big( 2^\omega ,s(\mathbb{T}_0)\big)\sqsubseteq_c\big( 2^\omega ,s(B)\big)}$. Replacing $B$ with its pre-image if necessary, we may assume that $B$ is a $D_2(\boraone )$ oriented graph on 
$2^\omega$ such that $s(B)\! =\! s(\mathbb{T}_0)$. Theorem \ref{insidesT0} gives $A'$ in 
$\{ G_{s(\mathbb{G}_0)},\mathbb{T}_0,\mathbb{U}_0,s(\mathbb{T}_0)\}$ such that 
$(2^\omega ,A')\sqsubseteq_c(2^\omega ,B)$. Proposition \ref{biganti10} shows that 
$A'\! =\!\mathbb{T}_0$, and we are done.\hfill{$\square$}

\begin{thm} \label{minU0} The set $\mathbb{U}_0$ is $\sqsubseteq_c$-minimal among 
non-$\mbox{pot}(\bormone )$ sets.\end{thm}

\noindent\bf Proof.\rm ~Assume that $A\!\subseteq\! X^2$ is not pot$(\bormone )$ and 
$(X,A)\sqsubseteq_c(2^\omega ,\mathbb{U}_0)$ with witness $g$. Then $A$ is a $D_2(\boraone )$ Acyclic digraph with locally countable closure. By Corollary \ref{situation}, 
$(2^\omega ,\mathbb{G}_0)\sqsubseteq_c(X,A)$ or 
$\big( 2^\omega ,s(\mathbb{G}_0)\big)\sqsubseteq_c(X,A)$ or there is a $D_2(\boraone )$ Acyclic digraph $B$ on $2^\omega$ with locally countable closure contained in 
$(N_0\!\times\! N_1)\cup (N_1\!\times\! N_0)$ such that 
${\mathbb{B}_0\! =\!\overline{\mathbb{B}_0}\cap B}$ and $(2^\omega ,B)\sqsubseteq_c(X,A)$. In particular, $(2^\omega ,B)\sqsubseteq_c(2^\omega ,\mathbb{U}_0)$ with witness $h$. As in the proof of Theorem \ref{minT0}, we may assume that $B$ is a $D_2(\boraone )$ digraph on 
$2^\omega$ such that $s(B)\! =\! s(\mathbb{T}_0)$. We conclude as in the proof of Theorem 
\ref{minT0}.\hfill{$\square$}\bigskip

\noindent\bf Proof of Theorem \ref{addicho}.\rm ~We set 
${\cal A}''\! :=\!\{\mathbb{B}_0,\mathbb{N}_0,\mathbb{M}_0,G_{s(\mathbb{G}_0)},\mathbb{U}_0\}$,  
${\cal B}''\! :=\!\{\mathbb{T}_0,s(\mathbb{B}_0),s(\mathbb{T}_0)\}$. By Proposition \ref{biganti10}, 
${\cal A}''\cup {\cal B}''$ is a $\sqsubseteq_c$-antichain made of $D_2(\boraone )$ Acyclic relations, with  locally countable closure contained in 
$(N_0\!\times\! N_1)\cup (N_0\!\times\! N_1)$, which are not pot$(\bormone )$. This implies that 
${\cal A}'$ is made of $D_2(\boraone )$ Acyclic relations, with locally countable closure, which are not pot$(\bormone )$. By Lemma \ref{anti7}, 
${\cal A}'''\! :=\!\big\{ A^e\mid A\!\in\! {\cal A}''\ \wedge\ e\!\in\!\{ =,\square ,\sqsubset ,\sqsupset\}\big\}\cup\big\{ A^e\mid A\!\in\! {\cal B}''\ \wedge\ e\!\in\!\{ =,\square ,\sqsubset\}\big\}$ is also a 
$\sqsubseteq_c$-antichain. The proof of Proposition \ref{biganti10} shows that 
$\{\mathbb{G}_0,s(\mathbb{G}_0)\}\cup {\cal A}'''\! =\! {\cal A}'$ is a $\sqsubseteq_c$-antichain.

\vfill\eject

 By Theorems \ref{generalpos}.(5).(ii), \ref{minGsG0}, \ref{minT0} and \ref{minU0}, the elements of the antichain in the statement of Proposition \ref{biganti10} are $\sqsubseteq_c$-minimal (among 
non-pot$(\bormone )$ relations). By Proposition \ref{minimalitysquare}, $A^\square$ is 
$\sqsubseteq_c$-minimal if $A\!\in\! {\cal A}''\cup {\cal B}''$. By Theorem \ref{generalpos}.(5).(ii), the elements of $\cal A$ are also minimal. It remains to see that the elements of 
$\big\{ A^e\mid A\!\in\!\{ G_{s(\mathbb{G}_0)},\mathbb{U}_0\}\ \wedge\ 
e\!\in\!\{\sqsubset ,\sqsupset\}\big\}\cup
\big\{ A^\sqsubset\mid A\!\in\!\{\mathbb{T}_0,s(\mathbb{T}_0)\}\big\}$ are minimal. Let us do it for 
$A\! :=\!\mathbb{T}_0$, the other cases being similar. Assume that 
$(X,S)\sqsubseteq_c(2^\omega ,A^\sqsubset )$ with witness $f$, where $X$ is Polish and $S$ is not 
$\mbox{pot}(\bormone )$. Then $f$ is also a witness for 
$$\big( X,S\!\setminus\!\Delta (X)\big)\sqsubseteq_c(2^\omega ,A).$$ 
Note that $S$ is the disjoint union of $S\!\setminus\!\Delta (X)$ and 
$\Delta (J)\!\in\!\mbox{pot}(\bormone )$, where $J$ is a Borel subset of $X$. Thus $S\!\setminus\!\Delta (X)$ is not $\mbox{pot}(\bormone )$. By Theorem \ref{minT0}, $A$ is minimal among 
non-$\mbox{pot}(\bormone )$ relations. Thus 
$(2^\omega ,A)\sqsubseteq_c\big( X,S\!\setminus\!\Delta (X)\big)$ with witness $h$. We set 
$S'\! :=\! (h\!\times\! h)^{-1}(S)$, so that $(2^\omega ,S')\sqsubseteq_c(X,S)$, $S'\! =\! A\cup\Delta (I)$ (where $I$ is a Borel subset of $2^\omega$). This means that we may assume that $X\! =\! 2^\omega$ and 
$S\! =\! A\cup\Delta (I)$, where $I$ is a Borel subset of $2^\omega$. We set, for $\varepsilon\!\in\! 2$, 
$S_\varepsilon\! :=\!\{\alpha\!\in\! 2^\omega\mid\varepsilon\alpha\!\in\! I\}$. This defines a partition 
$\{ S_0\cap S_1,S_0\!\setminus\! S_1,S_1\!\setminus\! S_0,(\neg S_0)\cap (\neg S_1)\}$ of 
$2^\omega$ into Borel sets. By Baire's theorem, one of these sets is not meager. Let 
$s\!\in\! 2^{<\omega}$ and $C$ be a dense $G_\delta$ subset of $2^\omega$ such that 
$N_s\cap C$ is contained in one of these sets.\bigskip

 We saw in the proof of Lemma \ref{propertiespi01} that $\mathbb{G}_0\cap (N_s\cap C)^2$ is not pot$(\bormone )$ if $s\!\in\! 2^{<\omega}$ and $C$ is a dense $G_\delta$ subset of $2^\omega$. In particular, there is no Borel countable coloring of $\mathbb{G}_0\cap (N_s\cap C)^2$. By Theorem \ref{G0}, 
$(2^\omega ,\mathbb{G}_0)\sqsubseteq_c\big( N_s\cap C,\mathbb{G}_0\cap (N_s\cap C)^2\big)$ with witness $g$. This implies that the map 
$g'\! :\!\varepsilon\alpha\!\mapsto\!\varepsilon g(\alpha )$ reduces $A^\sqsubset$ to 
$A^\sqsubset\cap\big ( 2\!\times\! (N_s\cap C)\big)^2$.\bigskip 

\noindent\bf Case 1\rm\ $S_0\cap S_1$ is not meager.\bigskip

 The map $g'$ is a witness for $(2^\omega ,A^\square )\sqsubseteq_c
\Big( 2\!\times\! (N_s\cap C),A^\square\cap\big( 2\!\times\! (N_s\cap C)\big)^2\Big)$. Now note that 
$S\cap\big( 2\!\times\! (N_s\cap C)\big)^2\! =\! A^\square\cap\big( 2\!\times\! (N_s\cap C)\big)^2$, 
so that $(2^\omega ,A^\square )\sqsubseteq_c(2^\omega ,S)\sqsubseteq_c
(2^\omega ,A^\sqsubset )$, which contradicts the fact that ${\cal A}'$ is a $\sqsubseteq_c$-antichain.\bigskip

\noindent\bf Case 2\rm\ $S_0\!\setminus\! S_1$ is not meager.\bigskip

 The map $g'$ is a witness for $(2^\omega ,A^\sqsubset )\sqsubseteq_c
\Big( 2\!\times\! (N_s\cap C),A^\sqsubset\cap\big( 2\!\times\! (N_s\cap C)\big)^2\Big)$. Now note that 
$S\cap\big( 2\!\times\! (N_s\cap C)\big)^2\! =\! A^\sqsubset\cap\big( 2\!\times\! (N_s\cap C)\big)^2$, 
so that $(2^\omega ,A^\sqsubset )\sqsubseteq_c(2^\omega ,S)$.\bigskip

\noindent\bf Case 3\rm\ $S_1\!\setminus\! S_0$ is not meager.\bigskip

 The map $g'$ is a witness for $(2^\omega ,A^\sqsupset )\sqsubseteq_c
\Big( 2\!\times\! (N_s\cap C),A^\sqsupset\cap\big( 2\!\times\! (N_s\cap C)\big)^2\Big)$. Now note that 
$S\cap\big( 2\!\times\! (N_s\cap C)\big)^2\! =\! A^\sqsupset\cap\big( 2\!\times\! (N_s\cap C)\big)^2$, 
so that $(2^\omega ,A^\sqsupset )\sqsubseteq_c(2^\omega ,S)\sqsubseteq_c
(2^\omega ,A^\sqsubset )$. It remains to note that 
$(2^\omega ,A^\sqsubset )\sqsubseteq_c(2^\omega ,A^\sqsupset )$ with witness 
$\varepsilon\alpha\!\mapsto\! (1\! -\! \varepsilon )\alpha$ if 
$A\!\in\!\{\mathbb{T}_0,s(\mathbb{T}_0)\}$.\bigskip

\noindent\bf Case 4\rm\ $(\neg S_0)\cap (\neg S_1)$ is not meager.\bigskip

 The map $g'$ is a witness for $(2^\omega ,A^=)\sqsubseteq_c
\Big( 2\!\times\! (N_s\cap C),A^=\cap\big( 2\!\times\! (N_s\cap C)\big)^2\Big)$. Now note that 
$S\cap\big( 2\!\times\! (N_s\cap C)\big)^2\! =\! A^=\cap\big( 2\!\times\! (N_s\cap C)\big)^2$, 
so that $(2^\omega ,A^=)\sqsubseteq_c(2^\omega ,S)\sqsubseteq_c
(2^\omega ,A^\sqsubset )$, which contradicts the fact that ${\cal A}'$ is a $\sqsubseteq_c$-antichain.
\hfill{$\square$}

\vfill\eject

\section{$\!\!\!\!\!\!$ References}

\noindent [Ca]\ \ J. Calbrix, Classes de Baire et espaces d'applications continues,~\it C. R. Acad. Sci. Paris 
S\'er. I Math.\ \rm 301, 16 (1985) 759-762

\noindent [C-L-M]\ \ J. D. Clemens, D. Lecomte and B. D. Miller, Dichotomy theorems for families of non-cofinal essential complexity,\ \it to appear in Adv. Math. (see arXiv:1412.8684)\ \rm 

\noindent [Ka]\ \ M. Kat\v{e}tov, Products of filters,~\it Comment. Math. Univ. Carolinae\ \rm 9 (1968), 173-189

\noindent [K]\ \ A. S. Kechris,~\it Classical descriptive set theory,~\rm Springer-Verlag, 1995

\noindent [K-S-T]\ \ A. S. Kechris, S. Solecki and S. Todor\v cevi\'c, Borel chromatic numbers,\ \it 
Adv. Math.\rm\ 141 (1999), 1-44

\noindent [L1]\ \ D. Lecomte, Classes de Wadge potentielles et 
th\'eor\`emes d'uniformisation partielle,\it ~Fund. Math.~\rm 143 (1993), 231-258

\noindent [L2]\ \ D. Lecomte, Complexit\'e des bor\'eliens~\`a coupes 
d\'enombrables,\ \it Fund. Math.~\rm 165 (2000), 139-174

\noindent [L3]\ \ D. Lecomte, On minimal non potentially closed subsets of the plane,\ \it Topology Appl.\rm\ 154, 1 (2007), 241-262

\noindent [L4]\ \ D. Lecomte, How can we recognize potentially $\bormxi$ subsets of the plane?,~\it J. Math. Log.\rm\ 9, 1 (2009), 39-62

\noindent [L5]\ \ D. Lecomte, Potential Wadge classes,~\it\ Mem. Amer. Math. Soc.,\rm ~221, 1038 (2013)

\noindent [L-M]\ \ D. Lecomte and B. D. Miller, Basis theorems for non-potentially closed sets and graphs of uncountable Borel chromatic number,\ \it J. Math. Log.\ \rm 8 (2008), 1-42

\noindent [L-Z]\ \ D. Lecomte and R. Zamora, Injective tests of low complexity in the plane,\ \it preprint (see arXiv:1507.05015)\ \rm 

\noindent [Lo1]\ \ A. Louveau, A separation theorem for $\Ana$ sets,\ \it Trans. A. M. S.\rm\  260 (1980), 363-378 

\noindent [Lo2]\ \ A. Louveau, Ensembles analytiques et bor\'eliens dans les 
espaces produit,~\it Ast\'erisque (S. M. F.)\ \rm 78 (1980)

\noindent [P]\ \ T. C. Przymusinski, On the notion of n-cardinality,\ \it Proc. Amer. Soc.\ \rm 69 (1978), 333-338

\end{document}